\documentclass[final]{siamart171218}
\usepackage[T1]{fontenc}
\usepackage[latin1]{inputenc}
\usepackage{graphicx}
\usepackage{fullpage}
\usepackage{epic}

\usepackage{comment}
\usepackage{tikz}
\usepackage{wrapfig}

\usepackage{soul}

\usepackage[font=small]{caption}

\usepackage{pgfgantt}
\usepackage{fancyhdr}
\usepackage{wrapfig,color}
\usepackage{tikz}
\usetikzlibrary{shapes,arrows}
\usepackage{pgfplots}

\usepackage{mythesis}
\usepackage{amssymb}
\newtheorem{remark}{Remark}

\newcommand{\ba}{\begin{array}}
\newcommand{\ea}{\end{array}}
\newcommand{\bd}{\begin{displaymath}}
\newcommand{\ed}{\end{displaymath}}
\newcommand{\bi}{\begin{itemize}}
\newcommand{\ei}{\end{itemize}}
\newcommand{\bn}{\begin{enumerate}}
\newcommand{\en}{\end{enumerate}}
\newcommand{\pa}{\partial}

\newcommand{\f}{\frac}

\newtheorem{prob}{Problem}

\usepackage{mathtools}
\DeclarePairedDelimiter\ceil{\lceil}{\rceil}

\usepackage[normalem]{ulem}

\usepackage{calc}
\usepackage{tikz}
\newlength{\tfwidth}
\newlength{\tfheight}
\newlength{\tfxa}
\newlength{\tfxb}
\newlength{\tfya}
\newlength{\tfyb}
\newcommand{\trimFigWithBox}[6]{%
\setlength\fboxsep{0pt}%
\setlength\fboxrule{1.0pt}
\fbox{\includegraphics[width=#2, clip, trim=#3 #4 #5 #6]{#1}}%
}
\newcommand{\trimFigNoBox}[6]{%
\setlength\fboxsep{1pt}
\setlength\fboxrule{0.0pt}
\fbox{\includegraphics[width=#2, clip, trim=#3 #4 #5 #6]{#1}}%
}
\newcommand{\trimFigHeightWithBox}[6]{%
\setlength\fboxsep{0pt}%
\setlength\fboxrule{1.0pt}
\fbox{\includegraphics[height=#2, clip, trim=#3 #4 #5 #6]{#1}}%
}
\newcommand{\trimFigHeightNoBox}[6]{%
\setlength\fboxsep{1pt}
\setlength\fboxrule{0.0pt}
\fbox{\includegraphics[height=#2, clip, trim=#3 #4 #5 #6]{#1}}%
}

\newsavebox\figBox

\newcommand{\trimw}[6]{%
\sbox\figBox{\includegraphics{#1}}
\setlength{\tfwidth}{\the\wd\figBox}
\setlength{\tfheight}{\the\ht\figBox}
\setlength{\tfxa}{\tfwidth*\real{#3}}%
\setlength{\tfxb}{\tfwidth*\real{#4}}%
\setlength{\tfya}{\tfheight*\real{#5}}%
\setlength{\tfyb}{\tfheight*\real{#6}}%
\trimFigNoBox{#1}{#2}{\tfxa}{\tfya}{\tfxb}{\tfyb}%
}

\newcommand{\trimwb}[6]{%

\sbox\figBox{\includegraphics{#1}}
\setlength{\tfwidth}{\the\wd\figBox}
\setlength{\tfheight}{\the\ht\figBox}
\setlength{\tfxa}{\tfwidth*\real{#3}}%
\setlength{\tfxb}{\tfwidth*\real{#4}}%
\setlength{\tfya}{\tfheight*\real{#5}}%
\setlength{\tfyb}{\tfheight*\real{#6}}%
\trimFigWithBox{#1}{#2}{\tfxa}{\tfya}{\tfxb}{\tfyb}%
}

\newcommand{\trimh}[6]{%
\sbox\figBox{\includegraphics{#1}}
\setlength{\tfwidth}{\the\wd\figBox}
\setlength{\tfheight}{\the\ht\figBox}
\setlength{\tfxa}{\tfwidth*\real{#3}}%
\setlength{\tfxb}{\tfwidth*\real{#4}}%
\setlength{\tfya}{\tfheight*\real{#5}}%
\setlength{\tfyb}{\tfheight*\real{#6}}%
\trimFigHeightNoBox{#1}{#2}{\tfxa}{\tfya}{\tfxb}{\tfyb}%
}

\newcommand{\trimhb}[6]{%

\sbox\figBox{\includegraphics{#1}}
\setlength{\tfwidth}{\the\wd\figBox}
\setlength{\tfheight}{\the\ht\figBox}
\setlength{\tfxa}{\tfwidth*\real{#3}}%
\setlength{\tfxb}{\tfwidth*\real{#4}}%
\setlength{\tfya}{\tfheight*\real{#5}}%
\setlength{\tfyb}{\tfheight*\real{#6}}%
\trimFigHeightWithBox{#1}{#2}{\tfxa}{\tfya}{\tfxb}{\tfyb}%
}
\newcommand{\revb}[2]{{{#1}} {}}
\newcommand{\reva}[2]{{{#1}} {}}

\title{WaveHoltz: Iterative Solution of the Helmholtz Equation via the Wave Equation}
\author{Daniel Appel\"{o}\thanks{Department of Applied Mathematics, University of Colorado, Boulder, CO 80309-0526
(\email{Daniel.Appelo@colorado.edu})}\funding{Supported in part by STINT initiation grant IB2019--8154 and NSF Grant DMS-1913076. Any conclusions or recommendations expressed in this paper are those of the author and do not necessarily reflect the views of the NSF. }
\and
Fortino Garcia\thanks{Department of Applied Mathematics, University of Colorado, Boulder, CO 80309-0526
(\email{Fortino.Garcia@colorado.edu})}\funding{Supported in part by STINT initiation grant IB2019--8154 and NSF Grant DMS-1913076. Any conclusions or recommendations expressed in this paper are those of the author and do not necessarily reflect the views of the NSF}
\and 
Olof Runborg\thanks{
Department of Mathematics, KTH, 100 44, Stockholm (\email{olofr@kth.se})}
\funding{Supported in part by STINT initiation grant IB2019--8154.}
}

\date{\today}
\begin{document}
\maketitle
\begin{abstract}
\reva{A new iterative method, the WaveHoltz iteration, for  solution of the Helmholtz equation is presented. WaveHoltz is a fixed-point iteration that filters the solution of wave equation with time-periodic forcing and boundary data. The WaveHoltz iteration corresponds to a linear and coercive operator which, after discretization, can be recast as a positive definite linear system of equations. The solution to this system of equations approximates the Helmholtz solution and can be accelerated by Krylov subspace techniques. Analysis of the continuous and the discrete cases are presented as are numerical experiments.}{A new idea for iterative solution of the Helmholtz equation is presented. We show that the iteration which we denote WaveHoltz and which  filters the solution to the wave equation with harmonic data evolved over one period, corresponds to a coercive operator or a positive definite matrix in the discretized case.}   
\end{abstract}	




\section{Introduction}
The defining feature of waves are their ability to propagate over large distances without changing their shape. It is this property that allows them to carry information which underpins all communication, be it through speech or electronic transmission of data. Waves can also be used to probe the interior of the earth, the human body or engineering structures like buildings or bridges. This probing can be turned into images of the interior by the means of solving inverse problems. Harnessing the nature of waves requires high-order accurate and efficient numerical methods that are able to simulate wave propagation in three dimensions and over long distances. For cutting edge problems in scientific and engineering research such simulations must be carried out on parallel high-performance computing platforms and thus the numerical methods must scale while being easy to implement and generally applicable. 

In this paper we focus on approximating solutions to the scalar wave equation in the frequency domain, i.e. the Helmholtz equation
\begin{equation}
  \nabla\cdot(c^2(x)\nabla u) + \omega^2u = f(x). \label{eq:helm_bg}
\end{equation}
However, to to obtain such solutions we will use time domain discretizations of the wave equation. The motivation for developing high order accurate and scalable Helmholtz solvers comes from both mathematics and applications. On the mathematics side the recent results by Engquist and Zhao \cite{no_low_rank_Engq} give sharp lower bounds on the number of terms in a separated representation approximation of the Green's function of the Helmholtz equation as a function of the frequency (wavenumber). These bounds limit the applicability of the state of the art sweeping preconditioners in the high frequency regime and, for example, for interior and wave guide problems. Motivation also comes from applications in seismology, optics and acoustics. For example in full waveform inversion the problems are very large and the robustness of the inversion process can be enhanced by combining frequency and time domain inversion in a multi-scale fashion to avoid getting trapped in local minima.

Designing efficient iterative solvers for the Helmholtz equation
(\ref{eq:helm_bg})
is notoriously difficult and has been the subject of much research (for detailed reviews see  Ernst and Gander, \cite{ernst2012difficult}, Gander and Zhang \cite{Gander_Zhang_SIAM_REV}, and Erlangga, \cite{erlangga2008advances}). The main two difficulties in solving the Helmholtz equation are the resolution requirements and the highly indefinite character of the discretized system of equations. 

Assuming that (\ref{eq:helm_bg}) has been scaled so that the mean of $c(x)$ is about 1 then the typical wavelength is $\lambda = 2 \pi / \omega$ and the typical wavenumber is $\omega / 2 \pi$. In order to numerically propagate solutions to the time dependent wave equation corresponding to (\ref{eq:helm_bg}) with small errors  it is crucial to control the dispersion by using high order methods. The basic estimate by Kreiss and Oliger \cite{KreOli72} shows that in order to propagate a wave over $J$ wavelengths with a $p$th order finite difference method and with an error no greater than $\epsilon$ one must choose the number of points per wavelength ${\rm PPW}(J,p)$ as 
\[
{\rm PPW}(J,p) \ge C(p,\epsilon) J^{\frac{1}{p}}.
\] 
Here $C(p,\epsilon)$ depends on the tolerance $\epsilon$ but decreases with increasing order of accuracy $p$.  
Consequently, for a problem in $d$-dimensions and with fixed physical size the number of wavelengths in the domain will scale as $\omega^d$ and  to maintain a fixed tolerance the total number of degrees of freedom needed, \mbox{$N_p(\omega) = \mathcal{O}(\omega^{d(1+\frac{1}{p})})$}, is very large for high frequencies.

The dependence on $p$ and $\omega$ in $N_p(\omega)$ immediately reveals two fundamental criteria for designing high frequency Helmholtz solvers:

\begin{itemize}
\item[1.] The solvers must be {\bf parallel, memory lean} and they must {\bf scale well}. In 3D the number of degrees of freedom representing the solution cannot be stored on a single computer, and even on a parallel computer it is important to preserve the sparsity of the discrete version of (\ref{eq:helm_bg}).     
\item[2.] The underlying discretizations must be {\bf high order accurate}. At high frequencies and in 3D the extra penalty due to pollution / dispersion errors becomes prohibitive.            
\end{itemize}  

Further, the linear system matrix, $A$, resulting from direct discretization of (\ref{eq:helm_bg}) is indefinite so that the robust and easy to implement preconditioned conjugate gradient (PCG) method cannot be used. Instead the method of necessity becomes the preconditioned generalized minimal residual method (GMRES). To efficiently precondition GMRES one must exploit the intrinsic properties of the wave equation. The oscillatory nature of the Helmholtz Green's function and its discrete counterpart $A^{-1}$ can only be well approximated if the (unconditioned) Krylov subspace is allowed to grow quite large (with ``large'' scaling adversely with the frequency $\omega$ , \cite{ernst2012difficult}). The slow growth of the ``spanning power'' of the Krylov vectors is due to the underlying local connectivity of the discretization, preventing information to propagate rapidly. Efficient preconditioners must thus accelerate the propagation of information or reduce the cost of each iteration. Without preconditioners the iteration typically stagnates.   

Perhaps the first contribution that aimed to improve the propagation of information was the Analytic Incomplete LU preconditioner (AILU) by Gander and Nataf \cite{GANDER2000261}. The AILU preconditioner finds an $LDL^T$ factorization from an approximation of the same pseudodifferential operators that are used to construct non-reflecting boundary conditions \cite{EngMaj77,Appthesis,Hag99} and sweeps forward then backward  along one of the coordinate directions in a structured grid. 

The pioneering works on sweeping preconditioners by Engquist and Ying \cite{engquist2011sweepingH,engquist2011sweepingPML} were major breakthroughs in the solution of the Helmholtz equation. Similar to the AILU, the preconditioners in \cite{engquist2011sweepingH,engquist2011sweepingPML} use a $LDL^T$ decomposition but exploit the low rank properties of off-diagonal blocks together with perfectly matched layers to obtain solvers that converge in a small number of GMRES iterations. The papers  \cite{engquist2011sweepingH,engquist2011sweepingPML} were the two first instances of iterative Helmholtz solvers that converge in a small number of iterations that is almost independent of frequency. 

Once it had been established that low rank approximations, combined with clever use of sweeping and perfectly matched layers (PML), could be used to find Helmholtz solvers with linear scaling then many extensions and specializations were constructed. For example, in \cite{stolk2013rapidly} Stolk introduced a domain decomposition method with transmission conditions based on the perfectly matched layer (PML) that is able to achieve near linear scaling. Chen and Xiang, \cite{Helmholtz_Chen_Xiang_2013}, and Vion and Geuzaine, \cite{vion2014double}, also considered  sweeping domain decomposition method combined with PML and showed that their methods could be used as  efficient preconditioners for the Helmholtz equation. The method of polarized traces by Zepeda-N\'{u}\~{n}ez, Demanet and co-authors, \cite{leonardo_1,ZEPEDANUNEZ2016347,2018arXiv180108655Z}, is a two step sweeping preconditioner that  compresses the traces of the Greens function in an offline computation and utilizes incomplete Green's formulas to propagate the interface data. See also the recent review by Gander and Zhang \cite{Gander_Zhang_SIAM_REV} for connections between sweeping methods.   

Alongside iterative methods there are some attractive direct and multigrid methods. Examples from the class of direct methods are the Hierarchically Semi-Separable (HSS) parallel multifrontal sparse solver by deHoop and co-authors, \cite{2011_deHoop_Xia}, \revb{}{and} the spectral collocation solver by Gillman, Barnett and Martinsson, \cite{Gillman2015}\revb{, and the $p$-FEM approach of B\'eriot, Prinn and Gabard, \cite{BePriGab16}, which utilizes an \textit{a priori} error indicator to choose the polynomial order of each element}{}. Notable examples of multigrid methods are the Wave-ray method by Brandt and Livshits \cite{brandt1997wave,Livshits:2006aa} and the shifted Laplacian preconditioner with multigrid by Erlangga et al. \cite{Erlangga2006}.

As mentioned previously, the invention of sweeping preconditioners was a breakthrough and it is likely that they will have lasting and continuing impacts for the solution of the Helmholtz equation in various settings. There are, however, some limitations. First, in the recent paper \cite{no_low_rank_Engq}, Engquist and Zhao provide precise lower bounds on how the number of terms that are needed to approximate the Helmholtz Green's function depends on the frequency. In particular, for the high frequency regime they show that for interior problems and waveguides the rank of the off-diagonal elements grows fast, rendering sweeping preconditioners less efficient. They also show that the situation is, in general, worse in 3D than in 2D. This lack of compressibility may, in cases of practical importance, increase the cost of both the factorization and compression as well as the application of the compressed preconditioner. We note that this loss of compressibility at high frequency will also prevent direct methods such as  \cite{2011_deHoop_Xia,Gillman2015,BePriGab16} from reaching their most efficient regimes. An additional drawback of direct methods is their memory consumption for 3D problems.    

Another potential drawback with the sweeping methods is the long setup times before the solve. Of course all of the algorithms above do not suffer from this deficiency but many of them do. This may not be problematic when considering a background velocity that does not change but this is not the case, for example, when inverting for material parameters. In this case the velocity model will change constantly, necessitating a costly factorization in each update.

Finally, the two criterions 1.) and 2.) above are not so easy to meet for sweeping preconditioners. The sweep itself is intrinsically sequential and although there have been at least partially successful attempts to parallelize the sweeping methods it is hard to say that they are easy to parallelize in a scalable way. In a similar vein most of the methods use (and some rely on) low order discretizations.  Although it is possible to use higher order accurate discretizations together with sweeping preconditioners, their scarcity in the literature is noticeable.  

Another approach that is somewhat popular in the engineering literature is to simply run the wave equation for a long time to get a Helmholtz solution, see e.g. \cite{hile1998hybrid}. The theoretical underpinning of this approach is the {\it limiting amplitude principle} which says that every solution to the wave equation with an oscillatory forcing,
 in the exterior of a domain with reflecting boundary conditions tends to the Helmholtz solution. However, since the limiting amplitude principle only holds for exterior problems this approach does not work for interior problems and becomes very slow for problems with trapping waves. See e.g. the articles by 
Ladyzhenskaya \cite{Ladyzhenskaya:57}, Morawetz \cite{morawetz1962limiting} and
Vainberg \cite{Vainberg:75}.

An alternative approach, the so called Controllability Method (CM), was    originally proposed by Bristeau et al. \cite{bristeau1998controllability}. In the CM the solution to the Helmholtz equation is found by solving a convex constrained least-squares minimization problem where the deviation from time-periodicity is minimized in the classic wave equation energy. The basic ingredients in an iteration step in CM are: a.) the solution of a forward wave and a backward wave equation over one time-period, and b.) the solution of a symmetric coercive elliptic (and wave number independent) problem. 

In \cite{bristeau1998controllability} and the later spectral element implementations of CM by Heikkola et al.  \cite{HEIKKOLA20071553,HEIKKOLA2007344} only sound-soft scatterers were considered. For more general boundary conditions the minimizer of  the cost functional of \cite{bristeau1998controllability} is not unique but alternative cost functionals that does guarantee uniqueness (and thus convergence to the Helmholtz solution) were recently proposed by Grote and Tang in \cite{grote2019controllability}. We also note that if the wave equation is formulated as a first order system it is possible to avoid solving the elliptic problem \cite{GlowRoss06,2019arXiv190312522G}.

In what follows we will present an alternative to the controllability method.  Our method, which we call the WaveHoltz Iteration method (WHI), only requires a single forward wave equation solve and no elliptic solves but produces a positive definite (and sometimes symmetric) iteration that can be accelerated by, e.g. the conjugate gradient method or other Krylov subspace methods. As the WaveHoltz iteration is built from a time domain wave equation solver we claim and hope to demonstrate that it meets both criterion 1. and 2. above.   

The rest of the paper is organized as follows. In Section 2 we present and analyze our method and its extensions, in Section 3 we briefly outline the numerical methods we use to solve the wave equation, in Section 4 we present numerical experiments, and in Section 5 we summarize and conclude.  
 
Before proceeding we would like to acknowledge that although our method is distinct from the controllability method, it was the work by Grote and Tang, \cite{grote2019controllability}, that introduced us to CM and inspired us to derive the method discussed below.


\section{WaveHoltz: A New Method for Designing Scalable Parallel Helmholtz Solvers}
We consider the Helmholtz equation in a bounded open smooth
domain $\Omega$,
\begin{equation}
  \nabla\cdot(c^2(x)\nabla u) + \omega^2u = f(x),\qquad
  x\in \Omega, \label{eq:helm}
\end{equation}
with boundary condtions of the type
\begin{equation}
i\alpha\omega u + \beta(c^2(x)\vec{n}\cdot \nabla u)=0,
\quad \alpha^2+\beta^2=1,\quad
\qquad x\in \partial\Omega.
 \label{eq:helmbc}
\end{equation}
We assume $f\in L^2(\Omega)$ and that $c\in L^\infty(\Omega)$ with
the bounds $0< c_{\min}\leq c(x)\leq c_{\max}<\infty$ a.e. in $\Omega$.
Away from resonances, this ensures that there is 
a unique weak solution $u\in H^1(\Omega)$
to (\ref{eq:helm}).
Due to the boundary conditions $u$ is in general complex valued.

We first note that the function $w(t,x) := u(x) \exp(i\omega t)$
is a $T=2\pi/\omega$-periodic (in time) 
solution to the forced scalar wave equation
\begin{eqnarray}
&&  w_{tt} =   \nabla\cdot(c^2(x)\nabla w)
  - f(x) e^{i\omega t}, \quad x\in\Omega, \ \ 0 \le t \le T, \nonumber \\
&&  w(0,x) = v_0(x), \quad w_t(0,x) = v_1(x), \nonumber \\
&&  \alpha w_t +\beta(c^2(x)\vec{n}\cdot \nabla w)=0, \quad x\in\partial\Omega
, \label{eq:apa}
\end{eqnarray}
where $v_0=u$ and $v_1=i\omega u$.
Based on this observation, our approach is to find this $w$
instead of $u$. We could thus look for
initial data  $v_0$ and $v_1$ such that $w$ is a $T$-periodic solution
to (\ref{eq:apa}). 
However, there may be several such $w$, see \cite{grote2019controllability}, 
and we therefore impose
the alternative constraint that a certain time-average of $w$
should equal the initial data. 
More precisely, we introduce the
following operator acting on the initial data $v_0\in H^1(\Omega)$, $v_1\in L^2(\Omega)$,
\be
{\Pi} \left[
\begin{array}{c}
v_0\\
v_1 
\end{array}
\right] = \frac{2}{T}\int_0^{T}\left(\cos(\omega t)-\frac14\right) \left[
\begin{array}{c}
w(t,x)\\
w_t(t,x) 
\end{array}
\right] dt,\quad
T=\frac{2\pi}{\omega},
\ee
where $w(t,x)$ and its time derivative $w_t(t,x)$ satisfies the wave equation (\ref{eq:apa}) with initial data $v_0$ and $v_1$.
The result of $\Pi [v_0,\ v_1]^T$ can thus be seen as
a filtering in time 
of $w(\cdot,x)$
around the $\omega$-frequency.
We will further 
motivate the choice of time averaging in the analysis below.
By construction, the solution $u$ of Helmholtz now satisfies the equation
\begin{equation}\lbeq{iteration_eq}
\left[
\begin{array}{c}
u\\
i\omega u 
\end{array}
\right] = {\Pi} \left[
\begin{array}{c}
u\\
i\omega u 
\end{array}
\right].
\end{equation}
The WaveHoltz method then amounts to solving
this equation with the fixed point iteration
\begin{equation}
\boxed{
\left[
\begin{array}{c}
v\\
v' 
\end{array}
\right]^{(n+1)} = {\Pi} \left[
\begin{array}{c}
v\\
v' 
\end{array}
\right]^{(n)}, \qquad \left[
\begin{array}{c}
v\\
v' 
\end{array}
\right]^{(0)} \equiv 0.
}
\end{equation}
Provided this iteration converges and the solution to \eq{iteration_eq}
is unique, we 
obtain the Helmholtz solution as
$u=\lim_{n\to\infty} v^n$.

\begin{remark}
Note that each iteration is inexpensive and that $T$ is reduced by the reciprocal of $\omega$ as $\omega$ grows. If we assume that the number of degrees of freedom in each dimension scales with $\omega$ and that we evolve the wave equation with an explicit method this means that the number of timesteps per iteration is independent of $\omega$. Also note that the iteration is trivial to implement (in parallel or serial) if there is already a time domain wave equation solver in place. The integral in the filtering is carried out independently for each degree of freedom and simply amounts to adding up a weighted  sum (e.g. a trapezoidal sum) of the solution one timestep at a time. Finally, note that WHI allows all the advanced techniques that have been developed for wave equations (e.g. local timestepping, non-conforming discontinuous Galerkin finite elements $h$- and $p$-adaptivity etc.) can be transferred to the Helmholtz equation and other time harmonic problems.       
\end{remark}

\subsection{Iteration for the Energy Conserving Case} \label{sec:iteration}

Here we consider boundary conditions of either
Dirichlet ($\beta=0$) or Neumann ($\alpha=0$) type. 
This is typically
the most difficult case for iterative Helmholtz solvers
when $\Omega$ is bounded.
The
wave energy is preserved in time and certain $\omega$-frequencies
in Helmholtz are resonant, meaning they equal an eigenvalue
of the operator $-\nabla\cdot(c^2(x)\nabla)$. Moreover,
the limiting amplitude principle does not hold, and one
can thus not obtain the Helmholtz solution by solving the
wave equation over a long time interval.

We start by introducing a simplified iteration for this case.
With the given boundary conditions
the solution to Helmholtz will
be real valued, since
$f$ is a real valued function.
Without loss of generality, we may then take $w_t(0,x)=0$ and $w(t,x) = u(x) \cos(\omega t)$, since for a $T$-periodic real valued 
solution there is a time when $w_t(0,x)=0$. We choose that time as the initial time so that (\ref{eq:apa}) becomes
\begin{eqnarray}
&&  w_{tt} =   \nabla\cdot(c(x)^2\nabla w)
  - f(x) \cos(\omega t), \quad x\in\Omega, \ \ 0 \le t \le T, \nonumber \\
&&  w(0,x) = v(x), \quad w_t(0,x) \equiv 0, \nonumber \\
&&  \alpha w_t +\beta(c^2(x)\vec{n}\cdot \nabla w)=0, \quad x\in\partial\Omega.
    \label{eq:wave}
\end{eqnarray}
The simplified iteration is then defined as
\begin{equation}\lbeq{iteration}
\boxed{  v^{n+1} = \Pi v^{n}, \qquad v^{0}\equiv 0,} 
\end{equation}
where
\be\lbeq{filterstep}
  \Pi v = \frac{2}{T}\int_0^{T}\left(\cos(\omega t)-\frac14\right)w(t,x) dt,\qquad
  T=\frac{2\pi}{\omega},
\ee
with $w(t,x)$ solving the wave equation (\ref{eq:wave}) with initial data $v=v^{n}\in H^1(\Omega)$.
We now analyze this iteration.

%
%
%
%
By the choice of boundary conditions the operator 
$-\nabla\cdot(c^2(x)\nabla)$ has a point spectrum with non-negative
eigenvalues \reva{with corresponding eigenfunctions that form an orthonormal basis of $L^2(\Omega)$}{}. Denote those eigenmodes $(\lambda_j^2,\phi_j(x))$, \reva{with $\|\phi_j\|_{L^2(\Omega)} = 1$}{}. We assume that the angular frequency $\omega$ is not a resonance, i.e. $\omega^2\neq \lambda_j^2$ for all $j$. The Helmholtz equation (\ref{eq:helm}) is then wellposed. 

\begin{figure}[]
\graphicspath{{figures/}}
\begin{center}
\includegraphics[width=0.48\textwidth]{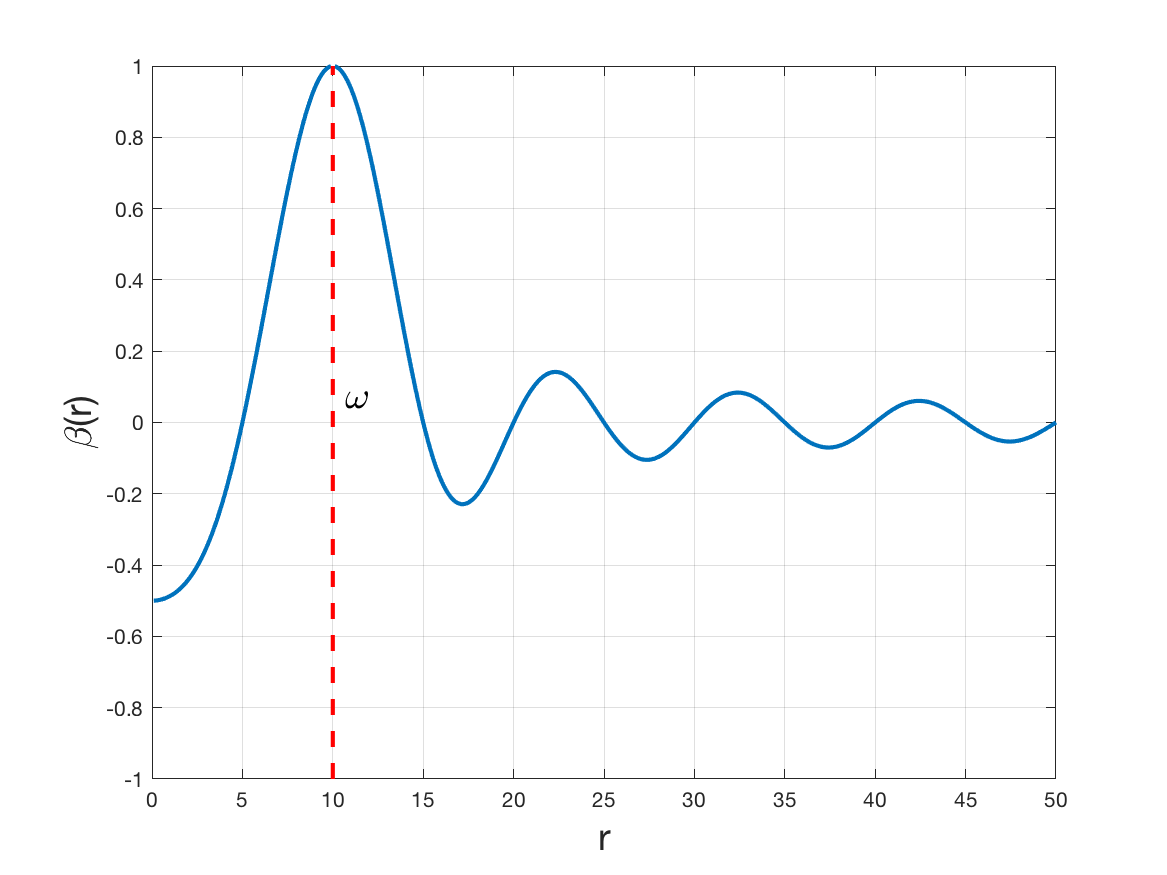}
\caption{The filter transfer function $\beta$ for $\omega=10$.
\label{fig:beta}}
\end{center}
\end{figure}

We recall that for any $q\in L^2(\Omega)$ we can
expand
$$
  q(x) = \sum_{j=0}^\infty \hat{q}_j\phi_j(x),
$$
for some coefficients $\hat{q}_j$ and
$$
   ||q||_{L^2(\Omega)}^2 = \sum_{j=0}^\infty |\hat{q}_j|^2,\qquad
   \reva{c_{\min}^2||\nabla q||_{L^2(\Omega)}^2}{c_{\min}^4||\nabla q||_{L^2(\Omega)}^2}
\leq   \sum_{j=0}^\infty \lambda_j^2|\hat{q}_j|^2\leq 
      \reva{ c_{\max}^2||\nabla q||_{L^2(\Omega)}^2}{c_{\max}^4||\nabla q||_{L^2(\Omega)}^2}.
$$
We start by
expanding the Helmholtz solution $u$, the initial data $v$ to the wave equation (\ref{eq:wave}), and the forcing $f$ in this way,
\[
  u(x) = \sum_{j=0}^\infty \hat{u}_j \phi_j(x), \ \ \ \  v(x) = \sum_{j=0}^\infty \hat{v}_j \phi_j(x), \ \ \ \    f(x) = \sum_{j=0}^\infty \hat{f}_j \phi_j(x).
\]
Then,
\[
-\lambda_j^2\hat{u}_j+\omega^2\hat{u}_j = \hat{f}_j \quad\Rightarrow\quad
\hat{u}_j = \frac{\hat{f}_j}{\omega^2-\lambda_j^2}.
\]
For the wave equation solution $w(t,x)$ with initial data $w=v$ and $w_t=0$ we have
\be
  w(t,x) = \sum_{j=0}^\infty \hat{w}_j(t) \phi_j(x),\qquad
  \hat{w}_j(t) 
  =\hat{u}_j\Bigl(\cos(\omega t) -\cos(\lambda_j t)\Bigr)
  + \hat{v}_j\cos(\lambda_j t). \label{eq:no_filter_yet}
\ee
The filtering step \eq{filterstep} then gives
$$
\Pi v = \sum_{j=0}^\infty \bar{v}_j \phi_j(x),\qquad
\bar{v}_j = \hat{u}_j(1-\beta(\lambda_j)) + \hat{v}_j\beta(\lambda_j), 
$$
where
\[
\beta(\lambda)
:= \frac{2}{T}\int_0^T \left(\cos(\omega t)-\frac14\right)\cos(\lambda t) dt.
 \]
 \reva{
We introduce the linear operator ${\mathcal S}:L^2(\Omega)\to L^2(\Omega)$,}{
Introducing the linear operator ${\mathcal S}:L^2(\Omega)\to L^2(\Omega)$,}
\be\lbeq{Sdef}
 {\mathcal S} \sum_{j=0}^\infty \hat{u}_j \phi_j(x) := \sum_{j=0}^\infty \beta(\lambda_j)\hat{u}_j \phi_j(x),
\ee
 \reva{
which gives the filtered solution of the wave equation with $f=0$, when applied
to the initial data $v$.
We can then write the iteration as
}
{
we can write the iteration as
}
\begin{equation}
  v^{n+1} =\Pi v^n= {\mathcal S}(v^n-u)+u. \label{eq:Siter}
\end{equation}
The operator ${\mathcal S}$ is self-adjoint and has the same eigenfunctions $\phi_j(x)$ as $-\nabla\cdot(c^2(x)\nabla)$ but with the (real)
eigenvalues  $\beta(\lambda_j)$. 
The convergence properties of the iteration depend on these
eigenvalues and
it is therefore of interest to study the 
range of the filter transfer function $\beta$.
Figure~\ref{fig:beta} shows a plot of $\beta$ which indicates
that the eigenvalues of $ {\mathcal S}$ 
are inside the unit interval, with 
a few of them being close to 1 (when $\lambda_j\approx \omega$), and
most of them being close to zero (when $\lambda_j\gg \omega$). 
In the appendix we show the following lemma about $\beta$.
\begin{lemma}\lblem{filterlemma}
The filter transfer function $\beta$ satisfies $\beta(\omega)=1$ and
\begin{align*}\lbeq{betabound}
0 \leq    \beta(\lambda) &\leq
   1 - \frac12 \left(\frac{\lambda-\omega}{\omega}\right)^2, 
   &{\rm when}\  \left|\frac{\lambda-\omega}{\omega}\right|\leq \frac12,\\
   |\beta(\lambda)| &\leq \frac12,
   &{\rm when}\  \left|\frac{\lambda-\omega}{\omega}\right|\geq \frac12,\\
   |\beta(\lambda)| & \leq 
   b_0 \frac{\omega}{\lambda-\omega}, &{\rm when}\  \lambda>\omega.
\end{align*}
where  $b_0 = \frac{3}{4\pi}$.
Moreover, 
close to $\omega$ we have the local expansion
\be\lbeq{locexp}
\beta(\omega+r) = 1- b_1 \left(\frac{r}{\omega}\right)^2 +
R(r/\omega)\left(\frac{r}{\omega}\right)^3 ,\qquad b_1 = \frac{2\pi^2}{3}-\frac14\approx 6.33,\quad
||R||_\infty\leq \reva{\frac{5\,\pi^{3}}{6}}{\frac{10\,\pi^{3}}{9}}.
\ee
\end{lemma}
\begin{remark}
It is easy to see that $\beta(\omega)=1$ for any
constant besides $1/4$. The particular choice $1/4$ is made to ensure
that $\beta'(\omega)=0$, which is necessary to keep $\beta\leq 1$ in a neighborhood of $\omega$.
We explore other possibilities in Section~\ref{sec:time_dep_fil}.
\end{remark}

From this lemma we can derive some results for the 
operator ${\mathcal S}$. To do this we first quantify
the non-resonance condition. We let
\[
 \delta_j =\frac{\lambda_j-\omega}{\omega},
\]
be the relative size of the gap between $\lambda_j$ and the Helmholtz frequency, and then denote the smallest gap (in magnitude) by $\delta$,
$$
   \delta =\delta_{j^*},\qquad 
   j^*={\rm argmin}_j |\delta_j|.
$$ 
Then we have
\begin{lemma}\lblem{Slemma}
Suppose $\delta>0$. The spectral radius $\rho$ of ${\mathcal S}$
is strictly less than one, and for small $\delta$,
\be\lbeq{rhoest}
  \rho = 1-b_1\delta^2 + O(\delta^3),
\ee
with $b_1$ as in \lem{filterlemma}. Moreover,
${\mathcal S}$ is a bounded linear map from
$L^2(\Omega)$ to $H^1(\Omega)$.
\end{lemma}
\begin{proof}
From \lem{filterlemma} we get
$$
   \rho = \sup_j |\beta(\lambda_j)|
   \leq \sup_j\max\left(1-\frac12\delta_j^2,\ \frac12\right)
   \leq \max\left(1-\frac12\delta^2,\ \frac12\right)<1.
$$ 
For the more precise estimate when $\delta$ is small
we will use \eq{locexp}. 
Since $1>\rho\geq\beta(\omega+\omega\delta)\to 1$ as $\delta\to 0$,
we can assume that
$\rho> 1-\eta^2/2$, with $\eta:=b_1/2||R||_\infty$,
for small enough $\delta$.
Then, since $|\beta(\omega+\omega\delta_j)|\leq 1-\eta^2/2$
for $|\delta_j|>\eta$ 
by \lem{filterlemma}, we have
$$
   \rho    =
   \sup_{|\delta_j|\leq \eta} \beta(\omega+\omega\delta_j) = 
   \beta(\omega+\omega\delta_{k^*}),
$$ 
for some $k^*$ with $|\delta_{k^*}|\leq \eta$. 
If $\delta_{k^*}=\delta_{j^*}$ (where $\delta=|\delta_{j^*}|$)
then \eq{locexp} gives \eq{rhoest}.
If not, we have
$\eta\geq |\delta_{k^*}|\geq\delta$ and by \lem{filterlemma}
$$
  0\leq \beta(\omega+\omega\delta_{k^*})-
  \beta(\omega+\omega\delta_{j^*})
  =-b_1(\delta^2_{k^*}-\delta^2)
  +R(\delta_{k^*})\delta^3_{k^*}
  -R(\delta_{j^*})\delta^3_{j^*}
  \leq
  -b_1(\delta^2_{k^*}-\delta^2)+
   \frac{b_1}2(\delta^2_{k^*}+\delta^2),
$$
which implies that $\delta_{k^*}^2\leq 3\delta^2$ and
\reva{that
$$
0\leq b_1(\delta^2_{k^*}-\delta^2)\leq 
  R(\delta_{k^*})\delta^3_{k^*}
  -R(\delta_{j^*})\delta^3_{j^*}\leq ||R||_\infty (1+3\sqrt{3})\delta^3.
$$
Therefore,}{therefore}
$$
{\rho} = 
1-b_1\delta_{k^*}^2 + O(\delta_{k^*}^3)
=
1-b_1\delta^2 + 
b_1(\delta^2-\delta_{k^*}^2)+
O(\delta_{k^*}^3)
=
1-b_1\delta^2 + 
O(\delta_{k^*}^3+\delta^3)
=
1-b_1\delta^2 + 
O(\delta^3).
$$
This shows \eq{rhoest}.
For the second statement, we note first that by \eq{betabound},
$$
  |\lambda_j\beta(\lambda_j)|\leq 
  \omega\ \begin{cases}
   1, & \lambda_j\leq \omega,\\
   \frac{b_0\lambda_j}{\lambda_j-\omega}, & \lambda_j> \omega,
  \end{cases}
  \ \ 
=  \omega\ \begin{cases}
   1, & \lambda_j\leq \omega,\\
   b_0(1+1/\delta_j), & \lambda_j> \omega,
  \end{cases}\ \ 
  \leq \omega \min(1,b_0(1+1/|\delta|))=:D.
$$
Suppose now that $g\in L^2(\Omega)$ and
$$
  g(x) = \sum_{j=0}^\infty \hat{g}_j \phi_j(x).
$$
Then
\begin{align*}
   ||{\mathcal S}g||_{H^1(\Omega)}^2
   &\leq \sum_{j=0}^\infty |\beta(\lambda_j)|^2|\hat{g}_j|^2 +
   \sum_{j=0}^\infty\frac{\lambda_j^2}{\reva{c^2_{\min}}{c^4_{\min}}}|\beta(\lambda_j)|^2|\hat{g}_j|^2
   \leq 
   \left(1 +
   \frac{D^2}{\reva{c^2_{\min}}{c^4_{\min}}}\right)
   \sum_{j=0}^\infty |\hat{g}_j|^2=
   \left(1 +
   \frac{D^2}{\reva{c^2_{\min}}{c^4_{\min}}}\right)
   ||g||_{L^2(\Omega)}^2.
\end{align*}
This proves the lemma.
\end{proof}
Letting $e^n:=u-v^n$ we can rearrange
(\ref{eq:Siter}) and obtain
$$
  e^{n+1}=  {\mathcal S}e^n\quad\Rightarrow\quad
  ||e^{n+1}||_{L^2(\Omega)}\leq \rho
  ||e^{n}||_{L^2(\Omega)}
  \quad\Rightarrow\quad ||e^{n}||_{L^2(\Omega)}
  \leq \rho^n ||e^{0}||_{L^2(\Omega)}\to 0,
$$
which shows that $v^n$ converges to $u$ in $L^2$. 
By \lem{Slemma} all iterates $v^n\in H^1(\Omega)$ since
$v^0=0$. We can therefore also get convergence in $H^1$.
Let 
$$
e^n(x) = 
\sum_{j=0}^\infty \hat{e}^{n}_j\phi_j(x),
$$
and consider similarly
\begin{align*}
\sum_{j=0}^\infty |\hat{e}^{n+1}_j|^2\lambda_j^2
&=
\sum_{j=0}^\infty \beta(\lambda_j)^2|\hat{e}^{n}_j|^2\lambda_j^2
\leq \rho^2\sum_{j=0}^\infty |\hat{e}^{n}_j|^2\lambda_j^2
\quad\Rightarrow\quad\\
||\nabla e^n||^2_{L^2(\Omega)}&\leq 
\frac{1}{\reva{c^2_{\min}}{c^4_{\min}}}
\sum_{j=0}^\infty |\hat{e}^{n}_j|^2\lambda_j^2
\leq
\frac{\rho^{2n}}{\reva{c^2_{\min}}{c^4_{\min}}}
\sum_{j=0}^\infty |\hat{e}^{0}_j|^2\lambda_j^2
\leq
\rho^{2n}\frac{\reva{c^2_{\max}}{c^4_{\max}}}{\reva{c^2_{\min}}{c^4_{\min}}}
||\nabla e^0||^2_{L^2(\Omega)}\to 0.
\end{align*}
We conclude that the iteration converges 
in $H^1$
with convergence rate $\rho$.
By \lem{filterlemma} we have
$\rho\sim 1-6.33\delta^2$ and, 
not surprisingly, the smallest gap, $\delta$, determines the convergence factor. 
We have thus showed
\begin{theorem}
The iteration in \eq{iteration} and \eq{filterstep} converges
in $H^1(\Omega)$
for the Dirichlet and Neumann problems
away from resonances
to the solution of the Helmholtz equation 
(\ref{eq:helm}). The convergence rate is $1-O(\delta^2)$,
where $\delta$ is the minimum gap between $\omega$ and 
the eigenvalues of $-\nabla\cdot(c^2(x)\nabla)$.
\end{theorem}
As discussed in the introduction,
the dependence of the convergence rate
on $\omega$ is often of interest. For the energy conserving
case, however,
this question is ambiguous as the problem
is not well-defined for all $\omega$. As soon as $\omega=\lambda_j$
there are either no or an infinite number of solutions.
In higher dimensions,
the eigenvalues $\lambda_j$ get denser as $j$ increases, meaning
that in general the problem will be closer and closer to resonance
as $\omega$ grows. Therefore, solving
the interior undamped Helmholtz equation for
high frequencies,
with pure Dirichlet or Neumann 
boundary conditions, may not be of great practical interest.

Nevertheless, we can make the following analysis. By the work of Weyl \cite{weyl:11} we know  that the eigenvalues grow asymptotically as $\lambda_j \sim j^{1/d}$ in $d$ dimensions. The average minimum gap $\delta$ when $\omega\approx \lambda_j$ is then 
$$
\delta \approx
\frac{1}{\lambda_{j+1}-\lambda_{j}}\int_{\lambda_j}^{\lambda_{j+1}}
\frac{\min(\lambda-\lambda_j,\lambda_{j+1}-\lambda)}{\omega}d\lambda
 =\frac{\lambda_{j+1}-\lambda_j}{4\omega}\sim \frac{(j+1)^{1/d}-j^{1/d}}{\omega}
\approx \frac{j^{1/d-1}}{d\omega}\sim
\frac{\omega^{1-d}}{\omega}\sim \omega^{-d}.
$$
When the convergence rate is $1-O(\delta^2)$, the number
iterations to achieve a fixed accuracy grows as $O(1/\delta^2)$.
This shows that the number of iterations 
would grow at the unacceptable rate $\omega^{2d}$
for the iteration.

Fortunately, one can accelerate the convergence by using the conjugate gradient method in the energy conserving case and with any other Krylov method in the general case. The linear system that we actually want to solve is 
$$
(I- {\mathcal S})v =:  {\mathcal A}v = b := \Pi 0.
$$
Moreover, with $b=\Pi 0$ pre-computed we can easily evaluate the action of ${\mathcal A}$ 
at the cost of a single wave solve. Precisely, since $ {\mathcal A}v=v-\Pi v+b$ we simply carry out the evaluation of ${\mathcal A}v$
by evolving the wave equation for one period in time with $v$ as the initial data and then subtract the filtered solution from the sum of the initial data and the right hand side $b$.


The operator ${\mathcal A}$ is self adjoint and positive, since 
$-1/2<\beta(\lambda_j)<1$, which implies that the eigenvalues of ${\mathcal A}$ lie in the interval $(0,3/2)$. The condition number of ${\mathcal A}$ is of the same order as $1-\rho$, where $\rho$ is the spectral radius of ${\mathcal S}$, i.e. 
by the simple analysis above,
${\rm cond}({\mathcal A})\sim \omega^{2d}$. 
If this system is solved using the (\emph{unconditioned}) conjugate gradient method the convergence rate is $1-1/\sqrt{{\rm cond}({\mathcal A})}\sim 1-1/\omega^d$, \cite{bjorck2015numerical}. Thus, then the method just requires $\sim \omega^d$ iterations for fixed accuracy. 

\begin{remark}
The operator $ {\mathcal A}$ is self-adjoint and coercive
when $\delta>0$ 
since
$$
   \ip{{\mathcal A}u,u} = \ip{(I-{\mathcal S})u,u} = 
\sum_{j=0}^\infty (1-\beta(\lambda_j))|\hat{u}_j|^2
\geq (1-\rho)||u||^2.
$$
This should be contrasted with
the original indefinite Helmholtz problem, which is not coercive. In fact, the eigenvalues satisfy the simple relation $\lambda_{\rm WHI} = 1-\beta(\lambda_{\rm Helmholtz} + \omega)$ > 0. The two formulations are however mathematically equivalent for the interior
Dirichlet and Neumann problems away from resonances,
as the analysis above shows. 

The coercivity also implies that the solution to
\eq{iteration_eq} for the simplified iteration
is unique since $w=\Pi w$ is equivalent to
${\mathcal A}(w-u)=0$.
\end{remark}

\begin{remark}
A discretization would have approximately a fixed number
of grid points per wavelength, leading to a (sparse) matrix
of size $N\times N$ with $N\sim \omega^d$. Hence, the number
of iterations for WHI would be $O(N^2)$
and the total cost $O(N^3)$ since each iteration costs $O(N)$.
This should be compared with a direct solution method which
is better than $O(N^3)$ when the matrix is sparse. 
\end{remark}

\begin{remark}
In the Krylov accelerated case this analysis suggests that the number of iterations would now be $O(N)$ and the total cost $O(N^2)$. However, in the experiments below we observe slightly better complexity for interior problems and significantly better complexity for open problems. In fact, for the open problems we find that, in both two and three dimensions, the number of iterations scale as $\sim \omega$ which is the required number of iterations  for the information to travel through the domain.    
\end{remark}



\subsection{Analysis of the Discrete Iteration} \label{sec:discrete}

To better understand the effects of discretizations 
we consider the following
discrete version of the
algorithm for the energy conserving case described 
above in Section~\ref{sec:iteration}.
We introduce the temporal grid points $t_n=n\Delta t$
and a spatial grid with $N$ points
together with the 
vector $w^n\in {\mathbb R}^N$ containing
the grid function values of the approximation at $t=t_n$.
We also let $f\in {\mathbb R}^N$ hold the corresponding values
of the right hand side.
The discretization of the
continuous spatial operator $-\nabla\cdot(c^2(x)\nabla)$, 
including the boundary conditions,
is denoted 
$L_h$ and it can be represented as an $N\times N$ matrix.
The values $-\nabla\cdot(c^2(x)\nabla w)$ are then approximated by $L_hw^n$.
As in the continuous case, 
we assume $L_h$ has
the eigenmodes $(\lambda_j^2,\phi_j)$,
such that $L_h\phi_j=\lambda_j^2\phi_j$
for $j=1,\ldots,N$, where all $\lambda_j$ are strictly positive and
ordered as
$0\leq \lambda_1\leq\ldots\leq \lambda_N$.

We let the Helmholtz solution $u$ be given
$$
  -L_hu + \omega^2 u = f.
$$
The numerical approximation of the iteration operator is
denoted $\Pi_h$, and it
is implemented as follows. Given $v\in\Real^N$, we use
the leap frog method to solve the wave equation as
\be\lbeq{leapfrog}
   w^{n+1} = 2w^n-w^{n-1} - \Delta t^2 L_h w^n - \Delta t^2 f\cos(\omega t_n),
\ee
with 
time step $\Delta t=T/M$ for some
integer $M$, and
initial data 
$$
   w^0 = v,\qquad w^{-1} = v-\frac{\Delta t^2}{2}(L_hv+f).
$$
The trapezoidal rule is then used to compute $\Pi_h v$,
\be\lbeq{trapzrule}
\Pi_h v
= 
   \frac{2\Delta t}{T}
   \sum_{n=0}^M \eta_n \left(\cos(\omega t_n)-\frac14\right)w^n,
   \qquad \eta_n = \begin{cases}
   \frac12,& \text{$n=0$ or $n=M$},\\
   1, & 0<n<M.
   \end{cases}
\ee

With these definitions we can prove
\begin{theorem}\label{thm:discrete_theorem}
Suppose there are no resonances, such that
$\delta_h=\min_j|\lambda_j-\omega|/\omega>0$.
Moreover, assume that $\Delta t$
satisfies the stability and accuracy requirements
\be\lbeq{cfl}
  \Delta t< \frac{2}{\lambda_N+2\omega/\pi},
  \qquad \Delta t\omega\leq \min(\delta_h,1).
\ee
Then the fixed point iteration $v^{(k+1)}=\Pi_h v^{(k)}$ with
$v^{(0)}=0$ converges to $v^\infty$ which 
is a solution to the discretized Helmholtz equation
with the \reva{modified}{modfied} frequency $\tilde\omega$,
$$
  -L_hv^\infty + \tilde\omega^2 v^\infty = f,\qquad 
  \tilde\omega=2\frac{\sin(\Delta t\omega/2)}{\Delta t}.
$$
The convergence rate is at least $\rho_h=\max(1-0.3\delta_h^2,0.6)$.
\end{theorem}
\begin{proof}
We expand all functions in eigenmodes of $L_h$,
$$
  w^n = \sum_{j=1}^N \hat{w}_j^n \phi_j,
  \qquad
  f = \sum_{j=1}^N \hat{f}_j \phi_j,\qquad
  u = \sum_{j=1}^N \hat{u}_j \phi_j,\qquad
  v = \sum_{j=1}^N \hat{v}_j \phi_j,\qquad
  v^\infty = \sum_{j=1}^N \hat{v}^\infty_j \phi_j.
$$
Then the Helmholtz eigenmodes of $u$ and $v^\infty$ satisfy
$$
\hat{u}_j = \frac{\hat{f}_j}{\omega^2-\lambda_j^2},\qquad
\hat{v}^\infty_j = \frac{\hat{f}_j}{\tilde\omega^2-\lambda_j^2}.
$$
We note that
$\tilde\omega$ is not resonant and
$\hat{v}^\infty_j$ is well-defined
for all $j$, since by
\eq{sincerror} and \eq{cfl}
$$
|\tilde\omega-\lambda_j|\geq |\omega-\lambda_j|-
|\tilde\omega-\omega|\geq \omega\delta_h - \frac{\Delta t^2\omega^3}{24}
\geq
\omega\left(\delta_h - \frac{1}{24}\min(\delta_h,1)^2\right)>0.
$$
The wave solution eigenmodes are given by the difference equation
\be\lbeq{diffeq}
\hat{w}_j^{n+1} - 2\hat{w}_j^n+\hat{w}_j^{n-1} + \Delta t^2\lambda_j^2 \hat{w}_j^n = -\Delta t^2 \hat{f}_j\cos(\omega t_n).
\ee
with initial data
$$
\hat{w}_j^{0} =
\hat{v}_j,\qquad
\hat{w}_j^{-1} =
\hat{v}_j\left(1 - \frac12\Delta t^2\lambda_j^2\right)-\frac12\Delta t^2\hat{f}_j.
$$
By \eq{cfl} 
$$
  |2-\Delta t^2\lambda_j^2|< 2,
$$
and the characteristic polynomial for the equation, $r^2+(\Delta t^2\lambda_j^2-2)r+1$, then has two roots \reva{on the boundary of}{inside} the unit circle.
The solution is therefore stable and is given
by (the verification of which is found in 
Appendix~\ref{sec:verification})
\be\lbeq{discretesol}
   \hat{w}^n_j = 
   (\hat{v}_j-\hat{v}^\infty_j)\cos(\tilde\lambda_j t_n) + \hat{v}^\infty_j\cos(\omega t_n),
\ee
where $\tilde\lambda_j$ is well-defined by the relation
$$
  2\frac{\sin(\Delta t\tilde{\lambda}_j/2)}
  {\Delta t} = \lambda_j.
$$
Now, let 
$$
  \Pi_hv = \sum_{j=1}^\infty    \bar{v}_j \phi_j.
$$
Then the numerical integration gives
\begin{align*}
  \bar{v}_j &= 
   \frac{2\Delta t}{T}
   \sum_{n=0}^M \eta_n\left(\cos(\omega t_n)-\frac14\right)
   \left(
   (\hat{v}_j-\hat{v}^\infty_j)\cos(\tilde\lambda_j t_n) + \hat{v}^\infty_j\cos(\omega t_n)\right)
   \\
   &= 
   (\hat{v}_j-\hat{v}^\infty_j) \beta_h(\tilde\lambda_j)+
   \hat{v}^\infty_j \beta_h(\omega)
=
   \hat{v}_j\beta_h(\tilde\lambda_j) +
   (1-\beta_h(\tilde\lambda_j))
   \hat{v}^\infty_j,
\end{align*}
where
$$
  \beta_h(\lambda) = 
     \frac{2\Delta t}{T}
   \sum_{n=0}^M\eta_n 
   \cos(\lambda t_n)\left(\cos(\omega t_n)-\frac14\right),
$$
and we used the fact that 
the trapezoidal rule is exact, and equal to one, when $\lambda=\omega$.
(Recall that for periodic functions
the trapezoidal rule is exact for all pure trigonometric
functions of order less than the number of grid points.)
Hence, if $|\beta_h(\tilde\lambda_j)|<1$ the
$j$-th mode in the fixed point iteration
 converges to
$\hat{v}_j^\infty$.
This is ensured by the following lemma, the proof
of which is found in Appendix~\ref{sec:filterlemma2}.
\begin{lemma}\lblem{filterlemma2}
Under the assumptions of Theorem~\ref{thm:discrete_theorem},
\begin{equation}\lbeq{betahest}
   \max_{1\leq j\leq N} |\beta_h(\tilde\lambda_j)|\leq \rho_h =:\max(1-0.3\delta_h^2, \reva{0.63}{0.6}).
\end{equation}
\end{lemma}
Since the bound 
$|\beta_h(\tilde\lambda_j)|\leq \rho_h<1$ 
in the lemma is
uniform for all $j$
the convergence
$v^{(k)}\to v^\infty$ with rate at least $\rho_h$ follows.
This concludes the proof of the theorem.

\end{proof}

\begin{remark}
The discretization above is used as an example
to illustrate the impact of going from the continuous
to the discrete iteration. 
For a particular discretization
we can improve the iteration further by using the knowledge
of how it approximates $\omega$ and the eigenvalues of the
continuous operator. Indeed, 
for the discretization above,
let us define $\bar\omega$
by the relation
$$
  \omega=2\frac{\sin(\Delta t\bar\omega/2)}{\Delta t}.
$$
Then if we use 
$f\cos(\bar\omega t_n)$ instead of
$f\cos(\omega t_n)$ in the time stepping \eq{leapfrog}, 
the limit will be
precisely the Helmholtz solution, $v^\infty=u$.
Furthermore, 
the condition
$\Delta t\omega \leq \min(\delta_h,1)$
can be quite restrictive for problems close 
to resonance. It is only important to ensure convergence
of the iterations. Another way to do that is to slightly
change the discrete filter by replacing the constant
$1/4$ in \eq{trapzrule}
by a $\Delta t$-dependent number such that
$|\beta_h(\lambda)|<1$ for $\lambda\neq \omega$.
Another option is to use a higher order quadrature rule,
which would mitigate the restriction on $\Delta t$.
\end{remark}


\subsection{Tunable Filters} \label{sec:time_dep_fil}
In \lem{filterlemma} we saw that the filter transfer function satisfies $\beta(\omega) = 1$ and $-1/2 < \beta(r) < 1$ when $r \neq \omega$ and that these conditions guaranteed convergence of the WaveHoltz iteration. To improve convergence when $r \approx \omega$ we now consider a more general filter transfer function 
\begin{align}
	\bar{\beta}(\lambda) = 
	\frac{2}{T}\int_0^T \left(\cos(\omega t) + \alpha(t)\right) \cos(\lambda t)\, dt,
	\quad
	\alpha(t) = a_0 + \sum_{n=1}^\infty a_n \sin(n\omega t)\reva{,}{.} \label{filter_case2}
\end{align}
\reva{where we refer to $\alpha(t)$ as a \textit{time-dependent shift}.}{}. As before \reva{,}{} necessary conditions for convergence are $\bar{\beta}(\omega) = 1$, $\bar{\beta}'(\omega) = 0$. Straightforward calculations reveal that these conditions require that the two first coefficients must satisfy 
\begin{align*}
 a_1 = \frac{1}{2\pi}(1 + 4a_0).
\end{align*}
The remaining terms in the sum are orthogonal to $\cos(\lambda t )$ when $\lambda = \omega$.
Carrying out the integration in full for each term yields the general 
form
  \begin{align*}
    \bar{\beta}(\lambda) =  
    \frac{\lambda \omega \sin(\lambda T)}{\pi (\lambda^2 - \reva{\omega^2}{\pi^2})} 
    + a_0 \frac{\omega \sin(\lambda T)}{\pi \lambda}  
    + \sum_{n=1}^\infty a_n\frac{n\omega^2}{\pi(\lambda^2 - n^2\omega^2)} 
    \left(\cos(\lambda T) - 1\right),
  \end{align*}
from which it follows that 
another necessary condition is
$|a_0| < 1/2$ since $|\bar{\beta}(r)|<1$ and 
\begin{align*}
 \bar{\beta}(0) = a_0 \lim_{\lambda \rightarrow 0}  \frac{\omega \sin(2 \pi \lambda/\omega )}{\pi \lambda} = 2 a_0. 
\end{align*}
We note that the standard filter, where $a_0=-1/4$ and $a_1=0$,
satisfies the necessary conditions.
\begin{remark}
For the remaining coefficients $a_n$ we only need to ensure that $|\bar{\beta}(r)|<1$ which leaves large freedom to design $\bar{\beta}$. For example we may try to maximize $|\bar{\beta}''(\omega)|$ 
(minimize $\bar{\beta}''(\omega)$)
so that $\bar{\beta}(r)$ is sharply peaked around $r=\omega$. We do not
pursue a systematic study of this here but illustrate the utility of the added flexibility of (\ref{filter_case2}) with numerical experiments below
in Section~\ref{sec:numexp}.
\end{remark}


\subsection{Multiple Frequencies in One Solve}\label{sec:multifreq}
We can use the WaveHoltz algorithm to solve for multiple frequencies at once. 
Suppose we look for the solutions $u_i$ of
$$
  \nabla\cdot(c^2(x)\nabla u_i) + \omega_i^2u_i = f_i(x),\qquad i=1,\ldots, N,
$$
with the same $c$ and boundary condition for all $i$.
To find those solutions we include all frequencies in the
wave equation part of the iteration (\ref{eq:apa}), and solve
  \begin{align}
    w_{tt} = \nabla\cdot (c(x)^2 \nabla w)
    - \sum_{i=1}^N f_i(x)\cos(\omega_i t).
\label{eqn:combined_wave}
  \end{align}
We then seek a decomposition 
\begin{equation}
w(x,t) \equiv \sum_{i=1}^N u_i(x) \cos(\omega_i t), \label{eq:decomp_mul_freq}
\end{equation}
of the solution. The 
filtering part of the
WaveHoltz iteration is also updated to reflect the multiple frequencies 
\begin{align*}
    v_{n+1}  = \frac{2}{T} \int_0^T \left(\sum_{i=1}^N\cos (\omega_i t) - \frac{1}{4} \right) w(x,t) \, dt.
\end{align*}
As before we take $v_0 = 0$ when we deal with energy conserving boundary conditions. To this end we assume that the frequencies are related by an integer multiple in a way so that  the period $T$ can be chosen based on the lowest frequency. 

The different $u_i(x)$ in (\ref{eq:decomp_mul_freq}) can be found as follows. Once we have found the time periodic solution to (\ref{eqn:combined_wave}) evolve one more period and sample $w(x,t)$ at $N$ distinct times $t_j, \, j = 1,\ldots,N$. We then have 
\[
u_i(x) = \sum_{j = 1}^{N} \beta_{ij} w(x,t_j),
\]
where the coefficients $\beta_{ij}$ are the elements of $A^{-1}$ with the elements of $A$ being  $a_{ij} = \cos(\omega_j t_i)$.


\subsection{WaveHoltz Iteration for Impedance Boundary Conditions} \label{sec:imp_iter}
For impedance and other boundary conditions that leads to a decreasing energy for the wave equation we cannot make the simplifying assumption in (\ref{eq:apa}) that $w_t(0,x) = 0$ but we must seek both $v_0(x)$ and $v_1(x)$ in (\ref{eq:apa}). To do so we define an extended iteration \eq{iteration} where we apply $\Pi$ to both the displacement and the velocity: 
\begin{equation}
\left[
\begin{array}{c}
v\\
v' 
\end{array}
\right]^{(n+1)} = \tilde{\Pi} \left[
\begin{array}{c}
v\\
v' 
\end{array}
\right]^{(n)}, \qquad \left[
\begin{array}{c}
v\\
v' 
\end{array}
\right]^{(0)} \equiv 0,
\end{equation}
where
\be
  \tilde{\Pi} \left[
\begin{array}{c}
v\\
v' 
\end{array}
\right] = \frac{2}{T}\int_0^{T}\left(\cos(\omega t)-\frac14\right) \left[
\begin{array}{c}
w(t,x)\\
w_t(t,x) 
\end{array}
\right] dt,\quad
T=\frac{2\pi}{\omega}.
\ee
Here $w(t,x)$ and its time derivative $w_t(t,x)$ satisfies the wave equation (\ref{eq:apa}) with initial data $v_0(x) \equiv v^{(n)}$ and $v_1(x) \equiv v'^{(n)}$.


\section{Wave Equation Solvers}
In this section we briefly outline the numerical methods we use in the experimental section below. We consider both discontinuous Galerkin finite element solvers and finite difference solvers. In all the experiments we always use the trapezoidal rule to compute the integral in the WaveHoltz iteration. 

\subsection{The Energy Based Discontinuous Galerkin Method} \label{sec:dG}
Our spatial discretization is a direct application of the formulation described for general second order wave equations in \cite{Upwind2,el_dg_dath}. Here we outline the spatial discretization for the special case of the scalar wave equation in one dimension and refer the reader to \cite{Upwind2} for the general case. 

The energy of the scalar wave equation is 
\[
H(t) = \int_{D} \frac{v^2}{2} + G(x,w_x) dx,
\]
where 
\[
 G(x,w_x) = \frac{c^2(x)w_x^2}{2},
\] 
is the potential energy density, $v$ is the velocity (not to be confused with the iterates $v^n$ above) or the time derivative of the displacement, $v = w_t$. The wave equation, written as a second order equation in space and first order in time then takes the form 
\begin{eqnarray*}
w_t &=& v, \\
v_t &=& - \delta G \reva{- f(x)\cos(\omega t)}{}, 
\end{eqnarray*}
where $\delta G$ is the variational derivative of the potential energy 
\[
\delta G = - (G_{w_x})_x = -(c^2(x)w_x)_x.
\]
For the continuous problem the change in energy is   
\begin{equation}
\f{d H(t)}{dt} = \int_{D} v v_t + w_t  (c^2(x)w_x)_x \,dx = \reva{-\int_{D} v f(x) \cos(\omega t) dx+}{} [ w_t  (c^2(x)w_x)]_{\partial D}, \label{eq:energy_derivative}
\end{equation}
where the last equality follows from integration by parts together with the wave equation. Now, a variational formulation that mimics the above energy identity can be obtained if the equation $v-w_t=0$ is tested with the variational derivative of the potential energy. Let $\Omega_j$ be an element and $\Pi^s(\Omega_j)$ be the space of polynomials of degree $s$, then the variational formulation on that element is:
\begin{prob}
Find $v^h \in \Pi^s(\Omega_j)$, $w^h \in \Pi^{r}(\Omega_j)$ such that for all
$\psi \in \Pi^s(\Omega_j)$, $\phi \in \Pi^{r}(\Omega_j)$ 
\begin{eqnarray}
\int_{\Omega_j} c^2 \phi_x \left( \f {\pa w^h_x}{\pa t}-v^h_x \right) dx & = & 
[c^2\phi_x \cdot n \left( v^{\ast}-v^h \right)]_{\pa \Omega_j}, 
\label{var1} \\
\int_{\Omega_j} \psi \f {\pa v^h}{\pa t} + c^2 \psi_x \cdot  w^h_x \reva{+ \psi f(x) \cos(\omega t)}{} \, dx& = &
 [\psi \, (c^2\,w_x)^{\ast}]_{\pa \Omega_j}. \label{var2} 
\end{eqnarray}
\end{prob}

Let $[[\zeta]]$ and $\{\zeta\}$ denote the jump and average of a quantity $\zeta$ at the interface between two elements, then, choosing the numerical fluxes as 
\begin{eqnarray*}
v^{\ast}  &=& \{v\} -\tau_1 [[c^2\,w_x]]\\
(c^2\,w_x)^{\ast} &=& \{ c^2\,w_x \}  -\tau_2 [[v]],
\end{eqnarray*}
will yields a contribution $ -\tau_1 ([[c^2\,w_x]])^2 -\tau_2 ([[v]])^2$ from each element face. \reva{To this end we choose $\tau_i > 0$ (so called upwind or Sommerfeld fluxes) which together with the choice that the approximation spaces be of the same degree $r=s$ result in methods that are $r+1$ order accurate in space and measured in the $L_2$ norm. We note that even in the case of energy conserving numerical fluxes the formulation does not lead to a {\em symmetric} matrix for the WaveHoltz iteration (it is of course positive definite though).}{ to the change of the discrete energy}
\[
\reva{}{\f{d H^h(t)}{dt} = \frac{d}{dt} \sum_{j} \int_{\Omega_j} \frac{(v^h)^2}{2} + G(x,w^h_x).}
\]

Physical boundary conditions can also be handled by appropriate specification of the numerical fluxes, see \cite{Upwind2} for details. The above variational formulation and choice of numerical fluxes results in an energy identity similar to (\ref{eq:energy_derivative}). However, as the energy is invariant to certain transformations the variational problem does not fully determine the time derivatives of $w^h$ on each element and independent equations must be introduced. In this case there is one invariant and an independent equation is $\int_{\Omega_j} \left( \f {\pa w^h}{\pa t}-v^h\right) = 0$. 

\reva{Denoting the degrees of freedom on element $\Omega_j$ by $v_j$ and $w_j$ the semi-discretization according to (\ref{var1})-(\ref{var2}) on element $\Omega_j$ can be written}{}
\begin{eqnarray}
\reva{S(\frac{\partial w_j}{\partial t} - v_j) = L_1(v_{j-1},v_j,v_{j+1},w_{j-1},w_j,w_{j+1}),}{}\\
\reva{M\frac{\partial v_j}{\partial t} + Sw_j + f_j \cos(\omega t) = L_2(v_{j-1},v_j,v_{j+1},w_{j-1},w_j,w_{j+1}), }{}
\end{eqnarray}
\reva{where the elements of the element matrices $M$ and $S$ are $M_{kl} = \int_{\Omega_j} \phi_x \phi_l dx$ and  $S_{kl} = \int_{\Omega_j} c^2 (\phi_k)_x (\phi_l)_x dx$ respectively and the lift operators $L_1$ and $L_2$ represents the numerical fluxes. Note that a convenient way to directly enforce the independent equation is to compute the time derivatives of $w_j$ according to}{}
\[
\reva{\frac{\partial w_j}{\partial t} =  v_j = S^{\dagger}L_1(v_{j-1},v_j,v_{j+1},w_{j-1},w_j,w_{j+1}),}{}
\]  
\reva{where $S^{\dagger}$ is the pseudo inverse of $S$.}{}

\reva{}{In this paper we always choose $\tau_i > 0$ (so called upwind or Sommerfeld fluxes) and we always choose the approximation spaces to be of the same degree $r=s$. These choices result in methods that are $r+1$ order accurate in space.} 

\subsection{Finite Difference Discretizations} \label{sec:FD}
For the finite difference examples we exclusively consider Cartesian domains $(x,y,z) \in [L_x,R_x] \times [L_y,R_y] \times [L_z,R_z]$ discretized by uniform grids $(x_i,y_j,z_k) = (L_x+ih_x,L_y+jh_y,L_z+kh_z)$, with $i = 0,\ldots,n_x$ and $h_x = (R_x-L_x) / n_x$, etc. 

When we have impedance boundary conditions on the form $w_t \pm \vec{n} \cdot \nabla w  = 0$ we evolve the wave equation as a first order system in time according to the semi-discrete approximation 
\begin{eqnarray}
\frac{d v_{ijk}(t)}{dt}  &=&  (D^x_+D^x_- + D^y_+D^y_-  + D^z_+D^z_- ) w_{ijk}, \\
\frac{d w_{ijk}(t)}{dt}  &=& v_{ijk},   
\end{eqnarray}
for all grid points that do not correspond to Dirichlet boundary conditions. On boundaries with impedance conditions we find the ghost point values by enforcing (here illustrated on the top of the domain) 
\begin{equation}
v_{ijn_z} - D_0^z w_{ijn_z} = 0.
\end{equation}
Here we have used the standard forward, backward and centered finite difference operators, for example \mbox{$h_x D^x_+ w_{i,j,k} = w_{i+1,j,k} - w_{i,j,k}$} etc. For problems with variable coefficients the above discretization is generalized as in \cite{AppPet09}. 

We note that in some of the examples where we require high order accuracy we use the summation by parts discretization for variable coefficients developed by Mattson in \cite{Mattsson2012} and described in detail there. 

\subsection{Time Discretization}
In most of the numerical examples we use either an explicit second order accurate centered discretization of $w_{tt}$ (for finite differences with energy conserving boundary conditions we eliminate $v$ and time discretize $w_{tt}$ directly as in the analysis in Section~\ref{sec:discrete}) or the classic fourth order accurate explicit Runge-Kutta method.      

For some of the DG discretizations we employ Taylor series time-stepping in order to match the order of accuracy in space and time.  Assuming that all the degrees of freedom have been assembled into a vector ${\bf w}$ we can write the semi-discrete method as ${\bf w}_t = Q {\bf w} $ with $Q$ being a matrix  representing the spatial discretization. Assuming we know the discrete solution at the time $t_n$ we can advance it to the next time step $t_{n+1} = t_n + \Delta t$ by the simple formula
\begin{eqnarray*}
{\bf w}(t_{n}+\Delta t) &=& {\bf w}(t_{n}) + \Delta t {\bf w}_t(t_{n}) +  \frac{(\Delta t)^2}{2!}{\bf w}_{tt}(t_{n}) \ldots \\
 &=& {\bf w}(t_{n}) + \Delta t Q {\bf w}(t_{n}) +  \frac{(\Delta t)^2}{2!} Q^2 {\bf w}(t_{n}) \ldots 
\end{eqnarray*} 
The stability domain of the Taylor series which truncates at time derivative number $N_{\rm T}$ includes \reva{part of}{} the imaginary axis if $ {\rm mod} (N_{\rm T},4) = 3$ or ${\rm mod} (N_{\rm T},4) = 0$ \reva{(see e.g. \cite{ketcheson2015absolute})}{}. However as we use a slightly dissipative spatial discretization the spectrum of our discrete operator will be contained in the stability domain of all sufficiently large choices of $N_{\rm T}$ (i.e. the $N_{\rm T}$ should not be smaller than the spatial order of approximation). Note also that the stability domain grows linearly with the number of terms.


\section{Numerical Examples}\label{sec:numexp}
In this section we illustrate the properties of the proposed iteration and its Krylov accelerated version by a sequence of numerical experiments in one, two and three dimensions.

\subsection{Examples in One Dimension}
We begin by presenting some very basic numerical experiments in one dimension.
\subsubsection{Convergence of Different Iterations / Solvers at a Fixed Frequency}
We start by repeating the example described in Section 3.5 in \cite{2019arXiv190312522G}. This example is used in \cite{2019arXiv190312522G} to illustrate that the original cost functional from \cite{bristeau1998controllability}  (denoted $J$ in \cite{2019arXiv190312522G}) does not yield the 
correct solution due to the existence of multiple minimizers.  

The example solves the Helmholtz equation with $c= 1$ and with the exact solution  
\[
u(x) = 16x^2(x-1)^2, \ \ 0 \le x \le 1.
\]
Here both $u$ (and $w_t$ for the time-dependent problem) and $u_x$ vanish at the endpoints so any boundary condition of the form 
\[
\alpha w_t + \beta (\vec{n}\cdot w_x) = 0, \ \ \alpha^2+\beta^2=1, 
\]
will be satisfied. Dirichlet boundary conditions correspond to $\alpha = 1$ and Neumann boundary conditions correspond to $\alpha = 0$, all other values will be an impedance boundary condition. Here, as in \cite{2019arXiv190312522G}, we take the frequency to be $\omega = \pi/4$. 

We discretize using the energy based DG method discussed above and use upwind fluxes which adds a small amount of dissipation. For this experiment we use 5 elements with degree $q = 7$ polynomials and we use an 8th order accurate Taylor series method in time.  We set  $\Delta t$ so that $n_{t} \Delta t = T = 2\pi/\omega$ while making the inequality $ \Delta t \le C_{\rm CFL} \Delta x /(q+1)$ as sharp as possible (in this experiment we fix $C_{\rm CFL} =1/2$). With this resolution in space and time the truncation errors are negligible and we expect that the observed convergence properties should match those of the continuous analysis. 
\begin{table}[]
\begin{center}
\begin{tabular}{|c|c|c|c|c|c|}
\hline
Method / b. c. & WHI & LSQR & QMR & CG & GMRES\\
\hline
\hline
 D-D & 94.5(-15)  &  76.1(-15) &   75.9(-15) &  151.5(-15)  &  97.9(-15)\\
N-N & 49.2(-15)  &   142.8(-15) &    144.4(-15)  &   158.5(-15)  &   144.1(-15) \\ 
 D-N  & 28.3(-15)    &  55.4(-15)    &  81.9(-15)   &  272.3(-15)    &  67.0(-15) \\
 \hline
 \end{tabular}
\caption{Maximum error for various combinations of boundary conditions and methods. \label{tab:1d_conv_err}}
\end{center}
\end{table}%

\begin{figure}[]
\graphicspath{{figures/}}
\begin{center}
\includegraphics[width=0.32\textwidth]{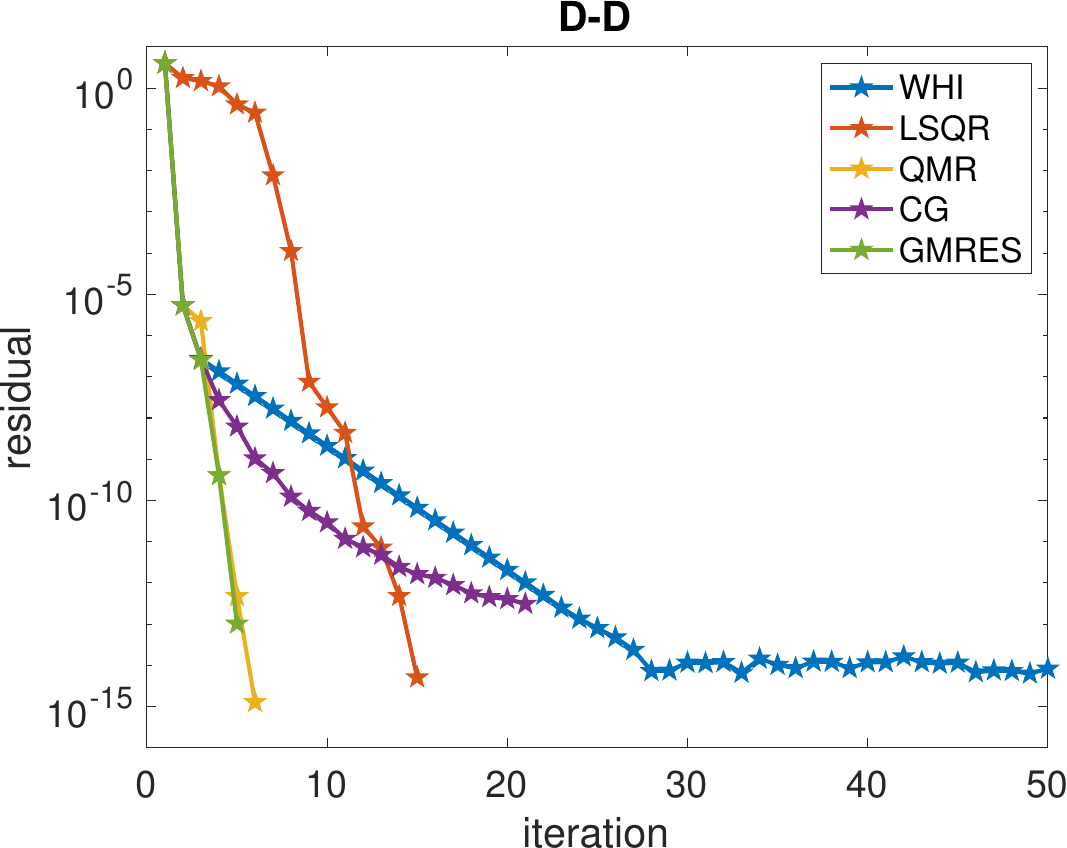}
\includegraphics[width=0.32\textwidth]{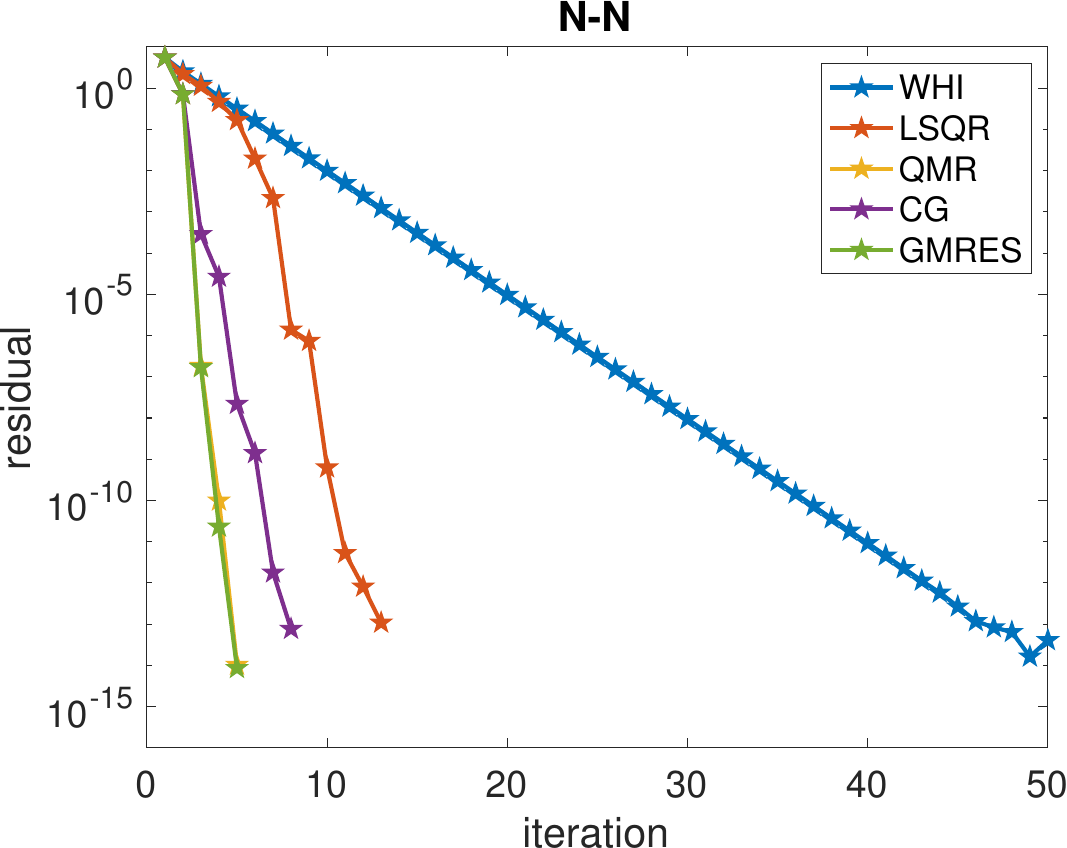}
\includegraphics[width=0.32\textwidth]{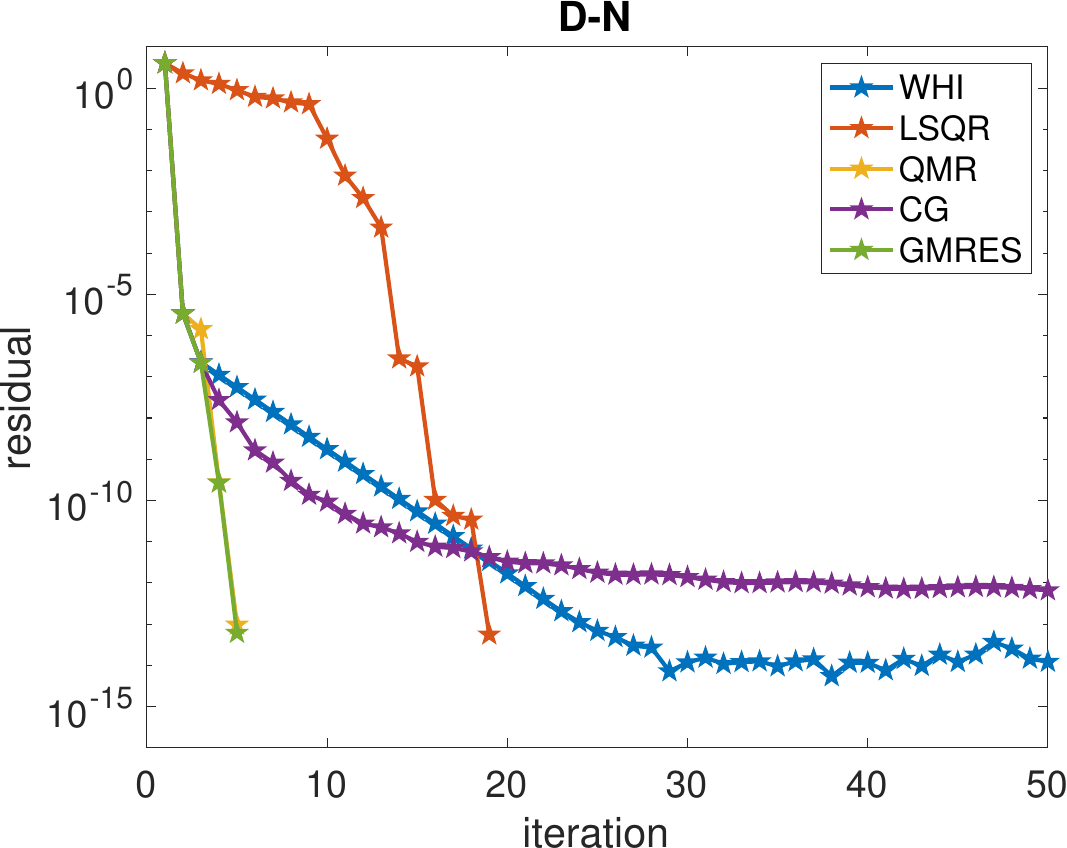}
\caption{Convergence of the residual for the plain WaveHoltz iteration and its accelerated versions using LSQR, QMR, CG and GMRES. The titles of the figures indicate the boundary conditions used to the left and right, e.g. D-N means Dirichlet on the left and Neumann on the right.  \label{fig:1d_bc_test}}
\end{center}
\end{figure}

As mentioned above we expect that our method works best when combined with a classical iterative Krylov subspace method. The energy based DG method will produce a matrix $A$ with real eigenvalues in $(0,3/2)$ but it will not yield a symmetric matrix $A$. We present results for the WaveHoltz iteration (denoted WHI in figures and tables), and its acceleration with Matlab implementations of  LSQR, QMR, CG and GMRES (we use the default unconditioned settings with a tolerance of $10^{-13}$). In Figure \ref{fig:1d_bc_test} we display the convergence histories for various combinations of boundary conditions. The residuals for the Krylov accelerated iterations are the ones returned by the Matlab functions and the residual for the WaveHoltz iteration is simply the $L_2$ norm of the difference between two subsequent iterations. As can be seen the convergence behavior for QMR and GMRES are uniformly the fastest and appears to be insensitive to the type of boundary condition used. Note that the numerical method used here does not yield a symmetric matrix and CG is not guaranteed to work. Evidence of this loss or stagnation of convergence can be found in the cases D-D and D-N in Figure \ref{fig:1d_bc_test}.

The actual errors in the converged solutions can be found in Table \ref{tab:1d_conv_err}, where it can be seen that the error for all of the iteration methods are close to the residual tolerance.   

\subsubsection{Convergence with Increasing Frequency}
To study how the number of iterations scale with the Helmholtz frequency $\omega$ we solve the wave equation on the domain $x\in[-6,6]$ \reva{with constant wave speed $c^2(x) = 1$}{} and with a forcing 
\[
f(x) = \omega^2 e^{-(\omega x)^2},
\]
that results in the solution being $\mathcal{O}(1)$ for all $\omega$. The solver is the same as in the previous example. We keep the number of degrees of freedom per wave length fixed by letting the number of elements be $5 \lceil \omega \rceil$. We always take the polynomial degree to be 7 and the number of Taylor series terms in the timestepping to be 8. As we now also consider impedance boundary conditions, with $\alpha = 1/2$, we use WHI accelerated by GMRES. 
\begin{figure}[]
\graphicspath{{figures/}}
\begin{center}
\includegraphics[width=0.31\textwidth]{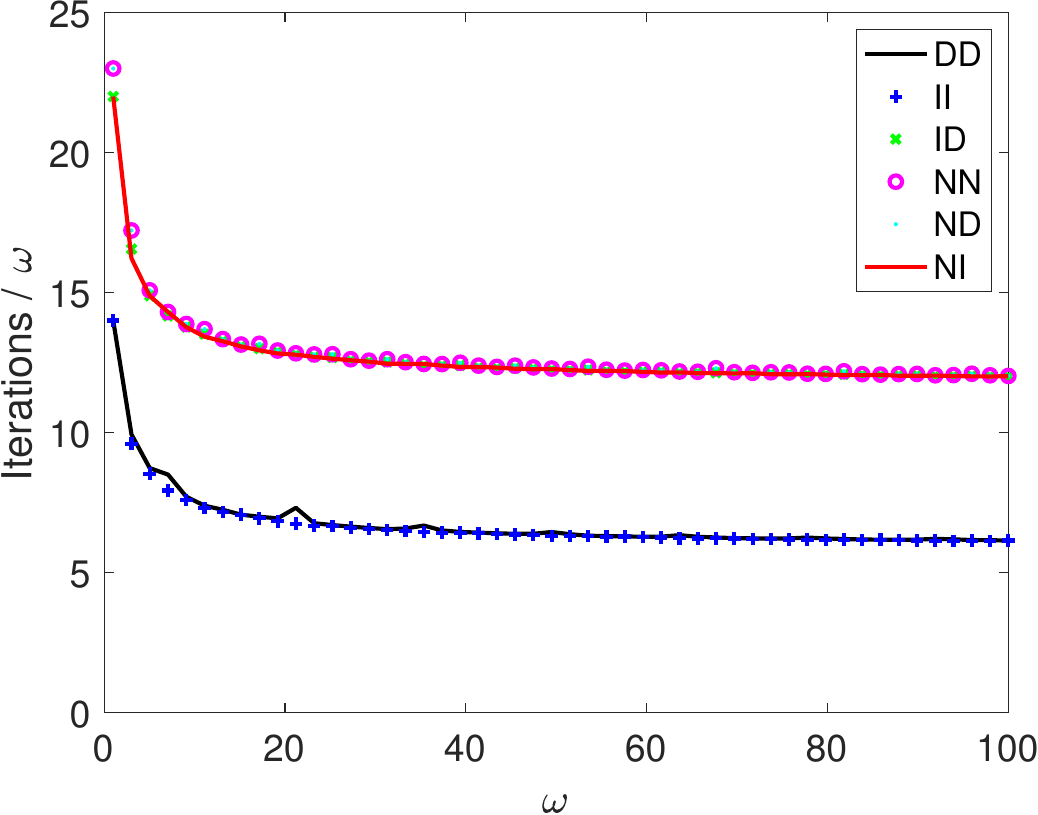}
\includegraphics[width=0.32\textwidth]{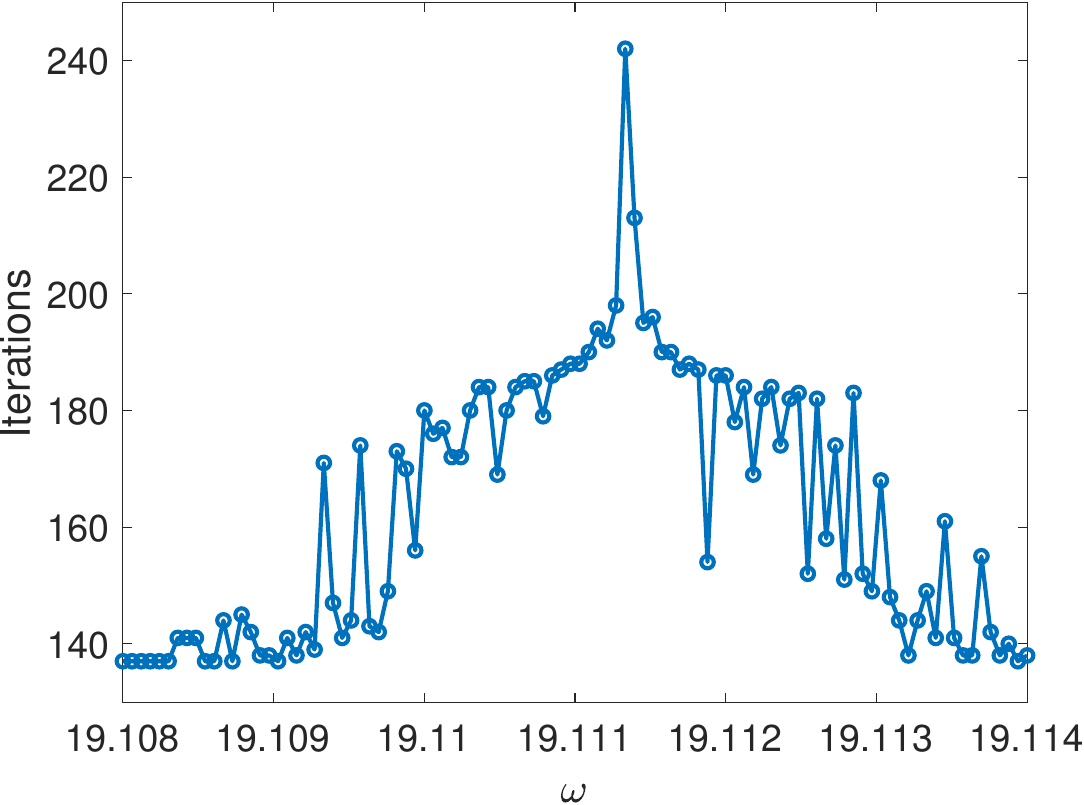}
\includegraphics[width=0.32\textwidth]{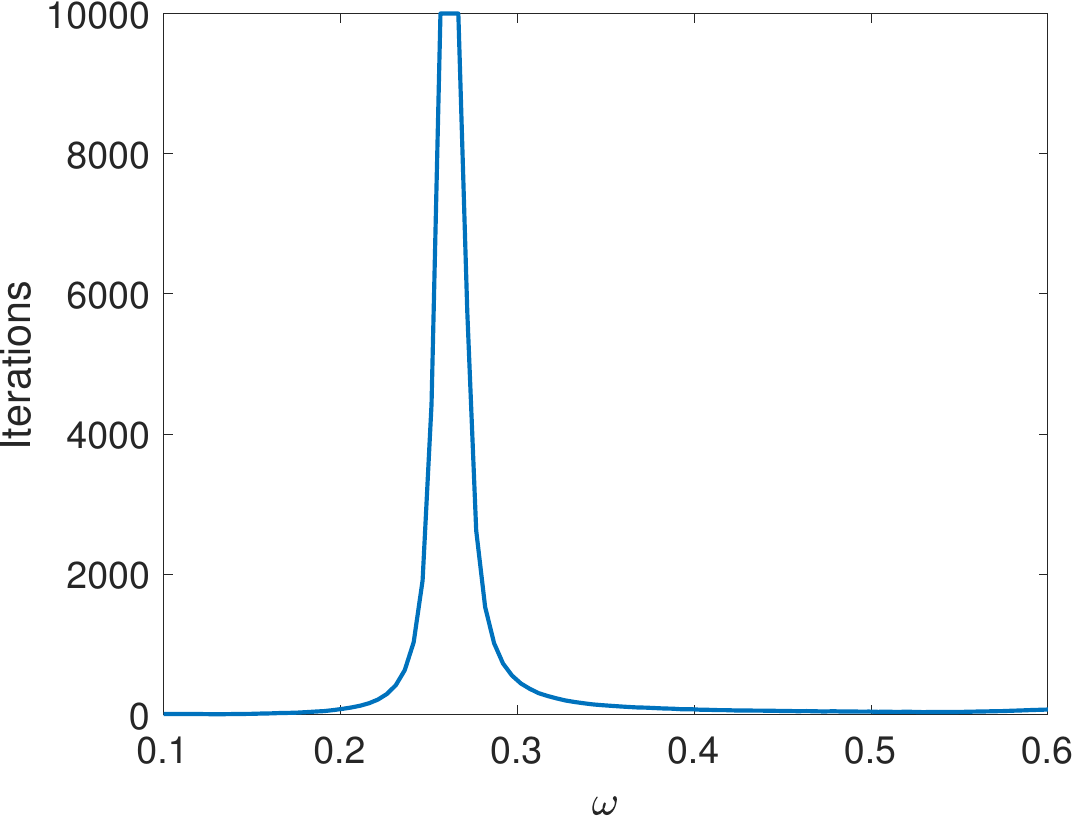}

\caption{Left: Number of iterations divided by $\omega$ as a function of $\omega$ for different boundary conditions. Middle and right: Zoom in around a resonance for the Dirichlet problem when using Krylov acceleration (middle) and when using WHI (right).   \label{fig:1d_iter_vs_freq}}
\end{center}
\end{figure}

We report the number of iterations it takes to reach a GMRES residual smaller than $10^{-10}$ for the six possible combinations of Dirichlet, Neumann and impedance boundary conditions for 50 frequencies distributed evenly between 1 and 100. The results are displayed to the left in Figure \ref{fig:1d_iter_vs_freq} where we plot the number of iterations divided by $\omega$ as a function of $\omega$. It is clear that the asymptotic scaling is linear with growing frequency. Interestingly all the combinations of boundary conditions collapse to two different curves with the Dirichlet-Dirichlet and impedance-impedance conditions converging the fastest.    

We know from the analysis in Section \ref{sec:iteration} that the rate of convergence of the WaveHoltz iteration deteriorates near resonant frequencies (for non-impedance problems) but from Figure \ref{fig:1d_iter_vs_freq} it appears that all frequencies converge more or less the same rate. To study the behavior of the accelerated algorithm    
for homogenous Dirichlet boundary conditions we zoom in around  $\omega \approx 19.114$ where the continuous problem has a resonance. In the middle graph in Figure  \ref{fig:1d_iter_vs_freq} we display the required number of iterations around the resonant frequency. As can be seen there is some deterioration but only in very narrow band around a frequency that is slightly less than $19.114$ and probably is the modified resonant frequency discussed in Section \ref{sec:discrete}. 
This behavior can be contrasted to the growth of the number of iterations for the WaveHoltz iteration without GMRES acceleration, 
see the right figure in Figure \ref{fig:1d_iter_vs_freq}.   
Clearly the acceleration of the WaveHoltz iteration by GMRES improves the robustness of the method near resonances.

\subsubsection{Multiple Frequencies in One Solve}
Here we illustrate the technique described in Section~\ref{sec:multifreq}
for finding solutions of multiple frequencies at once.
We set $\omega = 1$ and $\omega_j = 2^{j-1} \omega$ for $j = 1,\ldots,4$, 
and consider the domain $x\in [0,1]$. We use the finite difference discretization discussed in Section \ref{sec:FD} \reva{with Dirichlet boundary conditions}{}. The time evolution is done by a second order centered discretization of $w_{tt}$\reva{, as was done in the discrete analysis in Section~\ref{sec:discrete}}{}. 
The problem is forced by a point source centered at $x=1/2$ for $ j= 1, \dots, 4$, \reva{and we consider a constant wave speed $c^2(x) = 1$}{}. We display the convergence with decreasing $h$ in Figure \ref{fig:1d_time_dep_filt_convergence} \reva{on the right,}{} where it can be seen that each solution $u_j$ converges at a rate of $h^2$. 

\subsubsection{Tunable Filters}
\begin{figure}[]
\graphicspath{{figures/}}
\begin{center}
\includegraphics[width=0.33\textwidth]{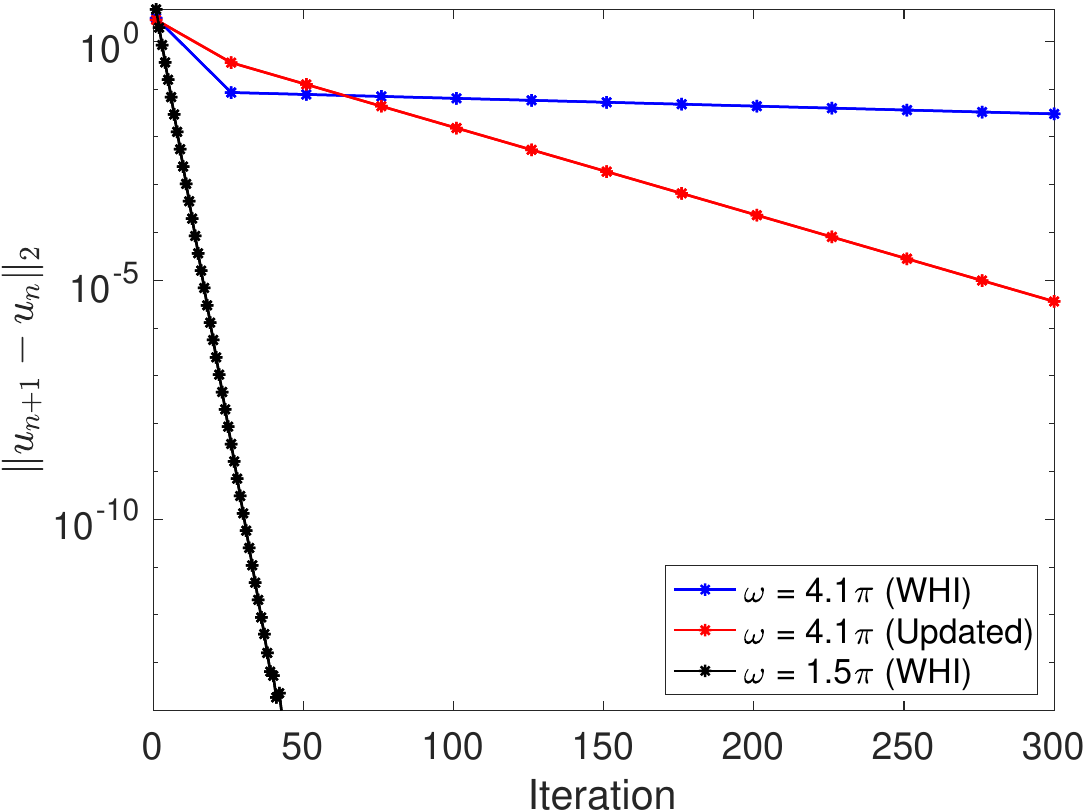}
\includegraphics[width=0.31\textwidth]{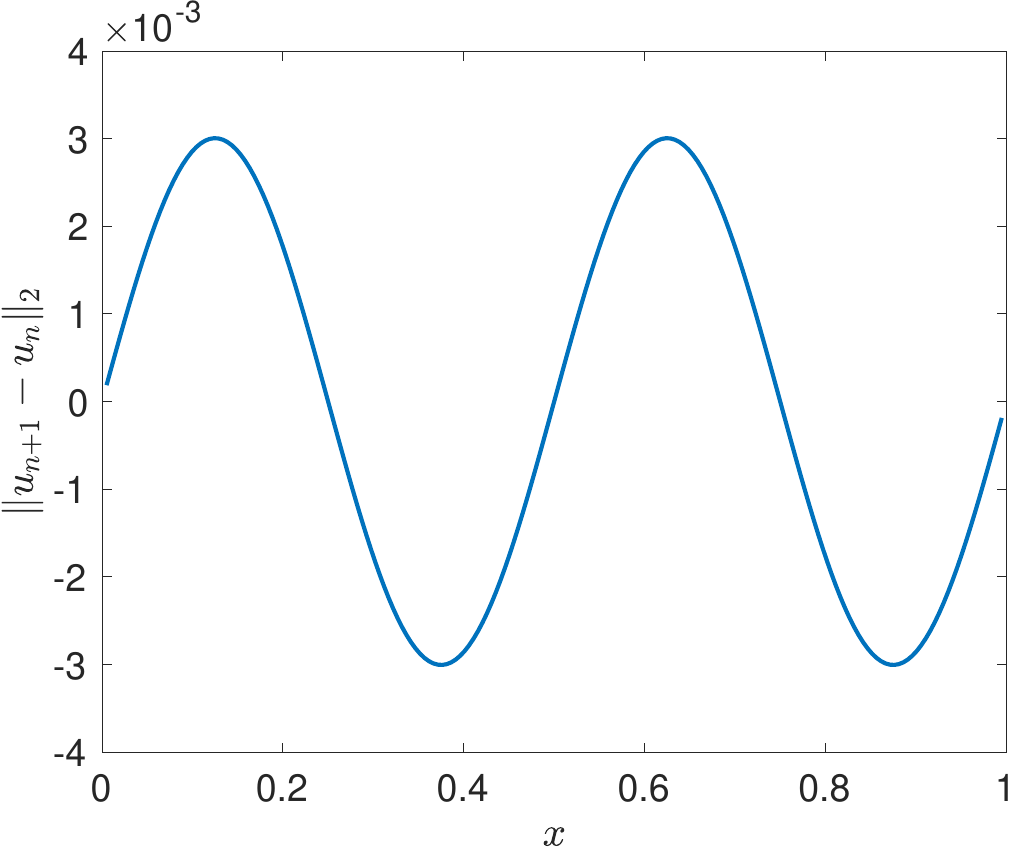}
\includegraphics[width=0.33\textwidth]{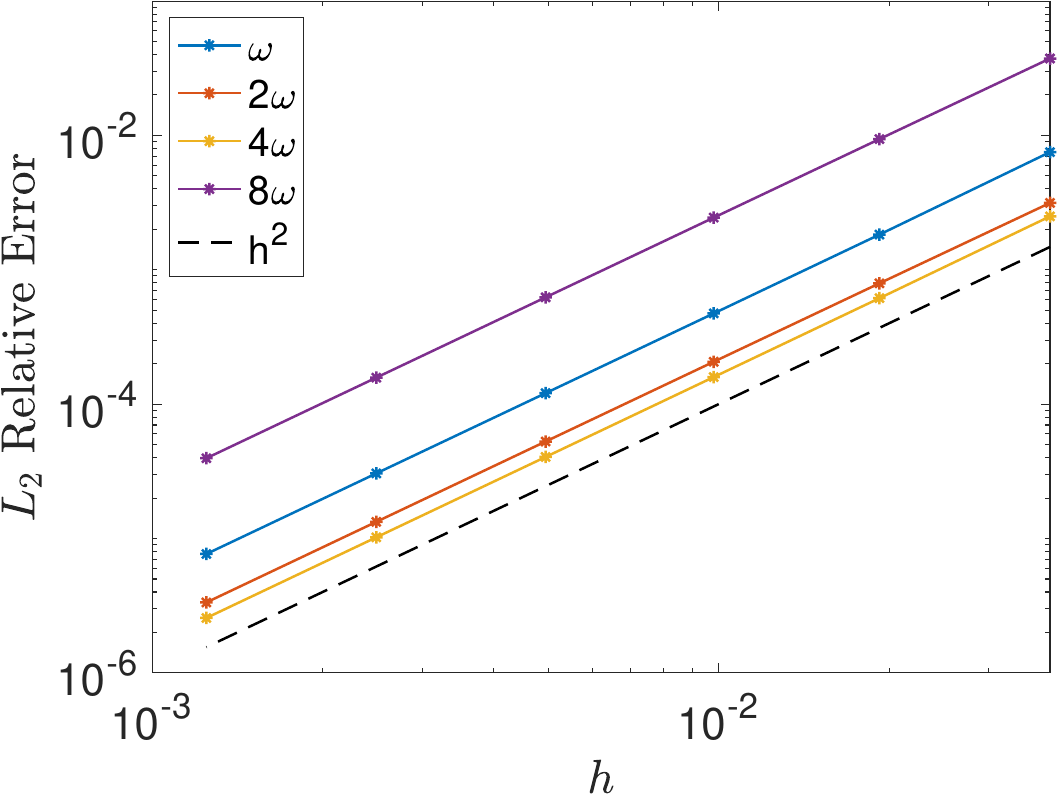}
\caption{Left: Convergence history of the near resonant frequency $4.1\pi$ for the WaveHoltz filter and a tunable filter, and that of the frequency $1.5\pi$ for reference. Middle: The error between successive WaveHoltz iterates with the usual WaveHoltz filter. Right: Convergence of the solution for the CG accelerated  WaveHoltz iteration with a point forcing. \label{fig:1d_time_dep_filt_convergence}}
\end{center}
\end{figure}

Here we consider solving a Helmholtz problem in the domain $x \in [0,1]$ with Dirichlet boundary conditions and 
constant wave speed
$c^2 = 1$. The discretization is the same as in the previous experiment and we use a point source centered at $x=1/2$. A straightforward calculation shows that the resonant frequencies of the problem are integer multiples of $\pi$ and we specifically consider solving the Helmholtz problem with frequency $\omega = 4.1\pi$, which has a minimum relative gap to resonance of $\delta = \reva{1/41 \approx 0.024}{0.025}$. As discussed previously, we expect that the convergence rate of the WaveHoltz iteration will stagnate since $\omega$ is close to resonance. We compare the convergence against the problem with frequency $\omega = 1.5 \pi$ which has a minimum relative gap to resonance of $\delta = \reva{1/3}{0.5}$. The iteration history is displayed in Figure \ref{fig:1d_time_dep_filt_convergence} \reva{on the left}{}.
It can be seen that the usual WaveHoltz iteration converges rapidly for the frequency $1.5\pi$ but that of $4.1\pi$ stagnates considerably. In \reva{the middle of}{} Figure \ref{fig:1d_time_dep_filt_convergence} we display the difference between successive WaveHoltz iterates for the Helmholtz problem with frequency $\omega$\reva{, from which it}{. It} is clear that the residual is a scaling of the resonant mode $\sin(4\pi x)$.

\begin{figure}[]
\graphicspath{{figures/}}
\begin{center}
\includegraphics[width=0.33\textwidth]{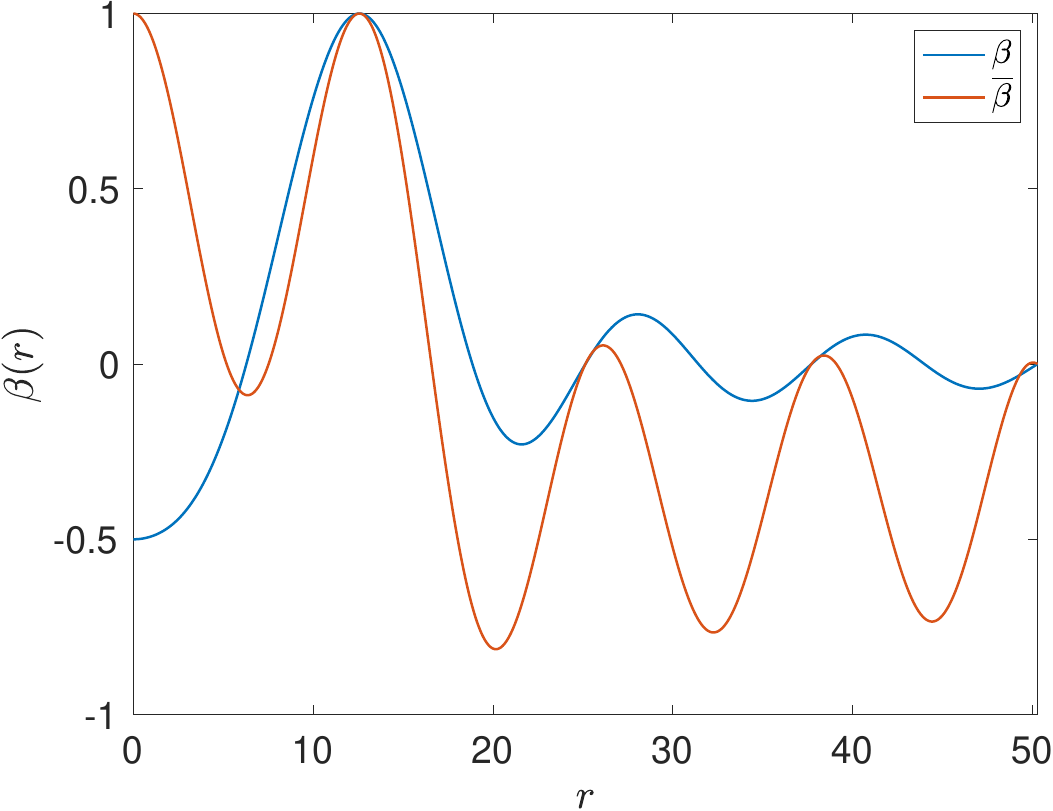}
\includegraphics[width=0.33\textwidth]{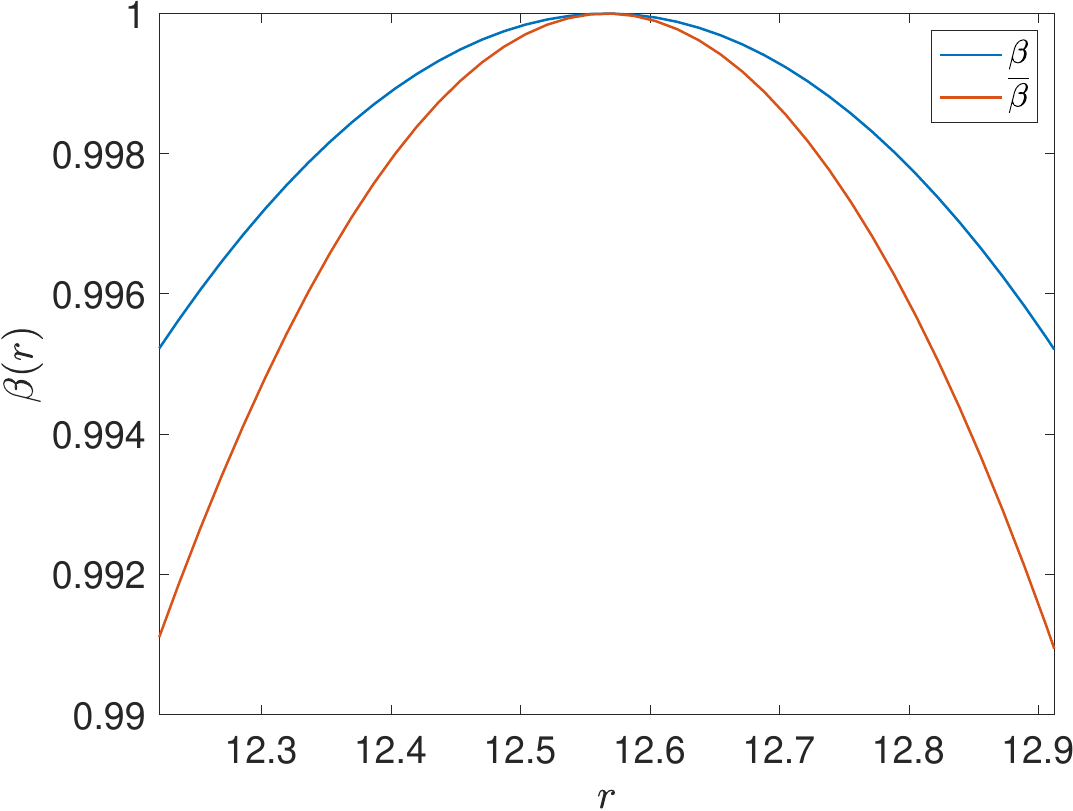}
\caption{ (Left) The usual WaveHoltz filter (in blue) and updated tunable filter (in red). (Right) Closeup of both the usual WaveHoltz filter and the updated tunable filter near 
the resonant frequency $4\pi$.\label{fig:1d_time_dep_filt_plot}}
\end{center}
\end{figure}

\reva{To improve the rate of convergence close to resonance we leverage a tunable filter as mentioned in Section \ref{sec:time_dep_fil}. To obtain this filter, we consider the filter transfer function \eqref{filter_case2} and truncate the $\sin$ expansion of the time-dependent shift $\alpha(t)$ to $12$ terms such that $a_n = 0$ for $n >11$. In this example we take the usual choice of $a_0 = -1/4$ for the constant term in the filter transfer function \eqref{filter_case2} which, as discussed in Section \ref{sec:time_dep_fil}, requires $a_1 = 0$. We then perform a minimization over a discrete set of $3000$ equispaced points $r_j \in [0,16\pi]$ of the empirically constructed functional
\begin{align}
    J(a_2,a_3,\dots,a_{11}) = 10.6 \bar{\beta}''(\omega) + 0.1 \sum_{|r_j - \omega| > 0.1} |\bar{\beta}(r_j)|^{20}, 
    \label{eqn:costFunctional}
\end{align}
via $100$ steepest descent iterations. The first term in the functional \eqref{eqn:costFunctional} minimizes the second derivative at the peak $\omega = 4.1\pi$, while the second weakly enforces that $|\beta(r)| \le 1$ 
for all $r > 0$ to ensure convergence of the fixed point iteration.
}
{To improve the rate of convergence close to resonance, we consider a tunable filter as mentioned in Section \ref{sec:time_dep_fil}. To obtain this filter, we consider a time dependent shift with $12$ terms and perform a minimization of the empirically constructed functional}
\reva{}{
\begin{align}
    J(\alpha_2,\alpha_3,\dots,\alpha_{12}) = 10.6 \beta''(4\pi) + 0.1 \sum_{j \in \mathcal{A}} |\beta(r_j)|^{20}, 
    \label{eqn:costFunctional2}
\end{align}
}
\reva{}{where $\mathcal{A}$ is the index set such that $|r_j - 4\pi| > 0.1$ where $r_j$ are equispaced points between $[0,16\pi]$. The first term in the functional \eqref{eqn:costFunctional} minimizes the second derivative at the peak, while the second weakly enforces that $|\beta(r)| \le 1$ 
for all $r > 0$.}

In Figure \ref{fig:1d_time_dep_filt_plot} \reva{on the right}{} we see that the \reva{updated}{} filter is steeper \reva{near $\omega = 4.1\pi$ so that repeated application of the updated filter will more quickly remove the resonant mode with frequency $4\pi$}{} and we \reva{thus}{} expect faster convergence. This is confirmed in the resulting iteration history of the updated filter, shown in Figure \ref{fig:1d_time_dep_filt_convergence} \reva{on the left}{}. The cost of improving convergence behavior near resonance, however, is a larger value of $\bar{\beta}$ for many other modes as shown in Figure \ref{fig:1d_time_dep_filt_plot} \reva{on the left}{}. A more careful investigation of optimized filters is left for the future.

\subsection{Problems in Two Dimensions}
In this section we present experiments in two space dimensions. 
\subsubsection{Convergence in Different Geometries}
In this example we solve the Helmholtz equation with a constant wave speed, $c^2=1$, in the domain $(x,y) \in [-1,1]^2$ and with forcing 
\[
f(x,y) = -\omega^2  e^{-\sigma [(x-0.01)^2+(y-0.015)^2]},
\] 
where $\sigma = \max(36,\omega^2)$\reva{}{`}. We vary the frequency according to $\omega =  1/2+ k, \, k = 1,\ldots,100,$ and keep the number of points per wavelength roughly constant by choosing  $n_x = n_y = 8 \ceil*{\omega}$. Here we use the finite difference method outlined in Section \ref{sec:FD} combined with the classic fourth order Runge-Kutta method in time with a timestep $\Delta t = h_x / c$. 

\begin{figure}[]
\graphicspath{{figures/fd_2d/}}
\begin{center}
\includegraphics[width=0.32\textwidth]{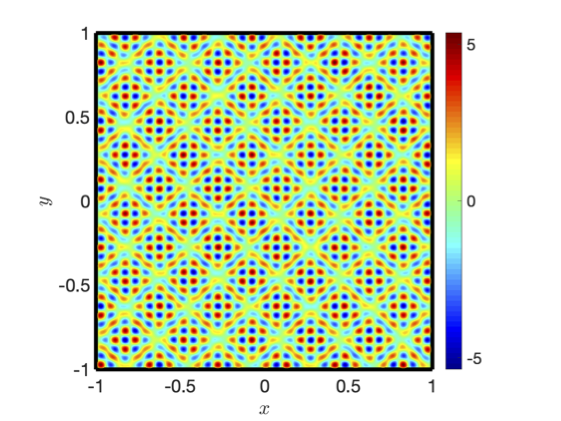}
\includegraphics[width=0.32\textwidth]{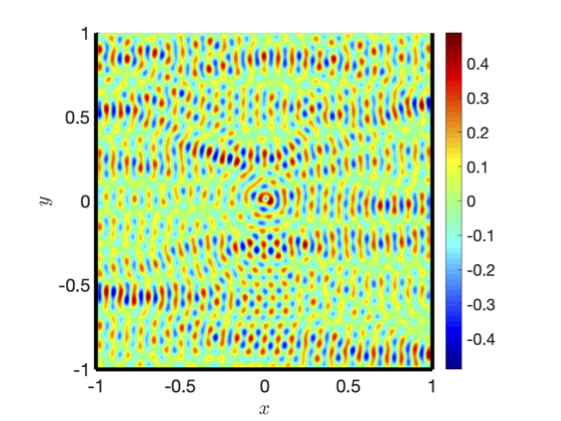}
\includegraphics[width=0.32\textwidth]{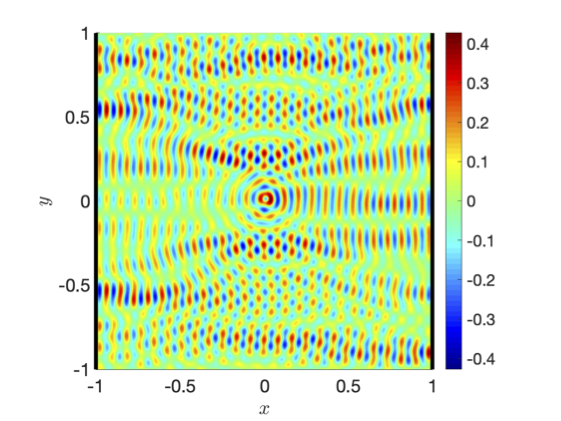}
\includegraphics[width=0.32\textwidth]{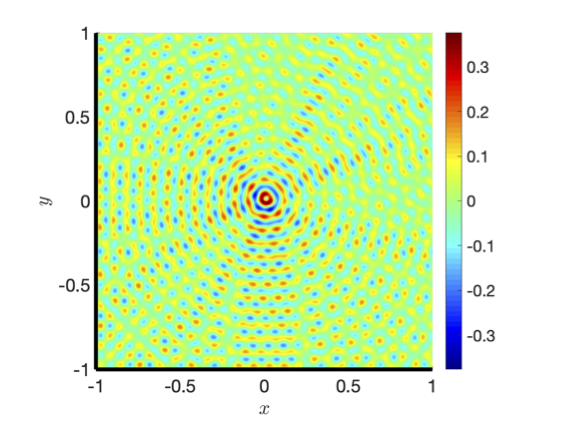}
\includegraphics[width=0.32\textwidth]{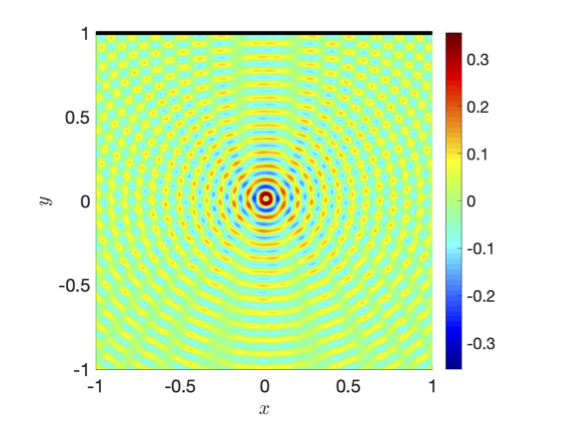}
\includegraphics[width=0.32\textwidth]{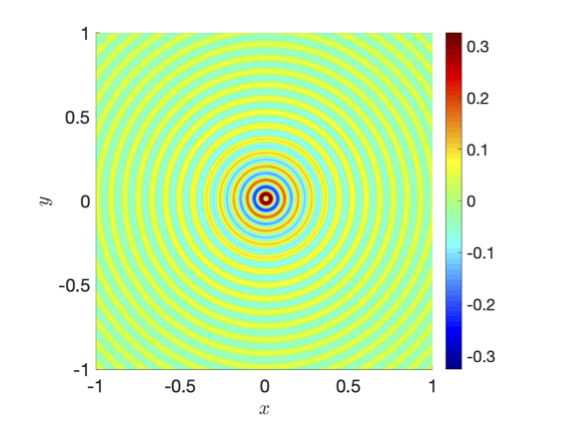}
\caption{Typical solutions computed with the GMRES accelerated WHI at $\omega = 77.5$. The thick lines indicate Dirichlet boundary conditions. \label{fig:2D_box_with_different_sides}}
\end{center}
\end{figure}

For each frequency we solve six different problems consisting of combinations of Dirichlet and impedance boundary conditions with zero to four open sides and with the two open boundary case forking into two cases: (1) the open boundaries are opposite each other, or (2) next to each other forming a corner. In Figure \ref{fig:2D_box_with_different_sides} we display the real part of the solution for the frequency $\omega = 77.5$ for the six different problems.

\begin{figure}[]
\graphicspath{{figures/fd_2d/}}
\begin{center}
\includegraphics[width=0.32\textwidth]{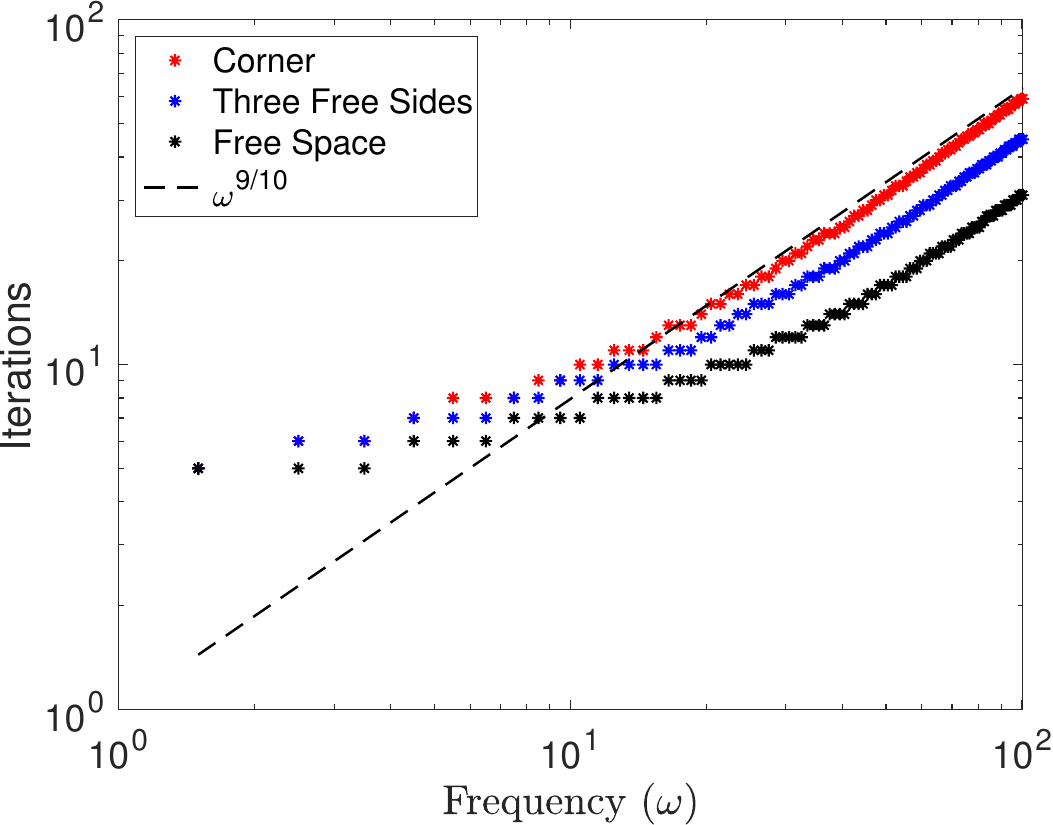}
\includegraphics[width=0.32\textwidth]{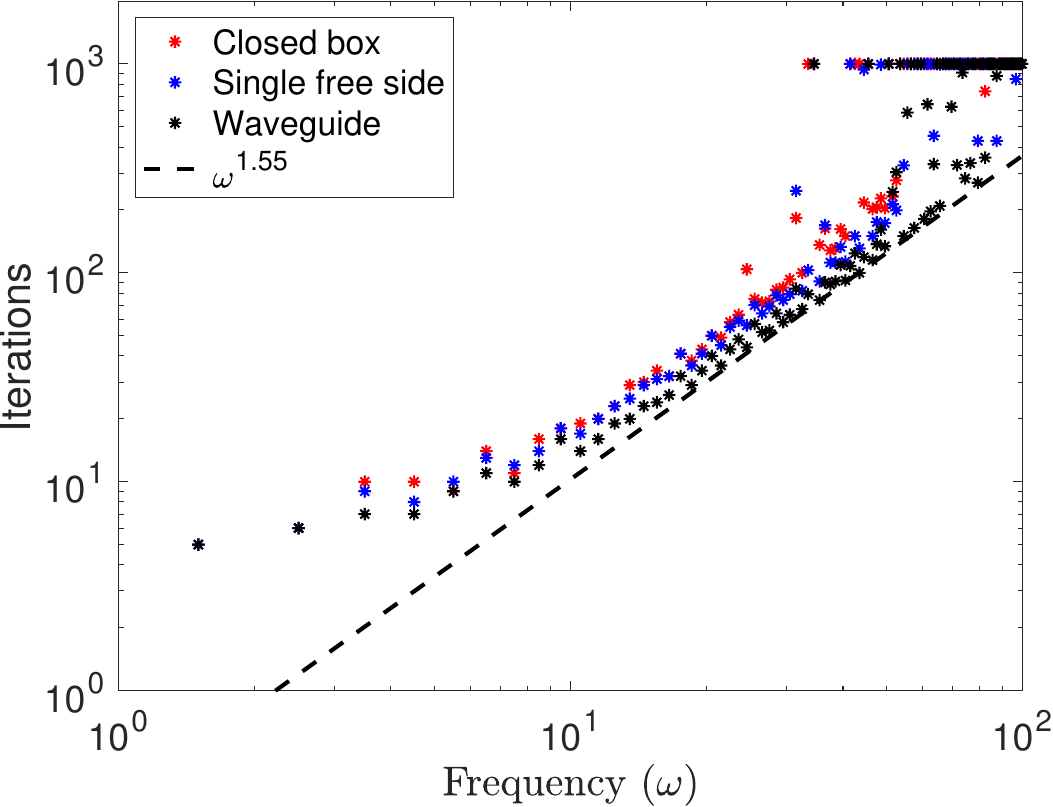}
\includegraphics[width=0.32\textwidth]{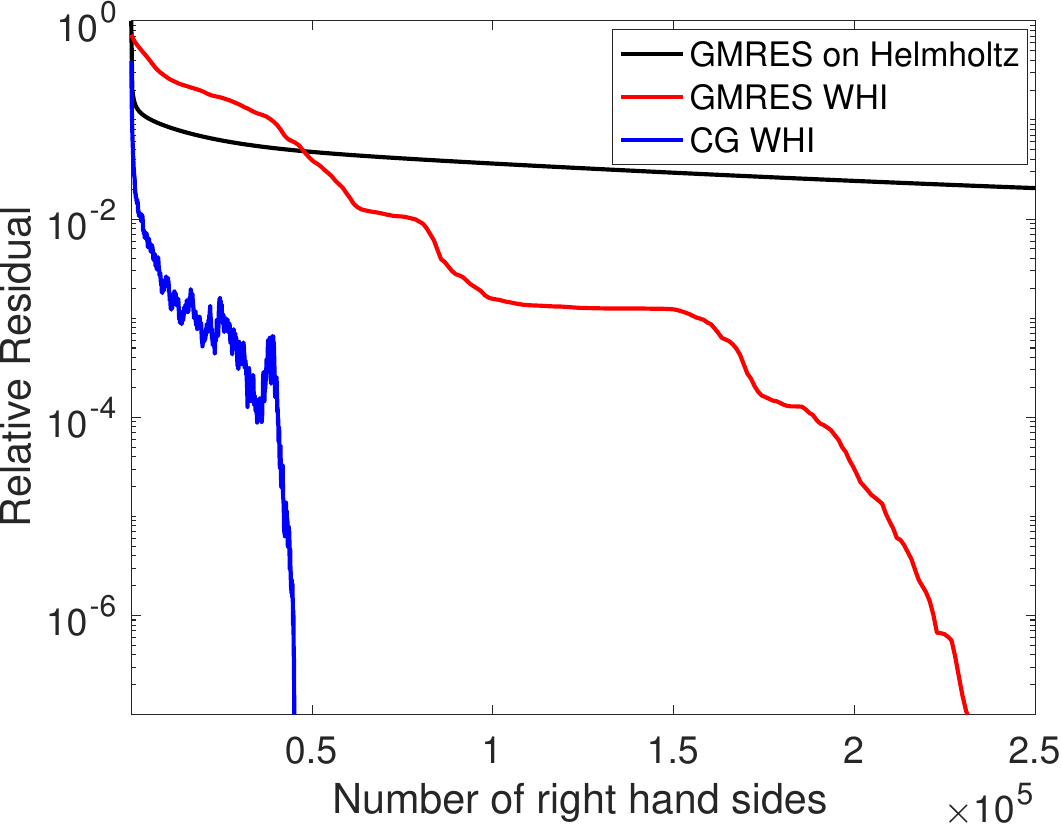}
\caption{To the left: number of iterations as a function of frequency to reduce the relative residual below $10^{-7}$ for problems with no trapped waves.  Middle: the same but for problems with trapped waves and for the interior problem. Both are with the \reva{GMRES accelerated WHI}{WHI and with GMRES}. To the right: Residuals for the \reva{GMRES accelerated WHI}{WHI with GMRES}, the \reva{CG accelerated WHI}{WHI with CG} and for GMRES solution of the directly discretized Helmholtz problem. \label{fig:2D_box_FD}}
\end{center}
\end{figure}

\reva{In this example, the WaveHoltz iteration is accelerated by GMRES without restarts. Given that the storage requirement for GMRES grows with the number of iterations, it is often beneficial (especially for high frequency problems) to integrate and average over several periods to allow further propagation of information within the domain while mitigating the rapid growth of the Krylov subspace. For this example we thus choose to perform the WaveHoltz iteration with an integration time of 10 periods (i.e. we choose $T = 10 \frac{2\pi}{\omega}$).}{In this example the WaveHoltz iteration is performed with an integration time of 10 periods (i.e. we choose $T = 10 \frac{2\pi}{\omega}$) and is accelerated by GMRES without restarts. 
} In Figure \ref{fig:2D_box_FD} we report the number of iterations needed to reduce the relative residual below $10^{-7}$. It is clear from the results that the geometries where the waves can get trapped are considerably more difficult and requires more iterations. The computational results appear to indicate that the number of iterations to reach the tolerance scale as $\omega^{1.55}$ for the inner Dirichlet problem and similarly for the waveguide and the case with three Dirichlet boundary conditions. As the frequency increases and the distance between resonant frequencies decreases the iteration is not able to reduce the relative residual below the tolerance $10^{-7}$ within the prescribed maximum 1000 iterations.  
On the other hand for geometries with no trapped waves we see faster convergence (see the left figure in Figure \ref{fig:2D_box_FD}) with the number of iterations scaling roughly as $\omega^{9/10}$. 

To the right in Figure \ref{fig:2D_box_FD} we display the residual as a function of the number of right hand side evaluations (for the wave equation this is equivalent to taking a timestep and for the direct discretization of Helmholtz this is equivalent to one application of the sparse system matrix, the cost of these are roughly equivalent) when $\omega = 51.5$ for the pure Dirichlet boundary condition problem. The three different results are for: 1. the  WaveHoltz acceleration, 2. the WaveHoltz iteration accelerated with conjugate gradient and based on the same spatial discretization but with a second order accurate centered discretization of $w_{tt}$ using $\Delta t = 0.7 h_x$ and, 3. a direct discretization of the Helmholtz equation (using the spatial discretization described in Section \ref{sec:FD}) combined with GMRES for solving the resulting system of equations. Precisely we use GMRES with restart every 100 iterations. For space reasons we only display this for one  frequency but note that although the results may differ a bit between frequencies the trend is similar in the problems we have investigated. 

It is clear from the residuals that both the GMRES and conjugate gradient  accelerated WaveHoltz iterations are radically faster than applying GMRES to the direct discretization of Helmholtz. As all the methods use the same spatial discretization this is an indication of the importance of changing the problem from an indefinite system of equations to a positive definite and to a symmetric positive definite system. 

\begin{remark}
We note that the problems considered in this experiment can be naturally solved with integral equation techniques since the the wave speed is constant. In addition as the problem is posed in two dimensions and can be stored in memory a good sparse solver will also be a very good alternative. What we want to demonstrate is: 1. The positive definiteness of the accelerated WHI makes it faster than standard iterative techniques for the direct discretization of Helmholtz, 2. The complexity is different for open and closed problems as predicted by the theory in \cite{no_low_rank_Engq}.    
\end{remark}

\subsubsection{Smoothly Varying Wave Speed in an Open Domain} \label{sec:lens}
\begin{figure}[]
\graphicspath{{figures/}}
\begin{center}
\includegraphics[width=0.36\textwidth]{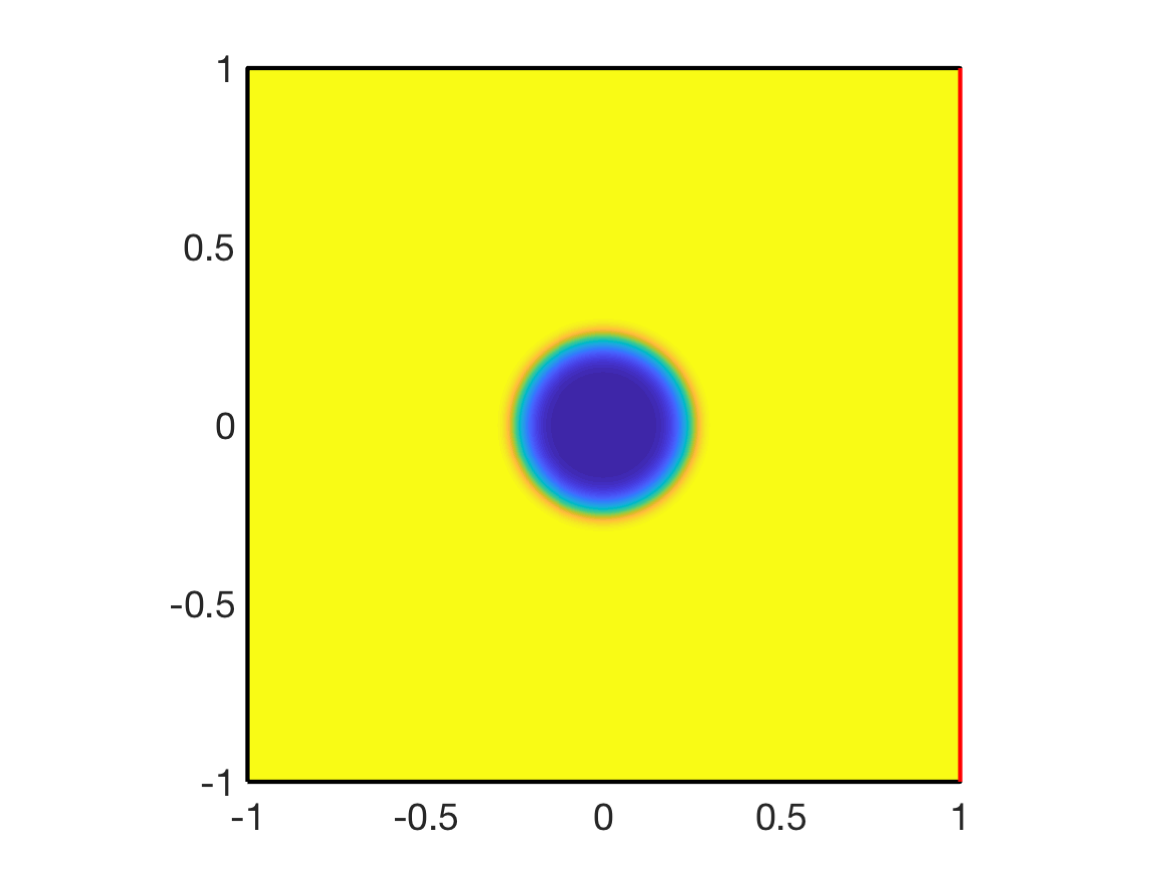}
\includegraphics[width=0.25\textwidth]{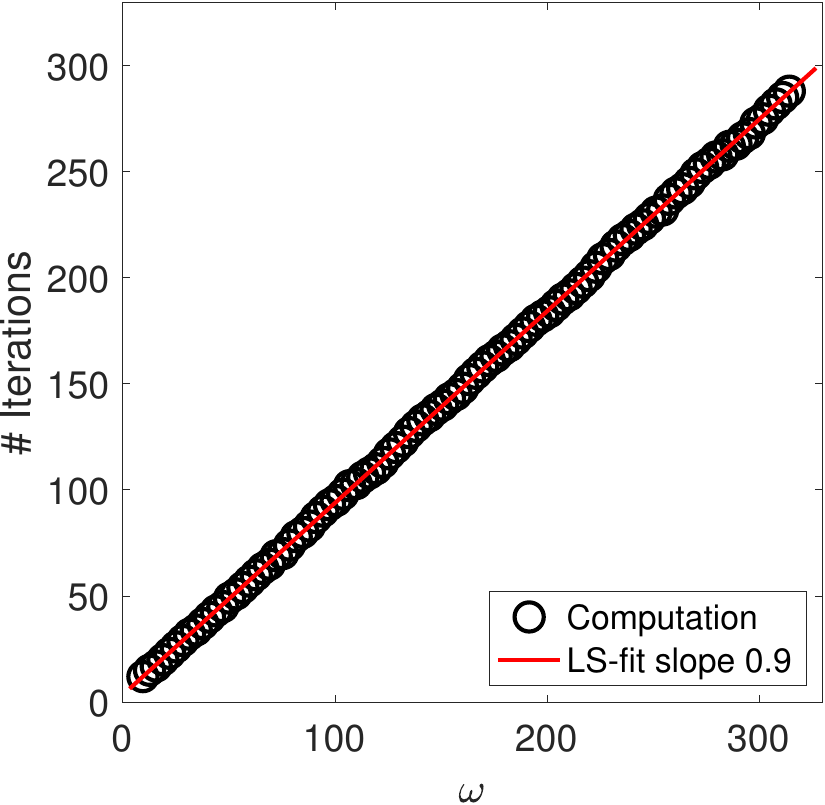}  \ \ \ \ \ \ \ \
\includegraphics[width=0.31\textwidth]{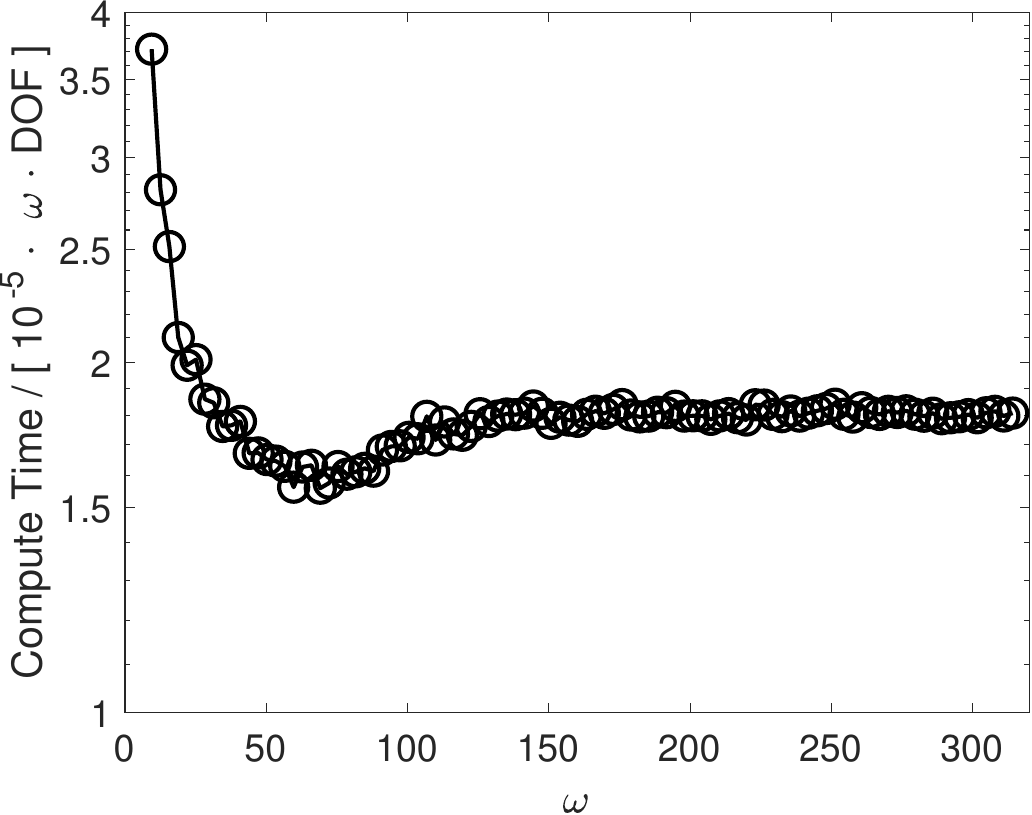}
\caption{Left: the speed of sound (squared) used in example \ref{sec:lens}. Red indicates a rigid wall and black indicates open walls. Middle: Number of iterations as a function of frequency. Right: Compute time normalized by the frequency times the number of degrees of freedom. \label{fig:2d_lens_complexity}}
\end{center}
\end{figure}

In this example we consider a smoothly varying medium in a box $(x,y) \in [-1,1]^2$. The wave speed is   
\[
c^2(x,y) = 1 - 0.4 e^{-\left( \frac{x^2+y^2}{0.25^2} \right)^4},
\]
and is also depicted in Figure \ref{fig:2d_lens_complexity}.

Here we use the energy based DG solver and impose a right going plane wave $e^{i \omega(t-x)}$ through impedance boundary conditions on the left, bottom and top faces of the domain. On the right boundary we impose a zero Dirichlet condition. 
In all the computations we use degree 5 polynomials and a 6th order Taylor series method. The elements used form a Cartesian structured grid and we scale the number of elements so that we have 8 degrees of freedom per wavelength. The  WHI is applied with an integration time of 5 periods and is accelerated by GMRES with a termination tolerance $10^{-7}$ on the relative residual. 
\begin{figure}[]
\graphicspath{{figures/}}
\begin{center}
\includegraphics[width=0.32\textwidth]{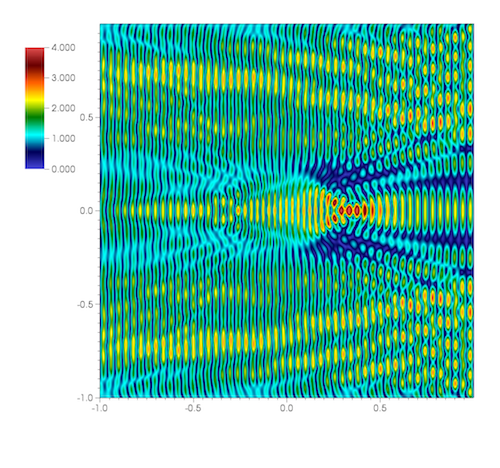}
\includegraphics[width=0.32\textwidth]{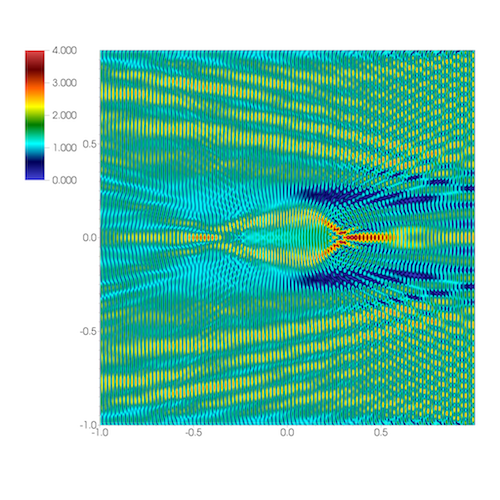}
\includegraphics[width=0.32\textwidth]{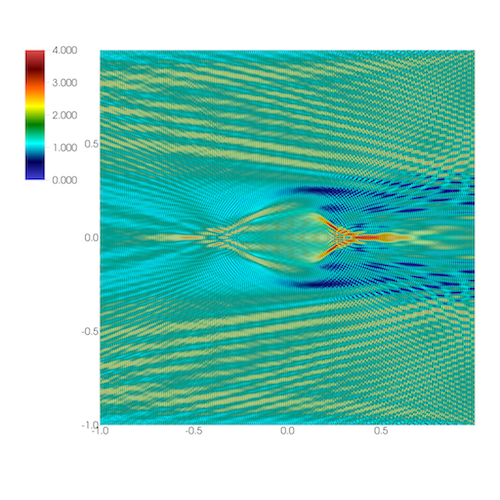}
\caption{The magnitude of the Helmholtz solution for, from left to right, $\omega = 25\pi, 50\pi$ and $100\pi$. \label{fig:2d_lens_fields}}
\end{center}
\end{figure}
We solve the Helmholtz problem with $\omega = k \pi, k=  3,4,\ldots,100$ and measure the total time from start to time of solution and we also measure the number of iterations needed to converge. The results, displayed in Figure \ref{fig:2d_lens_complexity}, again show that for this type of open problem the iteration appears to require $N_{\rm iter} \sim \mathcal{O}(\omega^{0.9})$ iterations to converge to a fixed tolerance. In terms of total computational time we observe $T_{\rm Total} \sim \mathcal{O}( \omega N_{\rm DOF})$ which is slightly higher than what would be expected from the $\mathcal{O}(\omega^{0.9})$ behavior.   

However, as the distance traveled by the wave solution is proportional to $cT = 2\pi c / \omega$ and the information must travel through the domain at least once the time to solution is as good as can be expected. To reduce the computational complexity further we would need to propagate the solution faster than the speed of sound by applying a preconditioner or some type of multi-level strategy. Although we believe this is possible we leave such attempts to future work.      

The magnitude of the solutions with $\omega = 25\pi, 50\pi$ and $100\pi$ are plotted in Figure \ref{fig:2d_lens_fields}. 

\subsubsection{Convergence of the Approximation Error and the Residual}
\begin{figure}[htb]
\graphicspath{{figures/}}
\begin{center}
\includegraphics[width=0.45\textwidth]{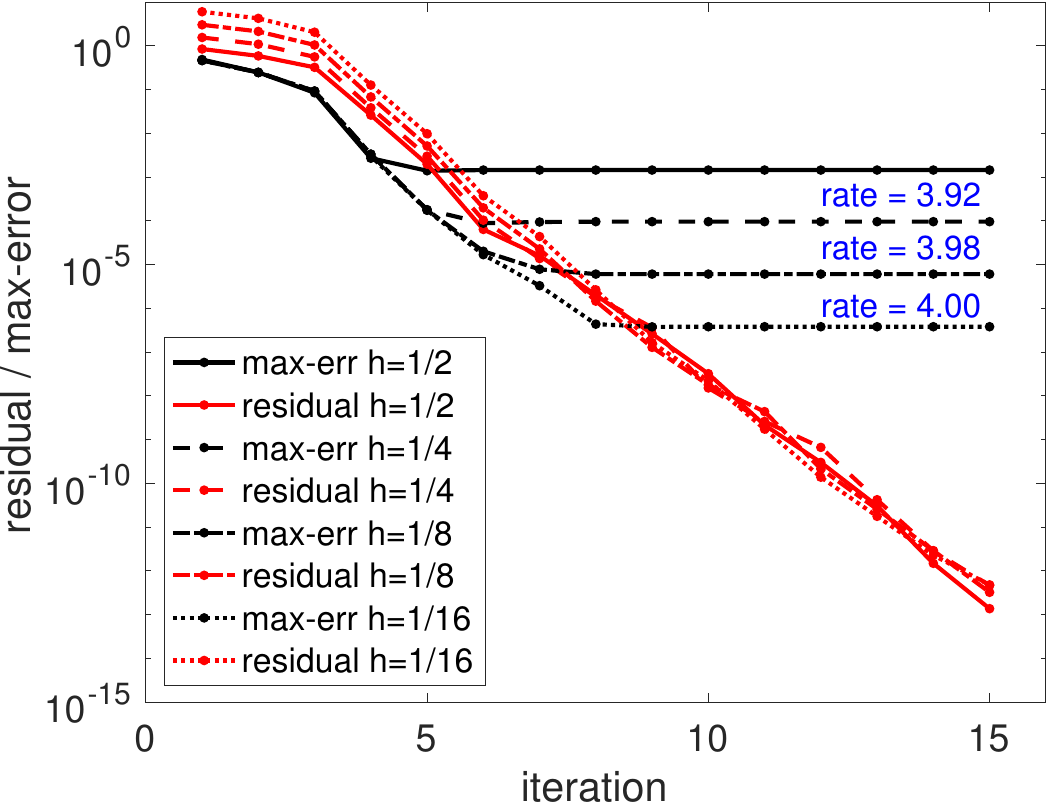}
\caption{The maximum error GMRES residuals as a function of number of iterations for four different mesh sizes. The rates of convergence agree with the order of the method.  \label{fig:err_and_res}}
\end{center}
\end{figure}

\revb{As our iteration leads to a linear system of equations (and consequently a different residual) we should check that the residual is still a suitable proxy for the discretization error. Although we have no reason to believe this would not be the case we note that we have not yet performed a detailed analysis and resort to checking this numerically.  We consider the same computational domain and method as above but with speed of sound $c=1$ and with zero Dirichlet boundary conditions. We set $\omega = 2$ and choose the forcing so that the solution is 
\[
u = -(x^2-1)^2(y^2-1)^2, 
\]
and compute the solution using polynomials of degree three in the energy DG method and a fourth order accurate Taylor time stepper. In Figure \ref{fig:err_and_res} we display the maximum errors in $u$ and the residuals for each GMRES iteration for Cartesian grids with grid spacings $1/2, 1/4, 1/8$ and $1/16$. As can be seen the residuals and the errors track well until the errors saturate. To the right in the figure we also indicate the rates of convergence based on the subsequent grid refinements. As expected they are very close to four.}{}

\subsubsection{The Marmousi2 Model}
In the last two examples in this section we use the sixth order summation-by-parts finite difference operators developed by Mattson in \cite{Mattsson2012}. Here we use the classic fourth order Runge-Kutta method for timestepping. In this example we simulate the solution caused by a point source placed in a material model where the speed of sound is taken from P-wave velocity in the Marmousi2 model\footnote{\url{http://www.agl.uh.edu/downloads/downloads.htm}}. We discretize the full model which consists of $13601 \times 2801$ grid points and covers a domain that is roughly $17\times 3.5$ kilometers. On the top surface we prescribe a zero Dirichlet condition and on the remaining three sides we add a 50 grid point wide supergrid layer (see \cite{AppCol08}) that is terminated by zero Dirichlet boundary conditions. We accelerate the  WHI by the transpose free quasi minimal residual (TFQMR) method and terminate the iteration when the relative residual is below $10^{-5}$.
We perform each iteration over 8 periods and take $500$ timesteps per iteration. The time periodic point forcing is applied near the surface in grid point $(6750,2600)$ and we perform computations with $\omega = 200, 400  $ and 800.

{\newcommand{\figWidth}{4.4cm}
\newcommand{\trimfig}[2]{\trimh{#1}{#2}{0.111}{0.08}{.09}{.12}}
\begin{figure}[]
\begin{center}
\begin{tikzpicture}[scale=1]
\useasboundingbox (0.,0.) rectangle (16.,10.5);  
\begin{scope}[xshift=-0.0cm,yshift=-0.9cm]
\begin{scope}[xshift=0.0cm]
\draw ( 0.0, 0.1) node[anchor=south west,xshift=-8pt,yshift=-4pt] {\trimfig{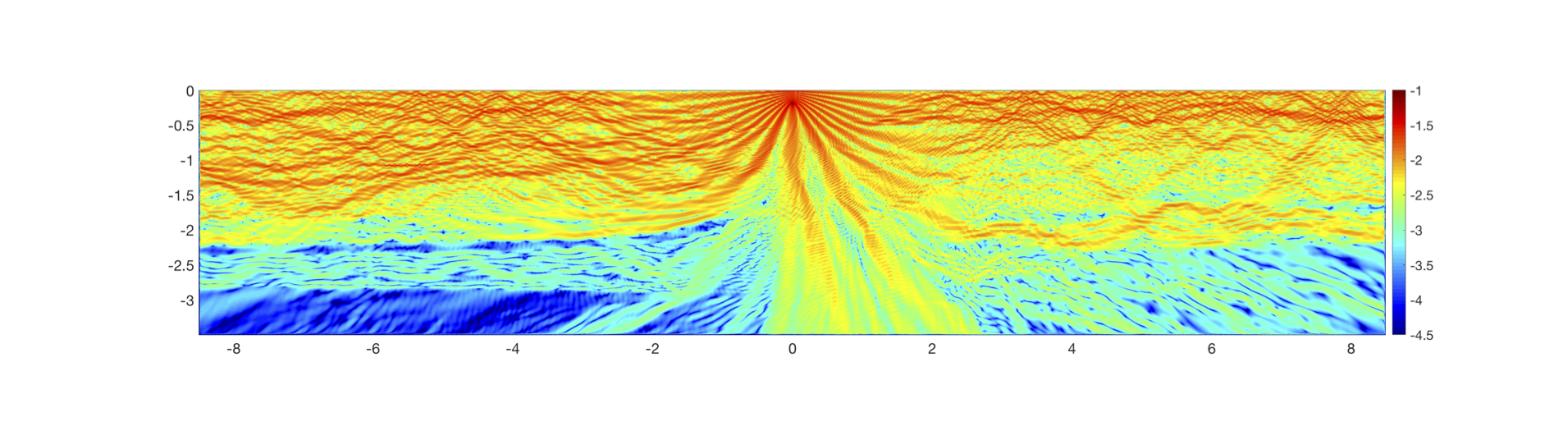}{\figWidth}};
\draw ( 0.0, 4.1) node[anchor=south west,xshift=-8pt,yshift=-4pt] {\trimfig{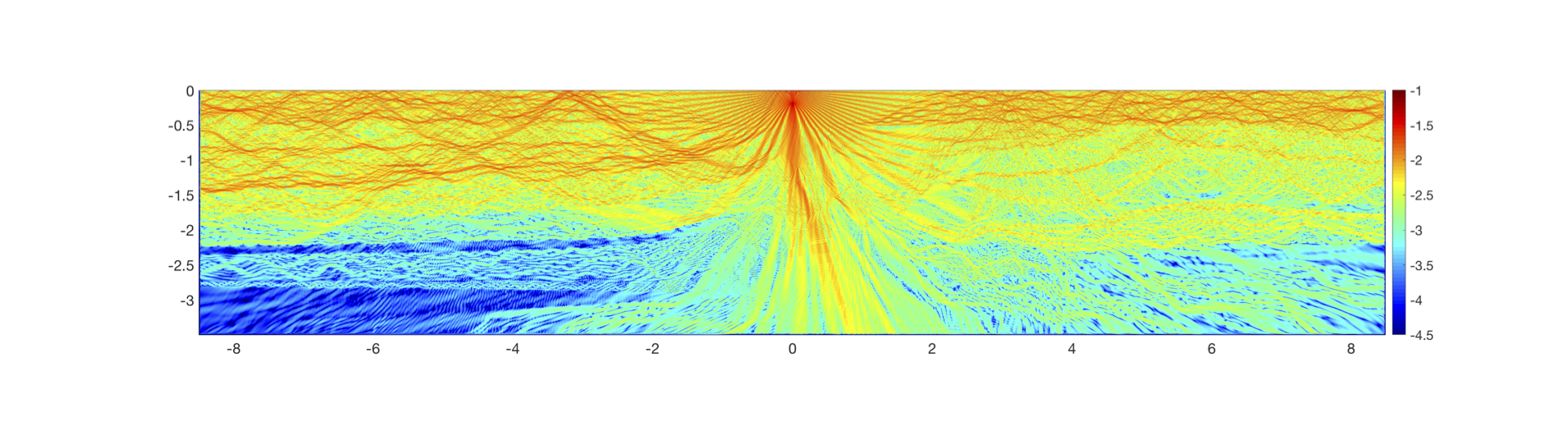}{\figWidth}};
\draw ( 0.0, 8.1) node[anchor=south west,xshift=-8pt,yshift=-4pt] {\trimfig{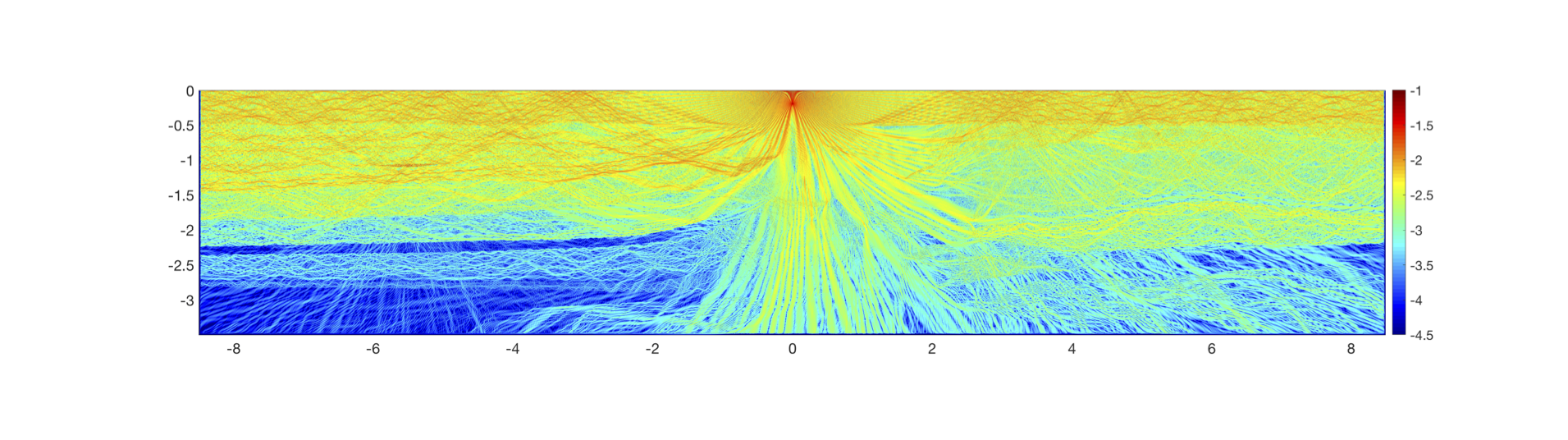}{\figWidth}};
\end{scope}
\end{scope}
\end{tikzpicture}
\end{center}
\caption{Displayed is the base 10 logarithm of the magnitude of the Helmholtz solution ($\log_{10} |u|$) caused by a point source near the surface. The results are, from top to botton, for $\omega = 800, 400  $ and 200. \label{fig:2d_marm} 
}
\end{figure}
}
As the number of unknowns is relatively large, $\sim 76 \cdot 10^6$, we parallelize the finite difference solver by a straightforward  domain decomposition with the communication handled by MPI. The simulations were carried out on Maneframe II at the Center for Scientific Computation at Southern Methodist University using 60 dual Intel Xeon E5-2695v4 2.1 GHz 18-core Broadwell processors with 45 MB of cache each and 256 GB of DDR4-2400 memory. The results displayed in Figure \ref{fig:2d_marm} and \ref{fig:2d_marm2} illustrate the ability of the method to find solutions to large problems and at high frequencies.     

{
\newcommand{\figWidth}{6.0cm}
\newcommand{\trimfig}[2]{\trimh{#1}{#2}{0.35}{0.33}{-0.16}{-0.13}}
\begin{figure}[]
\begin{center}
\begin{tikzpicture}[scale=1]
\useasboundingbox (0.,0.) rectangle (16.,4);  
\begin{scope}[xshift=-0.4cm,yshift=-0.9cm]
\begin{scope}[xshift=0.0cm]
\draw ( 0, 0.) node[anchor=south west,xshift= 0pt,yshift= 0pt] {\trimfig{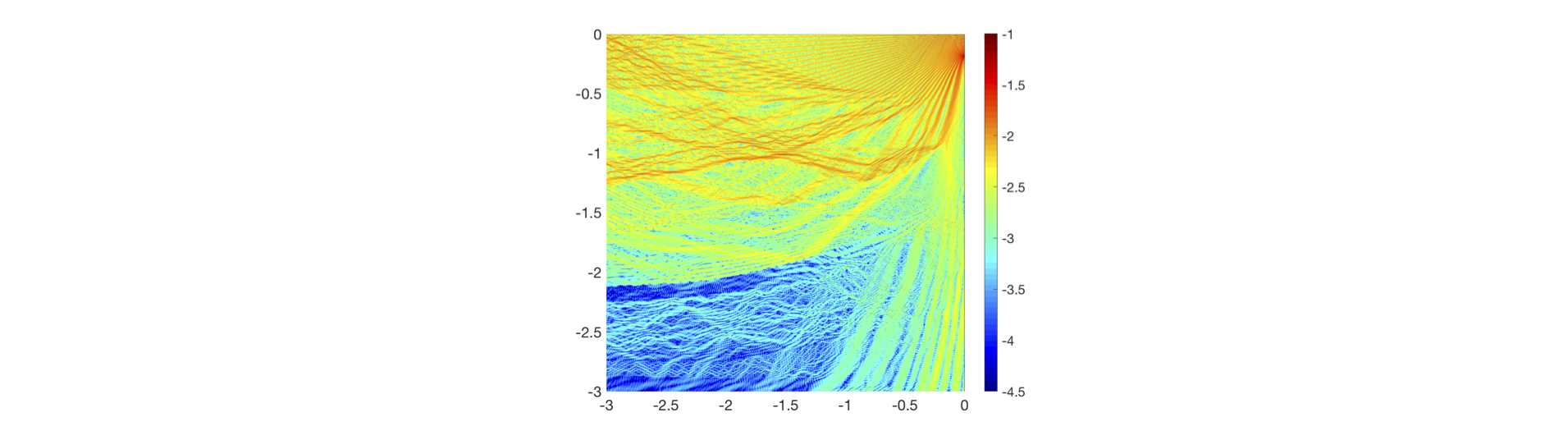}{\figWidth}};
\draw ( 5.2, 0.) node[anchor=south west,xshift=0pt,yshift=0pt] {\trimfig{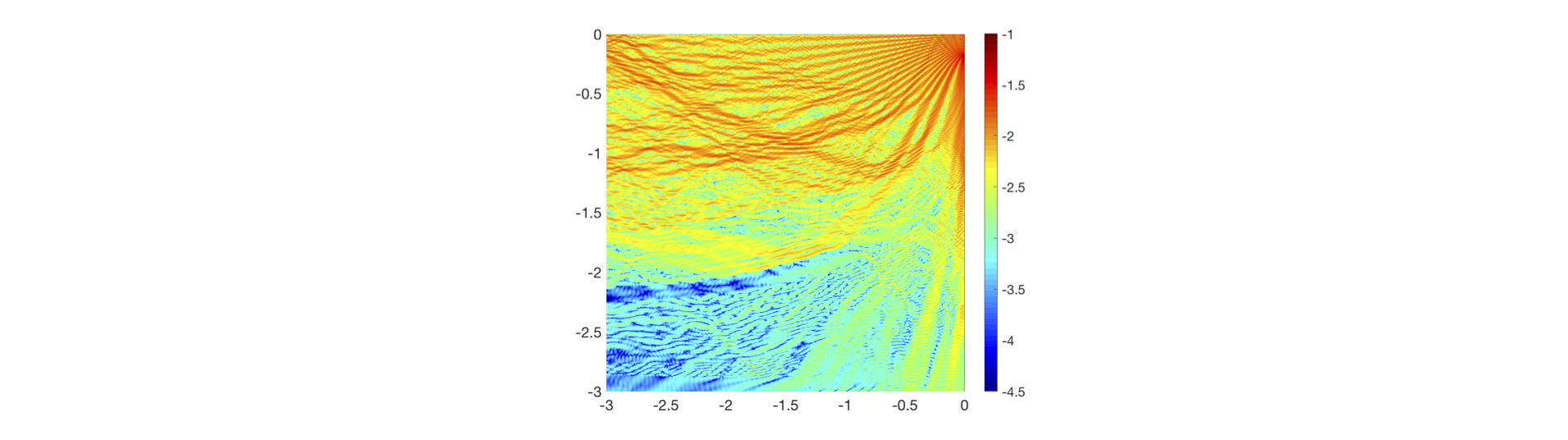}{\figWidth}};
\draw ( 10.4, 0.) node[anchor=south west,xshift=0pt,yshift=0pt] {\trimfig{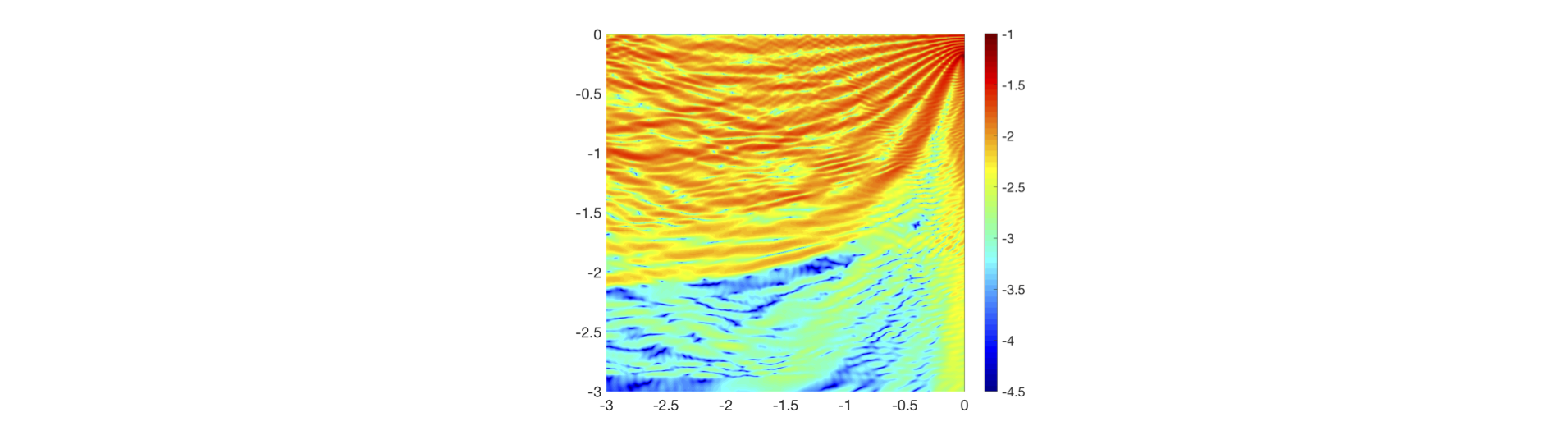}{\figWidth}};
\end{scope}
\end{scope}
\end{tikzpicture}
\end{center}
\caption{Zoom in of the base 10 logarithm of the magnitude of the Helmholtz solution ($\log_{10} |u|$) caused by a point source near the surface. The results are, from left to right, for $\omega = 800, 400  $ and 200.  \label{fig:2d_marm2} 
}
\end{figure}
}

\subsubsection{Multiple Frequencies}
{
\newcommand{\figWidth}{6.0cm}
\newcommand{\trimfig}[2]{\trimh{#1}{#2}{0.2}{0.17}{-0.15}{-0.15}}
\begin{figure}[]
\begin{center}
\begin{tikzpicture}[scale=1]
\useasboundingbox (0.3,0.) rectangle (15.,3.6);  
\begin{scope}[xshift=-0.8cm,yshift=-1.3cm]
\begin{scope}[xshift=0.0cm]
\draw ( 0, 0.) node[anchor=south west,xshift= 0pt,yshift= 0pt] {\trimfig{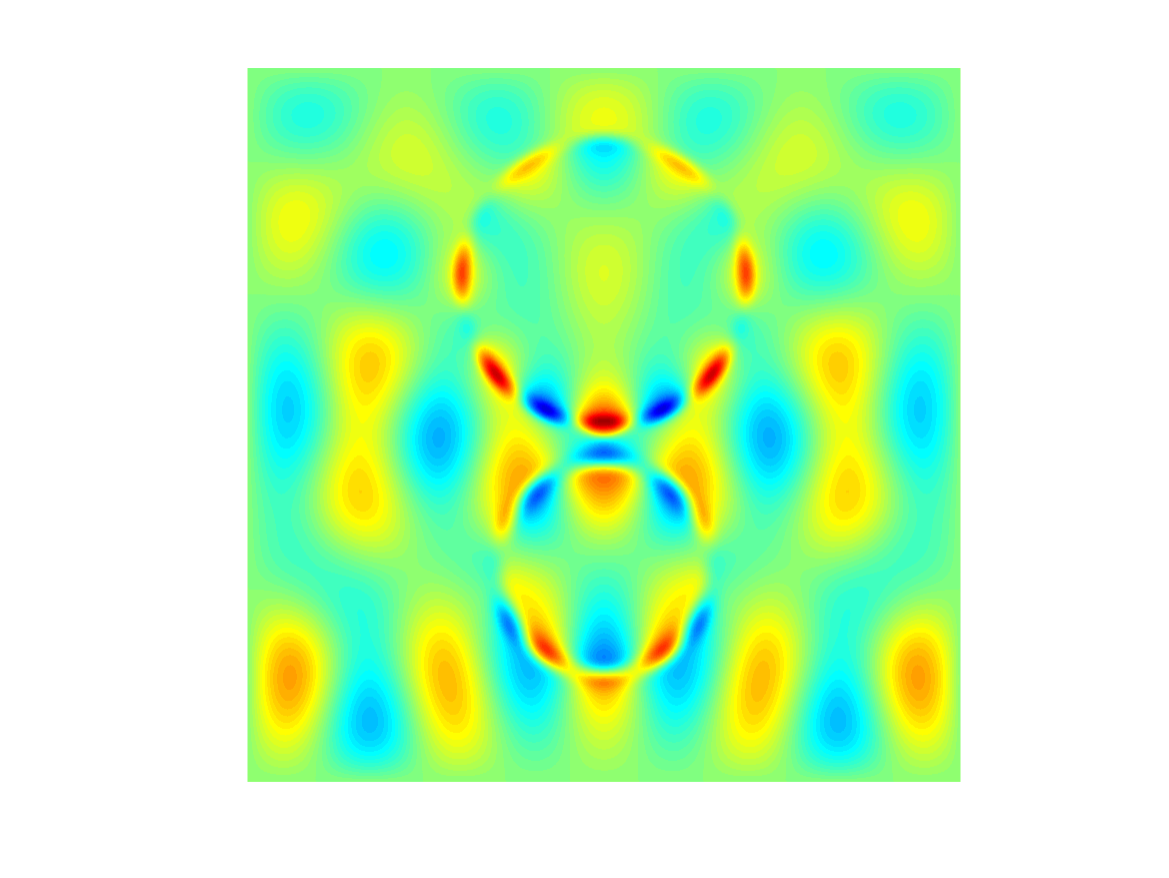}{\figWidth}};
\draw ( 4., 0.) node[anchor=south west,xshift=0pt,yshift=0pt] {\trimfig{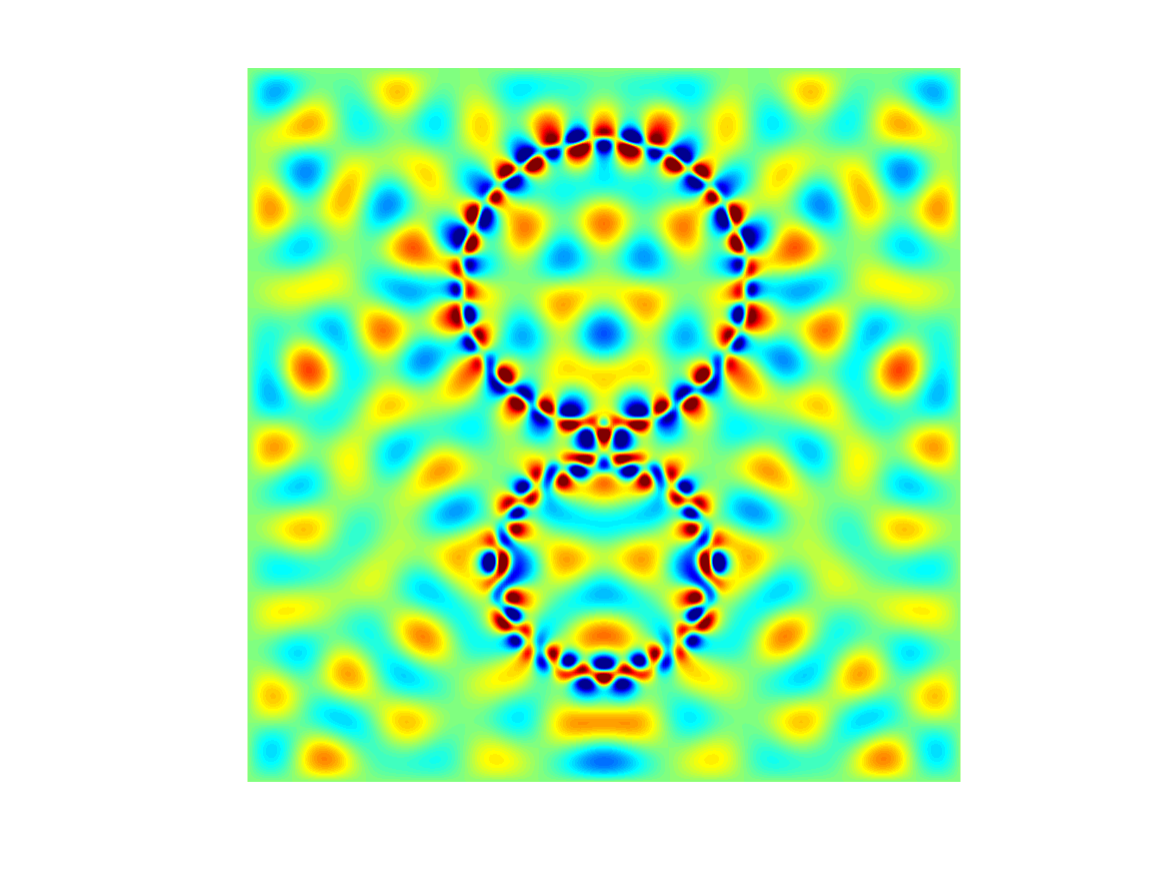}{\figWidth}};
\draw ( 8., 0.) node[anchor=south west,xshift=0pt,yshift=0pt] {\trimfig{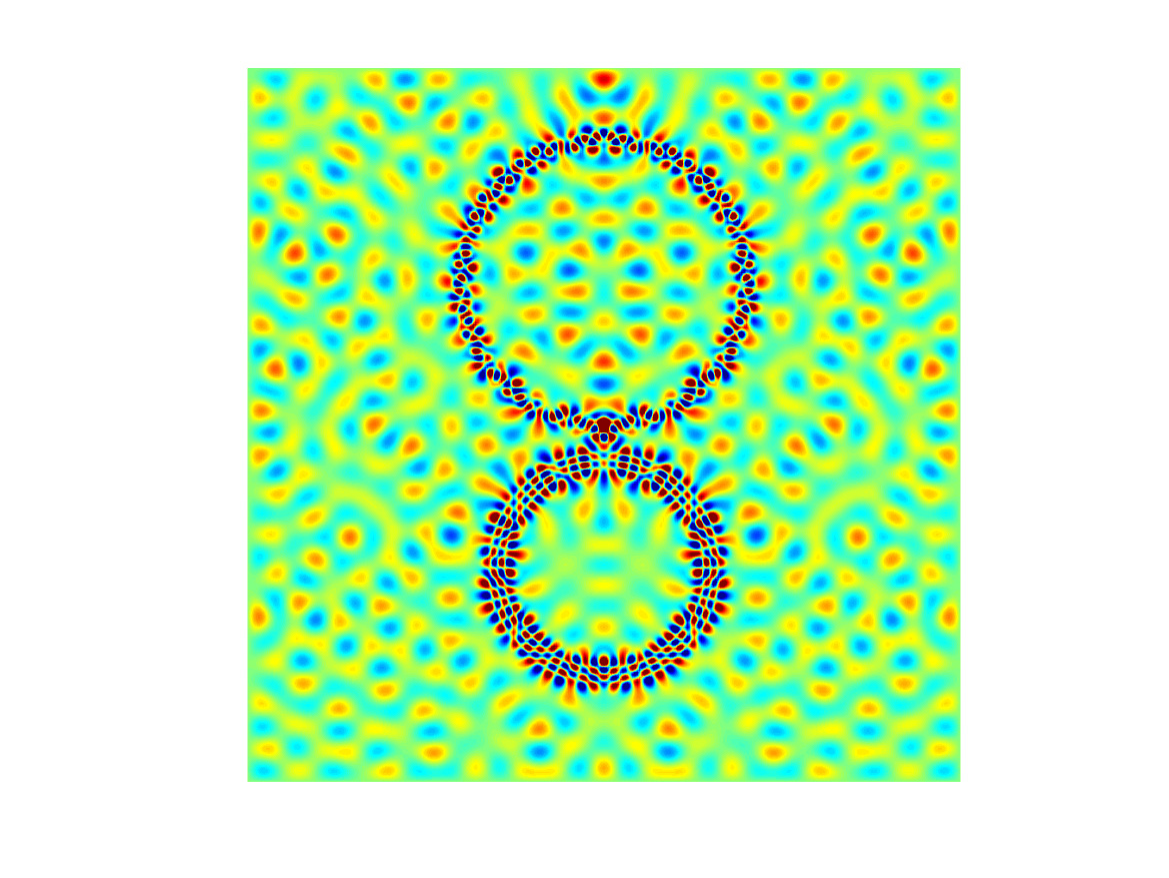}{\figWidth}};
\draw ( 12., 0.) node[anchor=south west,xshift=0pt,yshift=0pt] {\trimfig{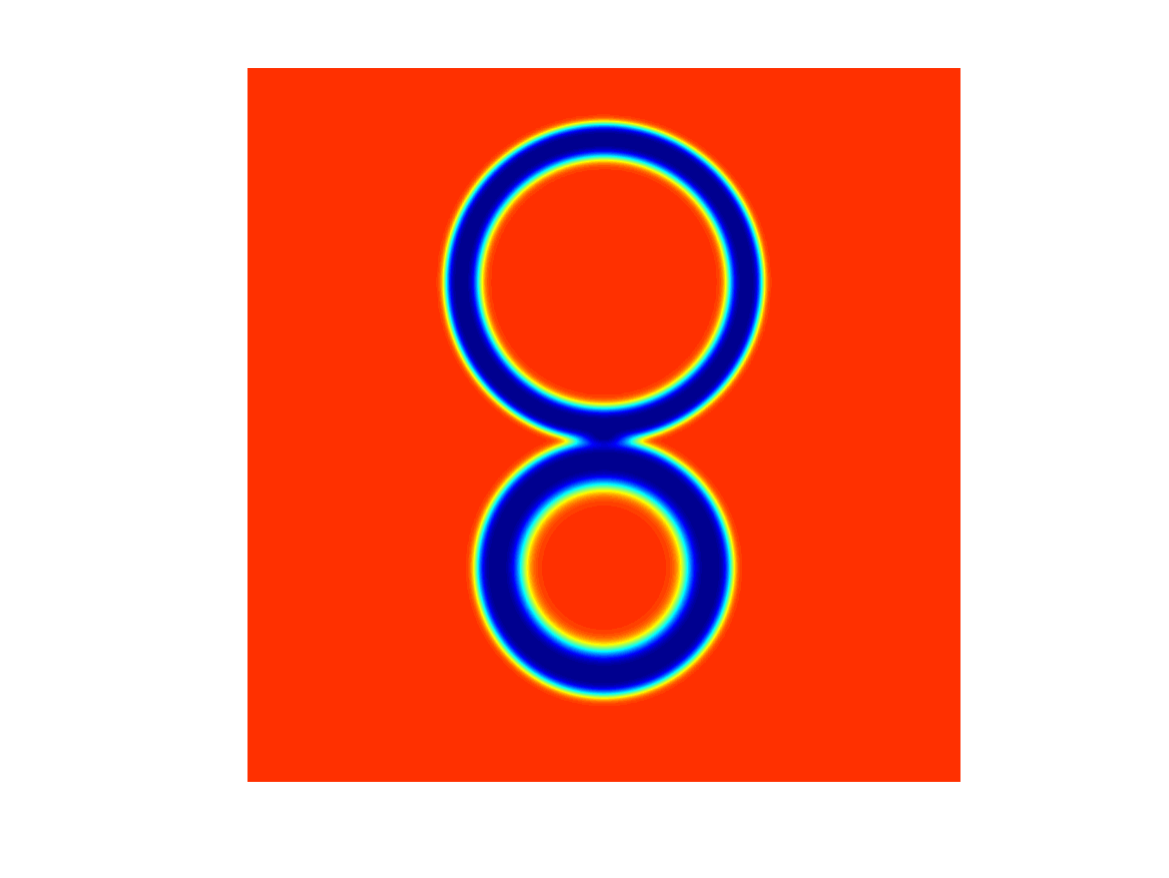}{\figWidth}};
\end{scope}
\end{scope}
\end{tikzpicture}
\end{center}
\caption{Computation of three Helmholtz problems by one solve. The frequencies are $\omega = 15, 30 $ and 60. The material model is also displayed, red is $c^2=1$ and dark blue is $c^2 = 0.1$. \label{fig:2d_variable_coeff} 
}
\end{figure}
}

In this final example in two dimensions we again use the sixth order accurate summation-by-parts discretization from \cite{Mattsson2012} with homogenous Dirichlet boundary conditions on the domain $(x,y) \in [-1,1]^2$. The spatial discretization size is the same in both coordinates and is taken to be 2/300. The velocity model is taken to be smoothly varying. Precisely we have that  
\begin{equation*}
     c^2(x,y)  = 1 - 0.9 \left(e^{-\left(\frac{(x^2+(y-0.4)^2-0.4^2)}{0.2^2}\right)^4}
     +e^{-\left(\frac{(x^2+(y+0.4)^2-0.3^2)}{0.2^2}\right)^4} \right),
\end{equation*}
see also Figure \ref{fig:2d_variable_coeff}. 
We consider three frequencies, $\omega=15,\,30,\,60$,
and use the same forcing in Helmholtz for all frequencies,
\begin{equation*}
f(x,y) = \frac{\sigma}{\pi}e^{-\sigma (x^2+y^2)}, \ \ \sigma = (4\omega)^2.
\end{equation*}
%
Here we use the WaveHoltz iteration over three periods
of the lowest frequency, 
accelerated by GMRES (with tolerance $10^{-8}$). We time step using a centered second order approximation to $w_{tt}$ with a timestep $\Delta t = 1/600$. Since we solve for three frequencies at once we adjust the filter as described 
in Section~\ref{sec:time_dep_fil}
 and extract all three solutions at once. Those solutions along with the material model are displayed in Figure \ref{fig:2d_variable_coeff}.

\begin{figure}[]
\graphicspath{{figures/}}
\begin{center}
\includegraphics[width=0.33\textwidth]{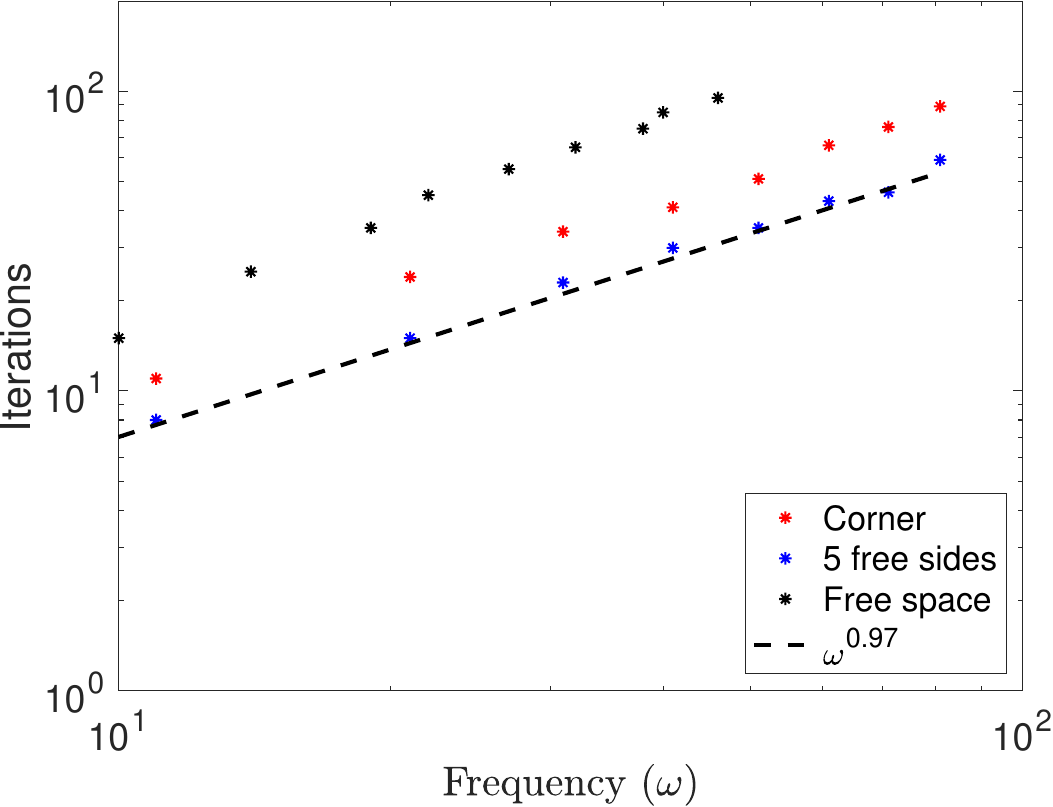}
\includegraphics[width=0.33\textwidth]{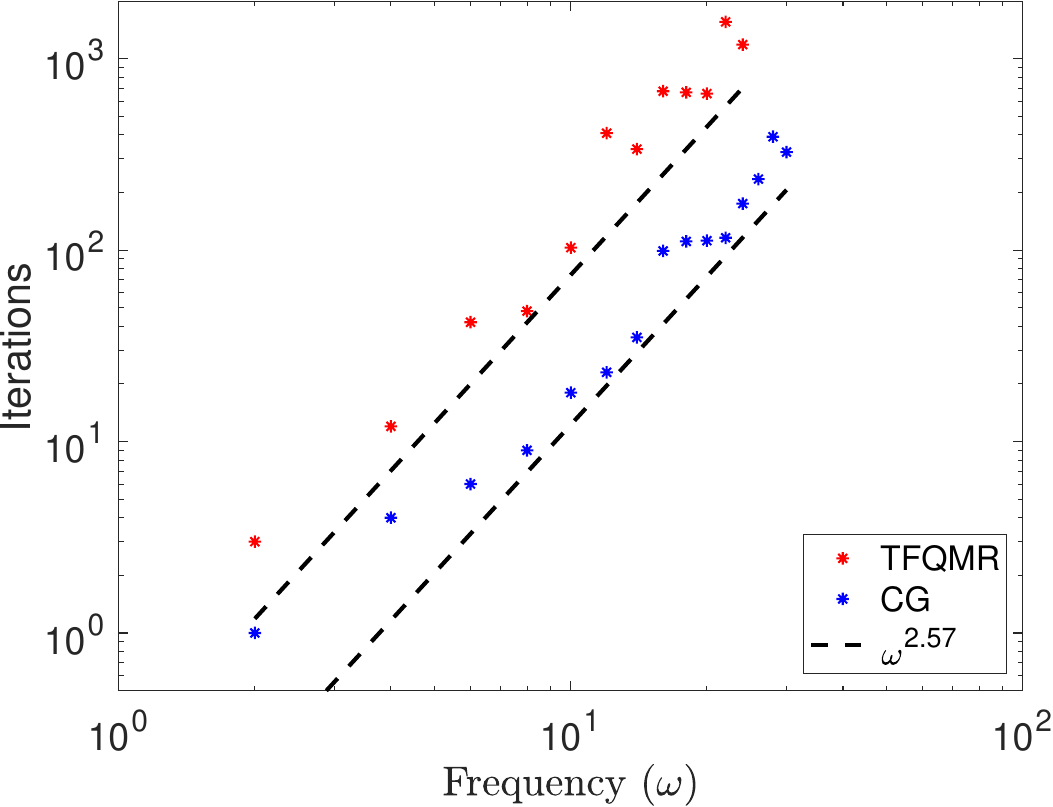}
\caption{To the left: number of iterations as a function of frequency to reduce the relative residual below $5\cdot 10^{-5}$ for problems with no trapped waves. Here WHI is accelerated by TFQMR. To the right: the same but for the interior problem. Here WHI is accelerated with either CG or TFQMR. \label{fig:3D_box_FD}
}
\end{center}
\end{figure}

\subsection{Problems in Three Dimensions}
In this section we present experiments in three dimensions. 
\subsubsection{Convergence in Different Geometries}
We solve the wave equation in a box $(x,y,z) \in [-1,1]^3$ with 
the smoothly varying medium 
\[
	c^2(x,y,z) = 1 + \frac{1}{10} e^{-\left( x^2+y^2+z^2 \right)}.
\]
We use a uniform grid $(x_i,y_j,z_k) = (-1+ih,-1+jh,-1+kh)$ with grid spacing $h=1/n$ and choose $n = \max (\lceil 10 \,\omega \rceil ,20)$ to keep the resolution fixed. The 
Helmholtz problem is forced by
\begin{equation}
      F(x,y,z) = \omega^3 e^{-36 \omega^2 \left( (x-x_0)^2 + (y-y_0)^2+(z-z_0)^2\right)} \label{eq:forcing},
\end{equation}
where $x_0 = 1/100$, $y_0 = 3/250$, and $z_0 = 1/200$.
We impose a mixture of boundary conditions consisting of homogenous Dirichlet and/or impedance boundary conditions: (1) impedance on all sides, (2) Dirichlet at $z = 1$ and impedance on all other sides, (3) Dirichlet at $z = -1$, $y = 1$, and $x = 1$ with impedance on all other sides, and (4) Dirichlet on all sides. We solve the equations in first order form in time and use the semi-discrete approximation described in Section \ref{sec:FD}.

In this example the WaveHoltz iteration is performed over $5$ periods, numerically integrated in time with the classic Runge Kutta method of order four, and accelerated by TFQMR method. For the pure Dirichlet problem we also use CG but note that although the spatial discretization leads to a symmetric WHI matrix when combined with a centered finite difference approximation in time the matrix is only close to symmetric when combined with the slightly dissipative Runge Kutta method. The experiments indicate that this slight non-symmetry does not destroy the convergence iteration of  CG in this case.

In Figure \ref{fig:3D_box_FD} we report the number of iterations needed to reduce the relative residual below $5\cdot 10^{-5}$. As was seen before in the $2$D case, the fully Dirichlet case is notably more difficult and requires more iterations than the other problems considered. The computational results indicate that the number of iterations to reach the tolerance scale as $\omega^{2.57}$ for the inner Dirichlet problem with either CG or TFQMR, with the former taking fewer overall iterations than the latter. By comparison, the set of problems with boundary conditions (1)-(3) listed above appear to converge in a number of iterations that scales as $\omega^{0.97}$, i.e. close to linear in the frequency $\omega$.
\begin{figure}[]
\graphicspath{{figures/}}
\begin{center}
\includegraphics[width=0.40\textwidth]{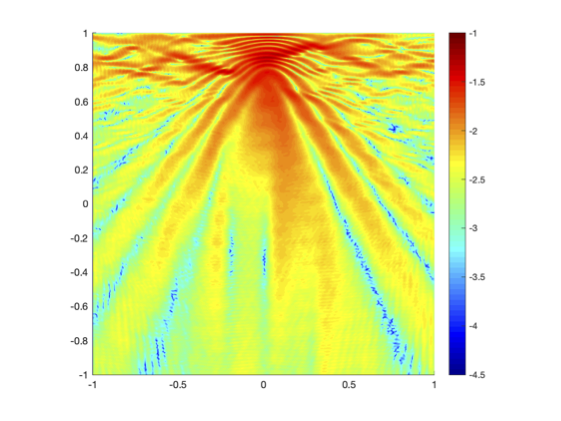}
\includegraphics[width=0.40\textwidth]{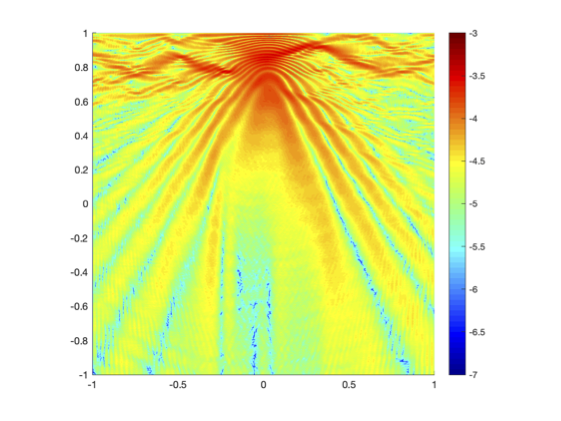}
\caption{Displayed is the base 10 logarithm of the magnitude of the Helmholtz solution ($\log_{10} |u|$) caused by a point source near the surface for $\omega = 200$ (left) and $\omega = 300$ (right) at the slice $x = 0.1$. \label{fig:3D_box}}
\end{center}
\end{figure}

\subsubsection{Scattering from a Plate}
For our final example, we again consider the box $(x,y,z) \in [-1,1]^3$ with smoothly varying medium 
\[
	c^2(x,y,z) = \frac{1}{2}\left[ 3 + \sin(16\pi z)\sin(4\pi(x + y))\right],
\]
and impose Dirichlet boundary conditions at $z = 1$ and impedance boundary conditions on all other sides. We use a uniform grid $(x_i,y_j,z_k) = (-1+ih,-1+jh,-1+kh)$ with grid spacing $h=1/n$ where $n$ is the number of gridpoints along a single dimension. The discretization in space and in time is exactly as in the previous example and the problem is forced by $ F(x,y,z)$ as in \eqref{eq:forcing} with $x_0 = 1/100$, $y_0 = 3/250$, and $z_0 = 4/5$. As in the Marmousi example in the previous section, we parallelize the finite difference solver by a straightforward domain decomposition with the communication handled by MPI. This simulation was carried out on Maneframe II at the Center for Scientific Computation at Southern Methodist University using 64 dual Intel Xeon E5-2695v4 2.1 GHz 18-core Broadwell processors with 45 MB of cache each and 256 GB of DDR4-2400 memory. The magnitude of the solution with $\omega = 200$ and $\omega = 300$ is plotted in Figure \ref{fig:3D_box}. We use $n = 1000$ for a total of $10^9$ gridpoints in the first case, and $n = 1500$ for a total of $3.375 \cdot 10^9$ gridpoints in the second for roughly $15\text{-}16$ points per wavelength.

\section{Summary and Future Work}
We have presented and analyzed the WaveHoltz iteration, a new iterative method for solving the Helmholtz equation. The iteration results in positive definite and sometimes symmetric matrices that are more amenable for iterative solution by Krylov subspace methods. \reva{In choosing a Krylov subspace method we note that CG is the most efficient and memory lean choice when the resulting system is symmetric positive definite, otherwise GMRES generally outperforms other methods such as QMR, LSQR, and TFQMR.}{} As the iteration is based on solving the wave equation it naturally parallelizes and can exploit techniques and spatial discretizations that have been developed for the time dependent problem. Numerical experiments indicate that our iteration appears to converge significantly faster than when the Helmholtz equation is discretized directly and solved iteratively with GMRES.

We believe that the numerical and theoretical results above are promising and note that there are many possible avenues for future exploration. For example we have exclusively used unconditioned Krylov solvers here but the spectral properties of the operator $\mathcal{S}$ indicate that preconditioning should be possible. Further, we have not tried to exploit adaptivity in space or time or any type of sweeping ideas here and we have only briefly touched on the possibilities for more advanced filter design. We hope to study both the numerical and theoretical properties of these in the future. 

Finally, here we only analyzed the energy conserving problem. We have carried out a preliminary but (we believe) suboptimal analysis of the impedance problem. We hope to improve our analysis of this case in the near future.     

\section{Acknowledgement}
We would like to acknowledge the careful reading of the manuscript by the reviewers that led to several improvements and clarifications. In particular one of the reviews was extremely careful and thorough and we want to extend an extra thank you for this.    

\appendix 

\section{Proof of \lem{filterlemma}}
\begin{proof}
We show the results for the rescaled transfer function
$$
  \bar{\beta}(r) := \beta(r\omega)
 = \frac{2}{T}\int_0^T\left(\cos(\omega t)-\frac14\right)\cos(r\omega t) dt = 
  \frac{1}{\pi}\int_0^{2\pi}\left(\cos(t)-\frac14\right)\cos(rt) dt.
$$
By direct integration we get
\begin{align}\lbeq{sincexp}
\bar\beta(r)& =
 \frac{1}{\pi} \int_0^{2\pi}
  \frac12(\cos((r+1)t)+\cos((r-1)t))-\frac14\cos (rt) dt = \\
 & \frac1{2\pi}\left(
  \frac{\sin(2\pi(r+1))}{r+1}+
  \frac{\sin(2\pi(r-1))}{r-1}
-  \frac12\frac{\sin(2\pi r)}{r}\right) 
= {\rm sinc}(r+1) + {\rm sinc}(r-1) -\frac12{\rm sinc}(r),
\end{align}
where 
$$
{\rm sinc}(r) =\frac{\sin(2\pi r)}{2\pi r}.
$$
We use the fact that $\sin(x)\leq x-\tilde\alpha x^3$ in the interval
$x\in[0,\pi]$ for any $\tilde\alpha\in [0,\pi^{-2}]$. This leads to
the following estimate for the sinc function
\begin{equation}\lbeq{sincest}
   0\leq  {\rm sinc}(r) \leq 1-\alpha r^2,\qquad r\in[-0.5,\ 0.5],\qquad \alpha\in[0,\ 4].
\end{equation}
We also note that ${\rm sinc}(r+n) = {\rm sinc}(r)r/(r+n)$ for all
integer $n$.

We now first consider $0\leq r\leq 0.5$ and use \eq{sincest} with $\alpha=4$ and $\alpha=0$,
$$
 |\bar\beta(r)|={\rm sinc}(r)\left|\frac{r}{r+1}+\frac{r}{r-1}-\frac{1}{2}
 \right|
 =\frac12 {\rm sinc}(r)\frac{1+3r^2}{1-r^2}
 \leq 
 \frac12 \frac{1-4r^2+3r^2}{1-r^2}=\frac12.
$$
For $0.5\leq r \leq 1.5$ we instead center around $r=1$ and get
for $|\delta|\leq 0.5$,
$$
\bar\beta(1+\delta)={\rm sinc}(\delta)
\frac{3(\delta+1)^2+1}{2(2+\delta)(1+\delta)}\geq 0,
$$
since sinc$(\delta)\geq 0$. Moreover,
using again \eq{sincest} with $\alpha=1$,
\begin{align*}
 \bar\beta(1+\delta)&
\leq
(1-\delta^2)\frac{3(\delta+1)^2+1}{2(2+\delta)(1+\delta)}
=
\frac{4+2\delta-3\delta^2-3\delta^3}{2(2+\delta)}
\leq 
\frac{4+2\delta-2\delta^2-\delta^3}{2(2+\delta)}
=1-\frac{\delta^2}{2}.
\end{align*}
Finally, for
$r>1$ we have $1/(r+1)-1/2r\geq 0$ and therefore
$$
  |\bar\beta(r)| =
  \frac{|\sin(2\pi r)|}{2\pi r}\left(\frac{r}{r+1} +
\frac{r}{r-1} -\frac1{2}\right)\leq
  \frac{1}{2\pi}\left(\frac{1}{r+1} +
\frac{1}{r-1} -\frac1{2r}\right)\leq
  \frac{3}{4\pi}\frac{1}{(r-1)}.
$$
We note also that for $r\geq 1.5$ this gives 
$|\bar\beta(r)| \leq 3/2\pi\leq 1/2$.

To prove \eq{locexp} we use the Taylor expansion
of
$\bar\beta$ around $r=1$,
$$
\bar\beta(1+\delta)= 1 +
\frac{\delta^2}{2}\bar\beta''(1)+\frac{\delta^3}{6} \bar{R}(\delta),
$$
where $\bar{R}(\delta)$ is the remainder term, which can be
bounded as
$$
  |\bar{R}(\delta)|\leq \sup_{r\geq 0} \left|\bar\beta^{(3)}(r)\right|\leq
  \frac{1}{\pi}\int_0^{2\pi}t^3\left(1+\frac14\right)dt
  =\reva{5\pi^{3}}{\frac{20\,\pi^{3}}{3}}.
$$
Hence, $|R(\delta)|\leq |\bar{R}(\delta)|/6\leq \reva{5\pi^3/6}{10\pi^3/9}$.
Finally, 
$$
\beta''(1) = {\rm sinc}''(2)+
{\rm sinc}''(0)
-\frac12{\rm sinc}''(1)=\frac{-1}{2}-\frac{(2\pi)^2}{3}
+1
=-2b_1.
$$
This shows \eq{locexp} 
and concludes the proof of the lemma.
\end{proof}

\section{Verification of Discrete Solution}\label{sec:verification}
Here we verify that \eq{discretesol} is indeed a solution to
the difference equation \eq{diffeq}. Direct substitution yields
\begin{align*}
\hat{w}_j^{n+1} - 2\hat{w}_j^n+\hat{w}_j^{n-1} + \Delta t^2\lambda_j^2 \hat{w}_j^n &=
   (\hat{v}_j-\hat{v}^\infty_j)\cos(\tilde\lambda_j t_n)
   \left(\cos(\tilde\lambda_j\Delta t)-2+\Delta t^2\lambda_j^2
   +\cos(\tilde\lambda_j\Delta t)
   \right)\\
   &\ \ 
    + \hat{v}^\infty_j\cos(\omega t_n)
    \left(\cos(\omega\Delta t)-2+\Delta t^2\lambda_j^2
   +\cos(\omega\Delta t)
   \right)\\
   &=
   (\hat{v}_j-\hat{v}^\infty_j)\cos(\tilde\lambda_j t_n)
   \left(-4\sin^2(\tilde\lambda_j\Delta t/2)+\Delta t^2\lambda_j^2
   )
   \right)\\
   &\ \ 
    + \hat{v}^\infty_j\cos(\omega t_n)
    \left(-4\sin^2(\omega\Delta t/2)+\Delta t^2\lambda_j^2
   \right)\\
   &=
    \hat{v}^\infty_j\cos(\omega t_n)
    \left(-\tilde\omega^2+\lambda_j^2
   \right)\Delta t^2 = -\Delta t^2f_j\cos(\omega t_n).
\end{align*}
Second, the initial conditions are satisfied since
\begin{align*}
  \hat{w}^0_j &= 
     (\hat{v}_j-\hat{v}^\infty_j) + \hat{v}^\infty_j=\hat{v}_j,\\
  \hat{w}^{-1}_j &= 
      (\hat{v}_j-\hat{v}^\infty_j)\cos(\tilde\lambda_j \Delta t) + \hat{v}^\infty_j\cos(\omega \Delta t) 
=    (\hat{v}_j-\hat{v}^\infty_j)\left(1-\frac12\Delta t^2\lambda_j^2\right)
+ \hat{v}^\infty_j\left(1-\frac12\Delta t^2\tilde\omega^2\right)\\
&=  \hat{v}_j\left(1-\frac12\Delta t^2\lambda_j^2\right)
+ \frac12\Delta t^2\hat{v}^\infty_j\left(\lambda_j^2-\tilde\omega^2\right)=  \hat{v}_j\left(1-\frac12\Delta t^2\lambda_j^2\right)
- \frac12\Delta t^2\hat{f}_j.
\end{align*}
This shows that \eq{discretesol} solves \eq{diffeq}.

\section{Proof of \lem{filterlemma2}}\label{sec:filterlemma2}

\begin{proof}

In general, we introduce the trapezoidal
rule applied to $\cos(\alpha t)$ in $[0,1]$,
$$
   \mathcal{T}_h(\alpha) := h \sum_{n=0}^M \eta_n \cos(\alpha t_n)
   \approx \int_0^1\cos(\alpha t) dt = \frac{\sin(\alpha)}{\alpha},
   \qquad h = 1/M,
$$
from which we attain the following lemma:
\begin{lemma}\lblem{traplemma}
The error in $\mathcal{T}_h(\alpha)$ satisfies\footnote{Note that this estimate is sharper than the standard error
estimate for the trapezoidal rule, which would have the factor
$\alpha^2$ from the second derivative of the integrand, not just $\alpha$.}
$$
  \left| \int_0^1\cos(\alpha t)dt - \mathcal{T}_h(\alpha)\right|
  \leq \frac{h^2|\alpha|}{\reva{\pi^2}{12}},\qquad \text{when $|h\alpha|\leq \pi$}.
$$
\end{lemma}
\begin{proof}
A direct calculation shows that
$$
\mathcal{T}_h(\alpha) = g(h\alpha)\int_0^1\cos(\alpha t)dt,
\qquad  g(x) = \frac{x}{2\tan(x/2)}. 
$$
The function $g(x)$ can be bounded as $\reva{1- x^2/\pi^2}{0}\leq g(x)\reva{\leq 1}{\leq 1-x^2/12}$ for $|x|\leq \pi$. This gives
\begin{align*}
\left| \int_0^1\cos(\alpha t)dt - \mathcal{T}_h(\alpha)\right|
=|1-g(h\alpha)|\left|\frac{\sin\alpha}{\alpha}\right|\leq
\frac{(h\alpha)^2}{\reva{\pi^2}{12}}\left|\frac{\sin\alpha}{\alpha}\right|
\leq\frac{h^2|\alpha|}{\reva{\pi^2}{12}}.
\end{align*}
\end{proof}

Since
$$
\left(\cos(\omega t)-\frac14\right)\cos(\lambda t)
  = \frac{1}{2}\left(
  \cos((\omega+\lambda) t)
  +\cos((\omega-\lambda) t)-\frac12\cos(\lambda t)\right),
$$
\reva{and $h=\Delta t/T$, we can write}{we can then let $h=\Delta t/T$ and write}
\begin{align*}
  \beta_h(\lambda) 
  &= 
     \frac{\Delta t}{T}
   \sum_{n=0}^M \reva{\eta_n}{\frac{\eta_n}{2}} \left(
\cos((\omega+\lambda) t_n)
  +\cos((\omega-\lambda) t_n)-\frac12\cos(\lambda t_n)
\right)
\\
  &= \reva{}{\frac{1}{2}} \left[\mathcal{T}_h(T(\omega+\lambda))+\mathcal{T}_h(T(\omega-\lambda))
  -\frac12 \mathcal{T}_h(T\lambda)\right].
\end{align*}
From \lem{traplemma}  we then get that
$$
|\beta(\tilde\lambda_j)-\beta_h(\tilde\lambda_j)|\leq 
\frac{h^2}{\reva{\pi^2}{12}}\left(\reva{}{\frac{1}{2}}|T(\omega+\tilde\lambda_j)|+\reva{}{\frac{1}{2}}|T(\omega-\tilde\lambda_j)|+\reva{\frac12}{\frac{1}{4}}|T\tilde\lambda_j|\right)
\leq \frac{5h^2T}{\reva{2\pi^2}{48}}(\omega+\tilde\lambda_j)
= \frac{5\Delta t^2}{\reva{2\pi^2}{48}T}(\omega+\tilde\lambda_j),
$$
when 
$$
\pi \geq
hT(\omega+\tilde\lambda_j) = \Delta t(\omega+\tilde\lambda_j),
$$
which is true by \eq{cfl} and the fact that arcsin$(x)\leq \pi x/2$
for $x\in[0,1]$:
\be\lbeq{tildelambda}
  \tilde\lambda_j = \frac{2}{\Delta t}{\rm arcsin}\left(\frac{\Delta t\lambda_j}{2}\right)
  \leq \frac{\pi}{2}\lambda_j.
\ee
Next, we use the inequality $|x-\sin(x)|\leq x^3/6$ for $|x|\leq \pi/2$
to show that
\begin{equation}\lbeq{sincerror}
  \left|\frac{\sin(xh)}{h}-x\right|
=  \frac{1}{h}\left|\reva{\sin(xh)}{\frac{\sin(xh)}{h}}-xh\right|
\leq \frac{(xh)^3}{6h} =\frac{h^2x^3}{6},\qquad |hx|\leq \frac{\pi}{2}.
\end{equation}
It gives us an estimate for $\tilde\lambda_j-\lambda_j$,
$$
|\lambda_j-\tilde\lambda_j|
=\left|  \frac{\sin\left(\Delta t\tilde{\lambda}_j/2\right)}{\Delta t/2}-\tilde\lambda_j\right|
  \leq 
  \frac{\Delta t^2}{24}\tilde\lambda_j^3,
$$
which is valid for all $j$ since $\Delta t\tilde\lambda_j/2\leq
 \Delta t\pi \lambda_j/4\leq  \Delta t\pi \lambda_N/4\leq\pi/2$,
 by \eq{cfl} and \eq{tildelambda}.
 
By \lem{filterlemma} 
$$
   |\beta(\omega+r)|\leq \begin{cases}
   1- \frac{r^2}{2\omega^2}, & |r/\omega|\leq \frac12,\\
   \frac12, & |r/\omega|\geq\frac12.
   \end{cases}
$$
We now claim that the statement \eq{betahest} in the
lemma holds
for all $j$ if $\Delta t\omega \leq \min(\delta_h,1)$.
On the one hand, if $|\omega-\tilde\lambda_j|\geq \omega/2$,
by \eq{cfl} and \eq{tildelambda}
$$
  |\beta_h(\tilde\lambda_j)|\leq 
  |\beta(\tilde\lambda_j)|+
  \frac{5\Delta t^2}{\reva{2\pi^2}{48}T}(\omega+\tilde\lambda_j)
  \leq \frac12
  +  \frac{5\Delta t^2}{\reva{2\pi^2}{48}T}(\omega+\tilde\lambda_j)
 \reva{=}{\leq} \frac12+
 \frac{5\Delta t\omega}{\reva{4\pi^3}{96\pi}}\Delta t(\omega+\tilde\lambda_j)\leq
 \frac12+\frac{5\Delta t\omega}{\reva{4\pi^2}{96}}\leq \reva{0.63}{0.6}.
$$
%
On the other hand, if $|\omega-\tilde\lambda_j|< \omega/2$,
$$
\frac{|\omega-\tilde\lambda_{j}|}{\omega}
\geq
\frac{|\omega-\lambda_{j}|}{\omega}
-\frac{|\tilde\lambda_{j}-\lambda_{j}|}{\omega}\geq 
\delta_h-
\frac{\Delta t^2}{24\omega}\tilde\lambda_j^3
\geq
\delta_h-\frac{\Delta t^2}{24\omega} (\omega+|\tilde\lambda_j-\omega|)^3
\geq
\delta_h-\frac{(3/2)^3}{24} \Delta t^2\omega^2\geq \frac{55}{64}\delta_h,
$$
since $\min(\delta_h,1)^2\leq \delta_h$. Then
\begin{align*}
  |\beta_h(\tilde\lambda_j)| &\leq 
  |\beta(\tilde\lambda_j)|+
  \frac{5\Delta t^2}{48T}(\omega+\tilde\lambda_j)
  \leq
  1-\frac12 \left(\frac{|\omega-\tilde\lambda_{j}|}{\omega}\right)^2
+
  \frac{5\Delta t^2\omega}{96\pi}(\omega+\omega+\omega/2)\\  
  &\leq 
  1-\frac{55^2}{2\cdot 64^2}\delta_h^2
+
  \Delta t^2\omega^2\frac{25}{192\pi}
  \leq 
    1-\left(\frac{55^2}{2\cdot 64^2}
-\frac{25}{192\pi}\right)\delta_h^2
\leq
  1-0.3\delta_h^2.
\end{align*}
This proves the lemma.

\end{proof}

\bibliographystyle{siam}
\bibliography{appelo}

\begin{thebibliography}{10}

\bibitem{Appthesis}
{\sc D.~Appel\"{o}}, {\em Absorbing {L}ayers and {N}on-{R}eflecting {B}oundary
  {C}onditions for {W}ave {P}ropagation {Problems}}, PhD thesis, {R}oyal
  {I}nstitute of {T}echnology, October 2005.

\bibitem{AppCol08}
{\sc D.~Appel\"{o} and T.~Colonius}, {\em A high-order super-grid-scale
  absorbing layer and its application to linear hyperbolic systems}, Journal of
  Computational Physics, 228 (2009), pp.~4200--4217.

\bibitem{Upwind2}
{\sc D.~Appel\"o and T.~Hagstrom}, {\em A new discontinuous {G}alerkin
  formulation for wave equations in second order form}, {SIAM} Journal On
  Numerical Analysis, 53 (2015), pp.~2705--2726.

\bibitem{el_dg_dath}
{\sc D.~Appel{\"{o}} and T.~Hagstrom}, {\em An energy-based discontinuous
  {G}alerkin discretization of the elastic wave equation in second order form},
  Comput. Meth. Appl. Mech. Engrg., 338 (2018), pp.~362--391.

\bibitem{AppPet09}
{\sc D.~Appel\"{o} and N.~A. Petersson}, {\em A stable finite difference method
  for the elastic wave equation on complex geometries with free surfaces},
  Communications in Computational Physics, 5 (2009), pp.~84--107.

\bibitem{BePriGab16}
{\sc Hadrien B{\'e}riot, Albert Prinn, and Gw{\'e}na{\"e}l Gabard}, {\em
  Efficient implementation of high-order finite elements for helmholtz
  problems}, International Journal for Numerical Methods in Engineering, 106
  (2016), pp.~213--240.

\bibitem{bjorck2015numerical}
{\sc A.~Bj\"{o}rck}, {\em Numerical methods in matrix computations}, vol.~59,
  Springer, 2015.

\bibitem{brandt1997wave}
{\sc A.~Brandt and I.~Livshits}, {\em Wave-ray multigrid method for standing
  wave equations}, Electron. Trans. Numer. Anal, 6 (1997), p.~91.

\bibitem{bristeau1998controllability}
{\sc M.O. Bristeau, R.~Glowinski, and J.~P{\'e}riaux}, {\em Controllability
  methods for the computation of time-periodic solutions; application to
  scattering}, Journal of Computational Physics, 147 (1998), pp.~265--292.

\bibitem{Helmholtz_Chen_Xiang_2013}
{\sc Z.~Chen and X.~Xiang}, {\em A source transfer domain decomposition method
  for {H}elmholtz equations in unbounded domain}, SIAM Journal on Numerical
  Analysis, 51 (2013), pp.~2331--2356.

\bibitem{EngMaj77}
{\sc B.~Engquist and A.~Majda}, {\em Absorbing boundary conditions for the
  numerical simulation of waves}, Math. Comp., 31 (1977), p.~629.

\bibitem{engquist2011sweepingH}
{\sc B.~Engquist and L.~Ying}, {\em Sweeping preconditioner for the {H}elmholtz
  equation: hierarchical matrix representation}, Communications on pure and
  applied mathematics, 64 (2011), pp.~697--735.

\bibitem{engquist2011sweepingPML}
\leavevmode\vrule height 2pt depth -1.6pt width 23pt, {\em Sweeping
  preconditioner for the {H}elmholtz equation: moving perfectly matched
  layers}, Multiscale Modeling \& Simulation, 9 (2011), pp.~686--710.

\bibitem{no_low_rank_Engq}
{\sc B.~Engquist and H.~Zhao}, {\em Approximate separability of the {G}reen's
  function of the {H}elmholtz equation in the high frequency limit},
  Communications on Pure and Applied Mathematics, 71 (2018), pp.~2220--2274.

\bibitem{erlangga2008advances}
{\sc Y.A. Erlangga}, {\em Advances in iterative methods and preconditioners for
  the {H}elmholtz equation}, Archives of Computational Methods in Engineering,
  15 (2008), pp.~37--66.

\bibitem{Erlangga2006}
{\sc Y.~Erlangga, C.~Oosterlee, and C.~Vuik}, {\em A novel multigrid based
  preconditioner for heterogeneous {H}elmholtz problems}, SIAM Journal on
  Scientific Computing, 27 (2006), pp.~1471--1492.

\bibitem{ernst2012difficult}
{\sc O.G. Ernst and M.J. Gander}, {\em Why it is difficult to solve {H}elmholtz
  problems with classical iterative methods}, in Numerical analysis of
  multiscale problems, Springer, 2012, pp.~325--363.

\bibitem{GANDER2000261}
{\sc M.~Gander and F.~Nataf}, {\em {AILU} for {H}elmholtz problems: a new
  preconditioner based on an analytic factorization}, Comptes Rendus de
  l'Acad{\'e}mie des Sciences - Series I - Mathematics, 331 (2000),
  pp.~261--266.

\bibitem{Gander_Zhang_SIAM_REV}
{\sc M.~Gander and H.~Zhang}, {\em A class of iterative solvers for the
  {H}elmholtz equation: Factorizations, sweeping preconditioners, source
  transfer, single layer potentials, polarized traces, and optimized {S}chwarz
  methods}, SIAM Review, 61 (2019), pp.~3--76.

\bibitem{Gillman2015}
{\sc A.~Gillman, A.H. Barnett, and P-G. Martinsson}, {\em A spectrally accurate
  direct solution technique for frequency-domain scattering problems with
  variable media}, BIT Numerical Mathematics, 55 (2015), pp.~141--170.

\bibitem{GlowRoss06}
{\sc R.~Glowinski and T.~Rossi}, {\em A mixed formulation and exact
  controllability approach for the computation of the periodic solutions of the
  scalar wave equation. (i): Controllability problem formulation and related
  iterative solution}, Comptes Rendus Math., 343 (2006), pp.~493--498.

\bibitem{grote2019controllability}
{\sc M.J. Grote and J.H. Tang}, {\em On controllability methods for the
  {H}elmholtz equation}, Journal of Computational and Applied Mathematics, 358
  (2019), pp.~306--326.

\bibitem{2019arXiv190312522G}
{\sc M.~J. {Grote}, F.~{Nataf}, J.~H. {Tang}, and P.-H. {Tournier}}, {\em
  {Parallel Controllability Methods For the Helmholtz Equation}}, arXiv
  e-prints,  (2019), p.~arXiv:1903.12522.

\bibitem{Hag99}
{\sc T.~Hagstrom}, {\em Radiation boundary conditions for the numerical
  simulation of waves}, Acta Numerica, 8 (1999), pp.~47--106.

\bibitem{HEIKKOLA2007344}
{\sc E.~Heikkola, S.~M\"{o}nk\"{o}l\"{a}, A.~Pennanen, and T.~Rossi}, {\em
  Controllability method for acoustic scattering with spectral elements},
  Journal of Computational and Applied Mathematics, 204 (2007), pp.~344--355.

\bibitem{HEIKKOLA20071553}
\leavevmode\vrule height 2pt depth -1.6pt width 23pt, {\em Controllability
  method for the {H}elmholtz equation with higher-order discretizations},
  Journal of Computational Physics, 225 (2007), pp.~1553--1576.

\bibitem{hile1998hybrid}
{\sc C.~V. Hile and G.~A Kriegsmann}, {\em A hybrid numerical method for loaded
  highly resonant single mode cavities}, Journal of Computational Physics, 142
  (1998), pp.~506--520.

\bibitem{ketcheson2015absolute}
{\sc D.~Ketcheson, L.~L{\'o}czi, and T.~Kocsis}, {\em On the absolute stability
  regions corresponding to partial sums of the exponential function}, IMA
  Journal of Numerical Analysis, 35 (2015), pp.~1426--1455.

\bibitem{KreOli72}
{\sc H.-O. Kreiss and J.~Oliger}, {\em Comparison of accurate methods for the
  integration of hyperbolic equations}, Tellus, 24 (1972), pp.~199--215.

\bibitem{Ladyzhenskaya:57}
{\sc O.~A. Ladyzhenskaya}, {\em On the limiting-amplitude principle}, Uspekhi
  Mat. Nauk, 12 (1957), pp.~161--164.

\bibitem{Livshits:2006aa}
{\sc I.~Livshits and A.~Brandt}, {\em Accuracy properties of the {W}ave-ray
  multigrid algorithm for {H}elmholtz equations}, SIAM Journal on Scientific
  Computing, 28 (2006), pp.~1228--24.

\bibitem{Mattsson2012}
{\sc K.~Mattsson}, {\em Summation by parts operators for finite difference
  approximations of second-derivatives with variable coefficients}, Journal of
  Scientific Computing, 51 (2012), pp.~650--682.

\bibitem{morawetz1962limiting}
{\sc C.S. Morawetz}, {\em The limiting amplitude principle}, Communications on
  Pure and Applied Mathematics, 15 (1962), pp.~349--361.

\bibitem{stolk2013rapidly}
{\sc C.C. Stolk}, {\em A rapidly converging domain decomposition method for the
  {H}elmholtz equation}, Journal of Computational Physics, 241 (2013),
  pp.~240--252.

\bibitem{Vainberg:75}
{\sc B.~R. Vainberg}, {\em On short-wave asymptotic behaviour of solutions to
  steady-state problems and the asymptotic behaviour as $t\to\infty$ of
  solutions of time-dependent problems}, Uspekhi Mat. Nauk, 30 (1975),
  pp.~1--58.

\bibitem{vion2014double}
{\sc A.~Vion and C.~Geuzaine}, {\em Double sweep preconditioner for optimized
  {S}chwarz methods applied to the {H}elmholtz problem}, Journal of
  Computational Physics, 266 (2014), pp.~171--190.

\bibitem{2011_deHoop_Xia}
{\sc S.~Wang, M.V. de~Hoop, and J.~Xia}, {\em On 3d modeling of seismic wave
  propagation via a structured parallel multifrontal direct {H}elmholtz
  solver}, Geophysical Prospecting, 59 (2011), pp.~857--873.

\bibitem{weyl:11}
{\sc H.~Weyl.}, {\em \"{U}ber die asymptotische verteilung der eigenwerte},
  Nachr. Konigl. Ges. Wiss.,  (1911), pp.~110--117.

\bibitem{2018arXiv180108655Z}
{\sc L.~{Zepeda-N{\'u}{\~n}ez}, A.~{Scheuer}, R.~J. {Hewett}, and
  L.~{Demanet}}, {\em {The method of polarized traces for the 3D {H}elmholtz
  equation}}, ArXiv e-prints,  (2018).

\bibitem{ZEPEDANUNEZ2016347}
{\sc L.~Zepeda-N{\'u}{\~n}ez and L.~Demanet}, {\em The method of polarized
  traces for the 2d {H}elmholtz equation}, Journal of Computational Physics,
  308 (2016), pp.~347--388.

\bibitem{leonardo_1}
\leavevmode\vrule height 2pt depth -1.6pt width 23pt, {\em Nested domain
  decomposition with polarized traces for the 2d {H}elmholtz equation}, SIAM
  Journal on Scientific Computing, 40 (2018), pp.~B942--B981.

\end{thebibliography}
\end{document}